\documentclass[11pt]{article}
\usepackage{amsmath, upgreek}
\usepackage{amssymb}
\usepackage{amscd}
\usepackage{xspace}
\usepackage{verbatim}
	
\usepackage{color}

\newcommand{\mysection}[1]{
\section{#1}\setcounter{equation}{0}}
\title{\bf Boundary singular solutions of a class of equations with mixed absorption-reaction}
\author{{\bf Marie-Fran\c{c}oise Bidaut-V\'eron\footnote{\noindent Laboratoire de Math\'{e}matiques et Physique Th\'{e}orique,
Universit\'e de Tours, 37200 Tours, France. E-mail: veronmf@univ-tours.fr},} \\{\bf Marta Garcia-Huidobro \footnote{\noindent
Departamento de Matematicas, Pontifica Universidad Catolica de Chile
Casilla 307, Correo 2, Santiago de Chile. E-mail: mgarcia@mat.puc.cl}}\\
 {\bf Laurent V\'eron \footnote{\noindent
Laboratoire de Math\'{e}matiques et Physique Th\'{e}orique, Universit\'e de Tours, 37200 Tours, France. E-mail: veronl@univ-tours.fr}}\\[2mm]
}

\date{}
\begin{document}
 \maketitle


\newcommand{\txt}[1]{\;\text{ #1 }\;}
\newcommand{\tbf}{\textbf}
\newcommand{\tit}{\textit}
\newcommand{\tsc}{\textsc}
\newcommand{\trm}{\textrm}
\newcommand{\mbf}{\mathbf}
\newcommand{\mrm}{\mathrm}
\newcommand{\bsym}{\boldsymbol}
\newcommand{\scs}{\scriptstyle}
\newcommand{\sss}{\scriptscriptstyle}
\newcommand{\txts}{\textstyle}
\newcommand{\dsps}{\displaystyle}
\newcommand{\fnz}{\footnotesize}
\newcommand{\scz}{\scriptsize}
\newcommand{\be}{\begin{equation}}
\newcommand{\bel}[1]{\begin{equation}\label{#1}}
\newcommand{\ee}{\end{equation}}
\newcommand{\eqnl}[2]{\begin{equation}\label{#1}{#2}\end{equation}}
\newcommand{\barr}{\begin{eqnarray}}
\newcommand{\earr}{\end{eqnarray}}
\newcommand{\bars}{\begin{eqnarray*}}
\newcommand{\ears}{\end{eqnarray*}}
\newcommand{\nnu}{\nonumber \\}
\newtheorem{subn}{\name}
\renewcommand{\thesubn}{}
\newcommand{\bsn}[1]{\def\name{#1}\begin{subn}}
\newcommand{\esn}{\end{subn}}
\newtheorem{sub}{\name}[section]
\newcommand{\dn}[1]{\def\name{#1}}   
\newcommand{\bs}{\begin{sub}}
\newcommand{\es}{\end{sub}}
\newcommand{\bsl}[1]{\begin{sub}\label{#1}}
\newcommand{\bth}[1]{\def\name{Theorem}
\begin{sub}\label{t:#1}}
\newcommand{\blemma}[1]{\def\name{Lemma}
\begin{sub}\label{l:#1}}
\newcommand{\bcor}[1]{\def\name{Corollary}
\begin{sub}\label{c:#1}}
\newcommand{\bdef}[1]{\def\name{Definition}
\begin{sub}\label{d:#1}}
\newcommand{\bprop}[1]{\def\name{Proposition}
\begin{sub}\label{p:#1}}
\newcommand{\R}{\eqref}
\newcommand{\rth}[1]{Theorem~\ref{t:#1}}
\newcommand{\rlemma}[1]{Lemma~\ref{l:#1}}
\newcommand{\rcor}[1]{Corollary~\ref{c:#1}}
\newcommand{\rdef}[1]{Definition~\ref{d:#1}}
\newcommand{\rprop}[1]{Proposition~\ref{p:#1}}
\newcommand{\BA}{\begin{array}}
\newcommand{\EA}{\end{array}}
\newcommand{\BAN}{\renewcommand{\arraystretch}{1.2}
\setlength{\arraycolsep}{2pt}\begin{array}}
\newcommand{\BAV}[2]{\renewcommand{\arraystretch}{#1}
\setlength{\arraycolsep}{#2}\begin{array}}
\newcommand{\BSA}{\begin{subarray}}
\newcommand{\ESA}{\end{subarray}}
\newcommand{\BAL}{\begin{aligned}}
\newcommand{\EAL}{\end{aligned}}
\newcommand{\BALG}{\begin{alignat}}
\newcommand{\EALG}{\end{alignat}}
\newcommand{\BALGN}{\begin{alignat*}}
\newcommand{\EALGN}{\end{alignat*}}
\newcommand{\note}[1]{\textit{#1.}\hspace{2mm}}
\newcommand{\Proof}{\note{Proof}}
\newcommand{\qeda}{\hspace{10mm}\hfill $\square$}
\newcommand{\qed}{\\
${}$ \hfill $\square$}
\newcommand{\Remark}{\note{Remark}}
\newcommand{\modin}{$\,$\\[-4mm] \indent}
\newcommand{\forevery}{\quad \forall}
\newcommand{\set}[1]{\{#1\}}
\newcommand{\setdef}[2]{\{\,#1:\,#2\,\}}
\newcommand{\setm}[2]{\{\,#1\mid #2\,\}}
\newcommand{\mt}{\mapsto}
\newcommand{\lra}{\longrightarrow}
\newcommand{\lla}{\longleftarrow}
\newcommand{\llra}{\longleftrightarrow}
\newcommand{\Lra}{\Longrightarrow}
\newcommand{\Lla}{\Longleftarrow}
\newcommand{\Llra}{\Longleftrightarrow}
\newcommand{\warrow}{\rightharpoonup}
\newcommand{
\paran}[1]{\left (#1 \right )}
\newcommand{\sqbr}[1]{\left [#1 \right ]}
\newcommand{\curlybr}[1]{\left \{#1 \right \}}
\newcommand{\abs}[1]{\left |#1\right |}
\newcommand{\norm}[1]{\left \|#1\right \|}
\newcommand{
\paranb}[1]{\big (#1 \big )}
\newcommand{\lsqbrb}[1]{\big [#1 \big ]}
\newcommand{\lcurlybrb}[1]{\big \{#1 \big \}}
\newcommand{\absb}[1]{\big |#1\big |}
\newcommand{\normb}[1]{\big \|#1\big \|}
\newcommand{
\paranB}[1]{\Big (#1 \Big )}
\newcommand{\absB}[1]{\Big |#1\Big |}
\newcommand{\normB}[1]{\Big \|#1\Big \|}
\newcommand{\produal}[1]{\langle #1 \rangle}

\newcommand{\thkl}{\rule[-.5mm]{.3mm}{3mm}}
\newcommand{\thknorm}[1]{\thkl #1 \thkl\,}
\newcommand{\trinorm}[1]{|\!|\!| #1 |\!|\!|\,}
\newcommand{\bang}[1]{\langle #1 \rangle}
\def\angb<#1>{\langle #1 \rangle}
\newcommand{\vstrut}[1]{\rule{0mm}{#1}}
\newcommand{\rec}[1]{\frac{1}{#1}}
\newcommand{\opname}[1]{\mbox{\rm #1}\,}
\newcommand{\supp}{\opname{supp}}
\newcommand{\dist}{\opname{dist}}
\newcommand{\myfrac}[2]{{\displaystyle \frac{#1}{#2} }}
\newcommand{\myint}[2]{{\displaystyle \int_{#1}^{#2}}}
\newcommand{\mysum}[2]{{\displaystyle \sum_{#1}^{#2}}}
\newcommand {\dint}{{\displaystyle \myint\!\!\myint}}
\newcommand{\q}{\quad}
\newcommand{\qq}{\qquad}
\newcommand{\hsp}[1]{\hspace{#1mm}}
\newcommand{\vsp}[1]{\vspace{#1mm}}
\newcommand{\ity}{\infty}
\newcommand{\prt}{\partial}
\newcommand{\sms}{\setminus}
\newcommand{\ems}{\emptyset}
\newcommand{\ti}{\times}
\newcommand{\pr}{^\prime}
\newcommand{\ppr}{^{\prime\prime}}
\newcommand{\tl}{\tilde}
\newcommand{\sbs}{\subset}
\newcommand{\sbeq}{\subseteq}
\newcommand{\nind}{\noindent}
\newcommand{\ind}{\indent}
\newcommand{\ovl}{\overline}
\newcommand{\unl}{\underline}
\newcommand{\nin}{\not\in}
\newcommand{\pfrac}[2]{\genfrac{(}{)}{}{}{#1}{#2}}

\def\ga{\alpha}     \def\gb{\beta}       \def\gg{\gamma}
\def\gc{\chi}       \def\gd{\delta}      \def\ge{\epsilon}
\def\gth{\theta}                         \def\vge{\varepsilon}
\def\gf{\phi}       \def\vgf{\varphi}    \def\gh{\eta}
\def\gi{\iota}      \def\gk{\kappa}      \def\gl{\lambda}
\def\gm{\mu}        \def\gn{\nu}         \def\gp{\pi}
\def\vgp{\varpi}    \def\gr{\rho}        \def\vgr{\varrho}
\def\gs{\sigma}     \def\vgs{\varsigma}  \def\gt{\tau}
\def\gu{\upsilon}   \def\gv{\vartheta}   \def\gw{\omega}
\def\gx{\xi}        \def\gy{\psi}        \def\gz{\zeta}
\def\Gg{\Gamma}     \def\Gd{\Delta}      \def\Gf{\Phi}
\def\Gth{\Theta}
\def\Gl{\Lambda}    \def\Gs{\Sigma}      \def\Gp{\Pi}
\def\Gw{\Omega}     \def\Gx{\Xi}         \def\Gy{\Psi}

\def\CS{{\mathcal S}}   \def\CM{{\mathcal M}}   \def\CN{{\mathcal N}}
\def\CR{{\mathcal R}}   \def\CO{{\mathcal O}}   \def\CP{{\mathcal P}}
\def\CA{{\mathcal A}}   \def\CB{{\mathcal B}}   \def\CC{{\mathcal C}}
\def\CD{{\mathcal D}}   \def\CE{{\mathcal E}}   \def\CF{{\mathcal F}}
\def\CG{{\mathcal G}}   \def\CH{{\mathcal H}}   \def\CI{{\mathcal I}}
\def\CJ{{\mathcal J}}   \def\CK{{\mathcal K}}   \def\CL{{\mathcal L}}
\def\CT{{\mathcal T}}   \def\CU{{\mathcal U}}   \def\CV{{\mathcal V}}
\def\CZ{{\mathcal Z}}   \def\CX{{\mathcal X}}   \def\CY{{\mathcal Y}}
\def\CW{{\mathcal W}} \def\CQ{{\mathcal Q}}
\def\BBA {\mathbb A}   \def\BBb {\mathbb B}    \def\BBC {\mathbb C}
\def\BBD {\mathbb D}   \def\BBE {\mathbb E}    \def\BBF {\mathbb F}
\def\BBG {\mathbb G}   \def\BBH {\mathbb H}    \def\BBI {\mathbb I}
\def\BBJ {\mathbb J}   \def\BBK {\mathbb K}    \def\BBL {\mathbb L}
\def\BBM {\mathbb M}   \def\BBN {\mathbb N}    \def\BBO {\mathbb O}
\def\BBP {\mathbb P}   \def\BBR {\mathbb R}    \def\BBS {\mathbb S}
\def\BBT {\mathbb T}   \def\BBU {\mathbb U}    \def\BBV {\mathbb V}
\def\BBW {\mathbb W}   \def\BBX {\mathbb X}    \def\BBY {\mathbb Y}
\def\BBZ {\mathbb Z}   \def\BBQ {\mathbb Q}

\def\GTA {\mathfrak A}   \def\GTB {\mathfrak B}    \def\GTC {\mathfrak C}
\def\GTD {\mathfrak D}   \def\GTE {\mathfrak E}    \def\GTF {\mathfrak F}
\def\GTG {\mathfrak G}   \def\GTH {\mathfrak H}    \def\GTI {\mathfrak I}
\def\GTJ {\mathfrak J}   \def\GTK {\mathfrak K}    \def\GTL {\mathfrak L}
\def\GTM {\mathfrak M}   \def\GTN {\mathfrak N}    \def\GTO {\mathfrak O}
\def\GTP {\mathfrak P}   \def\GTR {\mathfrak R}    \def\GTS {\mathfrak S}
\def\GTT {\mathfrak T}   \def\GTU {\mathfrak U}    \def\GTV {\mathfrak V}
\def\GTW {\mathfrak W}   \def\GTX {\mathfrak X}    \def\GTY {\mathfrak Y}
\def\GTZ {\mathfrak Z}   \def\GTQ {\mathfrak Q}

\font\Sym= msam10 
\def\SYM#1{\hbox{\Sym #1}}
\newcommand{\bdw}{\prt\Gw\xspace}
\maketitle\medskip

{\abstract We study properties of positive functions satisfying (E)$\;-\Gd u+u^p-M|\nabla u|^q=0$ is a domain $\Gw$ or in $\BBR^{_N}_+$ when $p>1$ and $1<q<\min\{p,2\}$. We concentrate our research on the solutions of (E) vanishing on the boundary except at one point. This analysis depends on the existence of separable solutions in $\BBR^{_N}_+$. We construct various types of positive solutions with an isolated singularity on the boundary. We also study conditions for the removability of compact boundary sets and the Dirichlet problem associated to (E) with a measure as boundary data.}\medskip

\nind {\it 2010 Mathematics Subject Classification:} 35J62-35J66-35J75-31C15\smallskip

\nind{\it  Keywords:} Elliptic equations, boundary singularities, Bessel capacities, measures, supersolutions, subsolutions.

\tableofcontents
\date{}



\mysection{Introduction}
The aim of this article is to study some properties of solutions of the following equation
\begin{equation}\label{Z1}
\BA{lll}
\CL_{q,M} u:=-\Gd u+|u|^{p-1}u-M|\nabla u|^{q}=0
\EA
\end{equation}
in a bounded domain $\Gw$ of $\BBR^{_N}$ or in the half-space $\BBR_+^N$, where $M>0$ and $p>q>1$. We are particularly interested in the analysis of boundary singularities of such solutions. If $M=0$ the boundary singularities problem has been investigated since thirty years, starting with the work of Gmira and V\'eron \cite{GmVe} who obtained an almost complete description of the solutions with isolated boundary singularities. When $M>0$  there is a balance between the absorption term $|u|^{p-1}u$ and the source term $M|\nabla u|^{q}$, a confrontation which can create very new effects. Furthermore, the scale of the two opposed reaction terms depends upon the position of $q$ with respect to $\frac{2p}{p+1}$. This is due to the fact that $(\ref{Z1})$ is equivariant with respect to the scaling transformation $T_\ell$ defined for $\ell>0$ by $T_\ell[u](x)=\ell^{\frac{2}{p-1}}u(\ell x)$.\smallskip

\nind If $q<\frac{2p}{p+1}$, the absorption term is dominant and the behaviour of the singular solutions is modelled by the equation studied in \cite{GmVe}
\begin{equation}\label{Z3}
\BA{lll}
-\Gd u+|u|^{p-1}u=0.
\EA
\end{equation}
If $q>\frac{2p}{p+1}$ , the source term is dominant and the behaviour of the singular solutions is modelled by positive separable solutions of the equation without diffusion
 \bel {SS15e}
u^p-M|\nabla u|^q=0.
\ee
Another associated equation which plays an important role in the construction of singular solutions since its positive solutions are supersolutions of $(\ref{Z1})$ is
\begin{equation}\label{Z4}
\BA{lll}
-\Gd u-M|\nabla u|^{q}=0.
\EA
\end{equation}
Note that in $(\ref{SS15e})$ and $(\ref{Z4})$, $M$ can be fixed to be $1$ by replacing $u$ by $\ell u$.\\
If $q=\frac{2p}{p+1}$, the coefficient $M>0$ plays a fundamental role in the properties of the set of solutions, in particular for the existence of singular solutions and removable singularities. This situation is similar in some sense to what happens for equation 
\begin{equation}\label{Z5}
\BA{lll}
-\Gd u=|u|^{p-1}u+M|\nabla u|^{q}
\EA
\end{equation}
which is studied thoroughfly in \cite {BVGHV1}, \cite{BVGHV2} in the case $M>0$ and in \cite {SeZ} in the case $M<0$. In this last article the opposition of a forcing term $|u|^{p-1}u$ and an absorption term $M|\nabla u|^{q}$ creates a very rich configuration of unexpected phenomena and new effects.  \medskip\medskip

In the present paper we will consider the case where $1<q<2$, with a special emphasis 
on the case $q=\frac{2p}{p+1}$ which allows to put into light the role of the value of $M$. 
We first analyze the following problem: given a smooth bounded domain $\Gw\subset\BBR^{_N}$ such that $0\in\prt\Gw$, under what conditions involving $p$, $q$ and  $M$  is the point $0$ a removable singularity for a solution of $(\ref{Z1})$ continuous in $\overline\Gw\setminus \{0\}$ and vanishing on $\prt\Gw\setminus \{0\} $ ? In the sequel we denote $\gr(x)=\dist(x,\prt\Gw)$ and for $1\leq s<\infty$, $L^s_\gr(\Gw):=L^s(\Gw;\gr dx)$ and the space of test functions in $\Gw$ is defined by
\begin{equation}\label{NZ4a}
\BA{lll}
\BBX(\Gw)=\left\{\gz\in C^1(\overline\Gw):\gz=0\text{ on }\prt\Gw,\,\Gd\gz\in L^\infty(\Gw)\right\}.
\EA
\end{equation}
If $\Gw$ is replaced by $\BBR^{_N}_+$, then 
\begin{equation}\label{NZ5}
\BA{lll}
\BBX(\BBR^{_N}_+)=\left\{\gz\in C^1(\overline\BBR^{_N}_+)\text{ with compact support in }\overline\BBR^{_N}_+,\,\Gd\gz\in L^\infty(\BBR^{_N}_+)\right\}.
\EA
\end{equation}
Our first result is the following:
\bth{remov} Assume $p\geq\frac{N+1}{N-1}$, $M>0$ and \smallskip

\noindent (i) either $p=\frac{N+1}{N-1}$ and $1<q< 1+\frac{1}{N}$.\smallskip

\noindent (ii) or  $p>\frac{N+1}{N-1}$ and $1<q\leq \frac{2p}{p+1}$.\smallskip

\noindent Then any nonnegative solution $u\in C^2(\Gw)\cap C^1(\overline\Gw\setminus \{0\})$ of 
\begin{equation}\label{NZ1}
\BA{lll}
-\Gd u+|u|^{p-1}u-M|\nabla u|^{q}=0\qquad&\text{in }\Gw\\
\phantom{-\Gd +|u|^{p-1}u-M|\nabla u|^{q}}u=0&\text{in }\prt\Gw\setminus \{0\}
\EA
\end{equation}
verifies $\nabla u\in L_\gr^q(\Gw)$, $u\in L_\gr^p\Gw)$ and is a weak solution of 
\begin{equation}\label{NZ2}
\BA{lll}
-\Gd u+|u|^{p-1}u-M|\nabla u|^{q}=0\qquad&\text{in }\Gw\\
\phantom{-\Gd +|u|^{p-1}u-M|\nabla u|^{q}}u=0&\text{in }\prt\Gw,
\EA
\end{equation}
in the sense that 
\begin{equation}\label{NZ3}
\BA{lll}
\myint{\Gw}{}\left(-u\Gd\gz+(|u|^{p-1}u-M|\nabla u|^{q})\gz\right)dx=0\quad\text{for all }\gz\in \BBX(\Gw).
\EA
\end{equation}
Furthermore, if we assume either (i), or \smallskip

\noindent (iii) $p>\frac{N+1}{N-1}$ and $1<q< \frac{2p}{p+1}$ or \smallskip

\noindent (iv) $p>\frac{N+1}{N-1}$, $q= \frac{2p}{p+1}$ and
\begin{equation}\label{Z6}
\BA{lll}
M<m^{**}:=(p+1)\left(\myfrac{(N-1)p-(N+1)}{2p}\right)^{\frac{p}{p+1}},
\EA
\end{equation}

\noindent then $u=0$.
\es

This result is optimal in the case $p=\frac{N+1}{N-1}$, $q= \frac{2p}{p+1}$ as we will see in Section 4. Combining the method used in proving \rth{remov} with the result of \cite{MV1} we prove the removability of compact boundary sets on $\prt\Gw$, provided they satisfy some zero Bessel capacity property. 
\bth{remov-2} Assume $p>\frac{N+1}{N-1}$ and $\frac{N+1}{N-1}<r<p$. If
 one of the following conditions is satisfied:\par
\noindent (i) either $q=\frac{2p}{p+1}$ and 
  \bel{mr}
M<m^{**}_r:=(p+1)\left(\myfrac{p-r}{p(r-1)}\right)^{\frac{p}{p+1}},
\ee
\noindent (ii) or $1<q<\frac{2p}{p+1}$, $r\leq 3$ and $M$ is arbitrary. \smallskip

\noindent Then if $K\subset \prt\Gw$ is a compact set  such that $cap^{\prt\Gw}_{\frac{2}{r}, r'}(K)=0$, any solution $u$ of 
  \begin{equation}\label{RW15}\BA {lll}
-\Gd u+|u|^{p-1}u-M|\nabla u|^q= 0\qquad\text{in }\,\Gw\\
\phantom{-\Gd +|u|^{p-1}u-M|\nabla u|^q}
u=0\qquad\text{on }\,\prt\Gw\setminus K,
\EA\end{equation}
is identically $0$. 
\es

Note that $m^{**}_\frac{N+1}{N-1}=m^{**}$. The capacitary framework allows to consider the Dirichlet problem for $(\ref{Z1})$ 
\begin{equation}\label{NY1}
\BA{lll}
-\Gd u+|u|^{p-1}u-M|\nabla u|^{q}=0\qquad&\text{in }\Gw\\
\phantom{-\Gd +|u|^{p-1}u-M|\nabla u|^{q}}u=\gm&\text{in }\prt\Gw,
\EA
\end{equation}
where $\gm$ is a Radon measure on $\prt\Gw$. By a weak solution of $(\ref{NY1})$ we understand a function $u\in L^1(\Gw)\cap L^p_\gr(\Gw)$ such that 
$|\nabla u|\in L^q_\gr(\Gw)$, which satisfies
  \begin{equation}\label{MB7}\BA {lll}
\myint{\Gw}{}\left(-u\Gd\gz+(|u|^{p-1}u-M|\nabla u|^{q})\gz\right)dx=-\myint{\prt\Gw}{}\myfrac{\prt\gz}{\prt {\bf n}}d\gm\quad\text{for all }\gz\in \BBX(\Gw).
\EA\end{equation}

When the two exponents are super-critical with respect to the equations $(\ref{Z3})$ and $(\ref{Z4})$, the admissibility condition on the measure for $(\ref{Z1})$ requires the introduction of two different Bessel capacities defined on Borel subsets of $\prt\Gw$. 
\bth{adm-meas} Let  $p>1$, $1<q<2$ and $\gm$ be a nonnegative Radon measure on $\prt\Gw$ which satisfies 
  \begin{equation}\label{MB6-0}
\gm(E)\leq C\min\left\{cap^{\prt\Gw}_{\frac {2-q}{q},q'}(E),cap^{\prt\Gw}_{\frac 2p,p'}(E)\right\}\quad\text{for all Borel set }\,E\subset\prt\Gw,
\end{equation}
for some $C>0$. Then there exists  $c_0>0$ such that for any $0< c\leq c_0$ there exists a nonnegative weak solution of $(\ref{NY1})$ with boundary data $c\gm$. Furthermore the boundary trace of $u$ is the measure $c\gm$.
\es

The proof is based upon a non-standard application of the sub and supersolutions technique since it relies of the dynamical (and more natural) aspect of the boundary trace as it is exposed in \cite{MVbook}. Another surprising fact is the use of the equation 
$$
\BA{lll}
-\Gd u=u^p\qquad&\text{in }\Gw\\
\phantom{-\Gd }u=\gm&\text{in }\prt\Gw,
\EA
$$
which yields key estimates for our construction. The theorem admits several corollaries the proof of which is based on properties of Bessel capacities as exposed in \cite{AdHe}. It is noticeable that the results therein cover the full range $(p,q)\in (1,\infty)\ti (1,2)$.
\bcor{adm-meas-cor1} Assume  $p\geq \frac{N+1}{N-1}$ and $\frac{2p}{p+1}\leq q<2$. If $\gm$ is a nonnegative Radon measure on $\prt\Gw$ which satisfies, for some $C>0$, 
  \begin{equation}\label{MB6-1}
\gm(E)\leq Ccap^{\prt\Gw}_{\frac {2-q}{q},q'}(E)\quad\text{for all Borel set }\,E\subset\prt\Gw,
\end{equation}
there the conclusions of \rth{adm-meas} hold.
\es
The condition on the measure is also fulfilled under the following conditions.
\bcor{adm-meas-cor2} Assume $\frac{N+1}{N}\leq q<\frac{2p}{p+1}$. If $\gm$ is a nonnegative Radon measure on $\prt\Gw$ such that for some constant $C>0$, there holds for any Borel set
$E\subset\prt\Gw$, 
  \begin{equation}\label{MB6}\BA {lll}
\gm(E)\leq Ccap^{\prt\Gw}_{\frac 2p,p'}(E),
\EA\end{equation}
then the conclusions of \rth{adm-meas} hold. 
\es

Since the exponents $p$ and $q$ can be separately super or sub-critical, or even both sub-critical, we have the following result in different configurations of exponents. 
\bcor{subcritmeas} Let $p>1$, $1<q<2$ and $\gm\in\mathfrak M_+(\prt\Gw)$. There exists a function $u\in L^1(\Gw)\cap L^p_\gr(\Gw)$ such that $\nabla u\in L^q_\gr(\Gw)$ which is a weak solution to $(\ref{NY1})$ in the following cases: \smallskip

\noindent (i) When $p<\frac{N+1}{N-1}$, $q<\frac{N+1}{N}$ and there exists some $c_1>0$ such that
$\norm\gm_{\mathfrak M}\leq c_1$.

\noindent (ii) When $p< \frac{N+1}{N-1}$, $q\geq\frac{N+1}{N}$ and $\gm$ satisfies 
$(\ref{MB6-1})$; in that case $\gm$ has to be replaced by $c\gm$ with $0<c\leq c_2$,  for some $c_2>0$, in problem $(\ref{NY1})$. \smallskip

\noindent (iii) When $p\geq \frac{N+1}{N-1}$, $q<\frac{N+1}{N}$, and $\gm$ 
satisfies $\norm\gm_{\mathfrak M}\leq c_3$ for some  $c_3>0$ and
 \begin{equation}\label{MB12-1}\BA {lll}
\gm(E)=0\quad\text{for all Borel set $E\subset\prt\Gw$ such that }cap^{\prt\Gw}_{\frac {2}p,p'}(E)=0.
\EA\end{equation}
\es

In \cite{BVGHV3} the same authors study the problem
\begin{equation}\label{NYxx}
\BA{lll}
-\Gd u+|u|^{p-1}u-M|\nabla u|^{q}=\gm\qquad&\text{in }\Gw\\
\phantom{-\Gd +|u|^{p-1}u-M|\nabla u|^{q}}u=0&\text{in }\prt\Gw,
\EA
\end{equation}
where $\gm$ is a bounded Borel measure in $\Gw$. There too sufficient conditions for solving the problem involves Bessel capacities, but since the boundary trace argument is no longer valid, an intensive utilization of potential theory with various kernels has to be used. \medskip 

In the sub-critical case (i) and when $\gm$ is a Dirac mass at $0$ on the boundary we have no restriction on its weight. 
\bth{weak-sing} Assume $1<p<\frac{N+1}{N-1}$ and $1<q<\frac{N+1}{N}$. Then for any
$k\geq 0$ there exists a minimal positive solution $u_k$ of 
\begin{equation}\label{Z8}
\BA{lll}
-\Gd u+|u|^{p-1}u-M|\nabla u|^{q}=0\qquad&\text{in }\BBR^{_N}_+\\
\phantom{-\Gd +|u|^{p-1}u-M|\nabla u|^{q}}u=0&\text{in }\prt\BBR^{_N}_+\setminus \{0\},
\EA
\end{equation}
satisfying 
\begin{equation}\label{Z9}
\BA{lll}\displaystyle
\lim_{x\to 0}\myfrac{u_k(x)}{P_{_N}(x)}=k
\EA
\end{equation}
where $P_{_N}(x)=c_{_N}x_{{_N}}|x|^{-N}$ is the Poisson kernel in $\BBR^{_N}_+$. Furthermore this solution is unique among the positve solutions of $(\ref{Z8})$-$(\ref{Z9})$ if $q\leq\frac{2p}{p+1}$. This function satisfies $u_k\in L^1_{loc}(\overline{\BBR^{_N}_+})\cap L^p_{loc}(\overline{\BBR^{_N}_+};x_{{_N}}dx)$, 
$\nabla u_k\in L^q_{loc}(\overline{\BBR^{_N}_+};x_{{_N}}dx)$ and 
\begin{equation}\label{Z9a}
\BA{lll}\displaystyle
\myint{\BBR^{_N}_+}{}\left(-u_k\Gd \gz+(u_k^p-M|\nabla u_k|^q)\gz\right) dx=k\frac{\prt\gz}{\prt x_{{_N}}}(0)\quad \text{for all }\gz\in \BBX(\BBR^{_N}_+).
\EA
\end{equation}
\es
The proof is completely different from the ones of \rth{adm-meas} and \rcor{subcritmeas} and is based upon a delicate construction of supersolutions and subsolutions. A similar result holds if $\BBR^{_N}_+$ is replaced by a bounded smooth domain $\Gw\subset \BBR^{_N}_+$ such that $0\in\prt\Gw$.
\bth{weak-sing2} Assume $1<p<\frac{N+1}{N-1}$ and $0<q< \frac{N+1}{N}$. Then for any $M>0$ and $k>0$ there exists a minimal solution $u_k\in C^1(\overline\Gw\setminus\{0\})$ of $(\ref{NZ1})$ satisfying 
\begin{equation}\label{Z10}
\BA{lll}\displaystyle
\lim_{x\to 0}\myfrac{u_k(x)}{P_\Gw(x)}=k,
\EA
\end{equation}
where $P_\Gw$ is the Poisson kernel in $\Gw$. Furthermore $u_k\in L^1(\Gw)\cap L_\gr^p(\Gw)$, 
$\nabla u_k\in L^q_\gr(\Gw)$, and 
\begin{equation}\label{Z9b}
\BA{lll}\displaystyle
\myint{\Gw}{}\left(-u_k\Gd \gz+(u_k^p-M|\nabla u_k|^q)\gz\right) dx=-k\frac{\prt\gz}{\prt {\bf n}}(0)\quad \text{for all }\gz\in \BBX(\Gw).
\EA
\end{equation}
\es

In order to study the behaviour of these solutions $u_k$ when $k\to\infty$ we have to introduce separable solutions of $(\ref{Z1})$ in the model case $\BBR^{_N}_+$. They are solutions of 
\begin{equation}\label{Z11}
\BA{lll}
-\Gd u+|u|^{p-1}u-M|\nabla u|^{\frac{2p}{p+1}}=0\qquad&\text{in }\BBR^{_N}_+\\
\phantom{-\Gd +|u|^{p-1}u-M|\nabla u|^{\frac{2p}{p+1}}}u=0&\text{on }\prt\BBR^{_N}_+\setminus \{0\},
\EA
\end{equation}
which have the following expression in spherical coordinates
$$u(r,\gs)=r^{-\frac{2}{p-1}}\gw(\gs)\qquad\text{for all }(r,\gs)\in (0,\infty)\ti S_+^{_{N-1}}.
$$
Put 
\begin{equation}\label{alpha}\ga=\myfrac{2}{p-1},
\end{equation}
and denote by $\Gd'$ and $\nabla'$ the Laplace-Beltrami operator and the spherical gradient, then  $\gw$ satisfies 
\begin{equation}\label{Z12}
\BA{lll}
-\Gd' \gw+\ga(N-2-\ga)\gw+|\gw|^{p-1}\gw-M\left(\ga^2\gw^2+|\nabla' \gw|^2\right)^{\frac{p}{p+1}}=0\qquad&\text{in }S_+^{_{N-1}}\\
\phantom{-\Gd' \gw+\ga(N-2-\ga)\gw+|\gw|^{p-1}-M\left(\ga^2\gw^2+|\nabla' \gw|^2\right)^{\frac{p}{p+1}}}\gw=0&\text{in }\prt S_+^{_{N-1}}.
\EA
\end{equation}
 \bth{exist} There exists a positive solution $\gw$ to problem $(\ref{Z12})$ if one of the following conditions is satisfied:\smallskip
 
 \nind (i) either $1<p<\frac{N+1}{N-1}$ and $M\geq 0$,\smallskip
 
 \nind (ii) or $p=\frac{N+1}{N-1}$ and $M> 0$,\smallskip 
 
 \noindent (iii)  or $1<p<3$ or $p> \frac{N+1}{N-1}$, and $M\geq M_{_{N,p}}$  for some explicit value $M_{_{N,p}}>0$.
 \es

The positive solutions of $(\ref{Z12})$ allow to characterize the limit $u_\infty$ of the solutions  $u_k$ constructed in \rth{weak-sing}. 

\bth{souscrit-i} Let $1<p<\frac{N+1}{N-1}$, $1<q<\frac{N+1}{N}$ and $M>0$, then 
\bel{Z13}
\lim_{x\to 0}\myfrac{u_\infty(x)}{P_{_N}(x)}=\infty.
\ee
Furthermore
\smallskip

\nind (i) If $1<q<\frac{2p}{p+1}$
\bel{Z14}
\lim_{r\to 0} r^\ga u_\infty(r,.)=\psi\quad\text{uniformly on }S^{_{N-1}}_+,
\ee
where $\psi$ is the unique positive solution of 
 \begin{equation}\label{Z15}
\BA{lll}
-\Gd' \psi+\ga(N-2-\ga)\psi+|\psi|^{p-1}\psi=0\qquad&\text{in }S_+^{_{N-1}}\\
\phantom{-\Gd' \psi+\ga(N-2-\ga)\psi+|\psi|^{p-1}}\psi=0&\text{in }\prt S_+^{_{N-1}},
\EA
\end{equation}
and $u_\infty$ is the unique positive function solution of $(\ref{Z8})$ and satisfying $(\ref{Z13})$.\smallskip

\nind (ii) If $q=\frac{2p}{p+1}$
\bel{Z14a}
\lim_{r\to 0} r^\ga u_\infty(r,.)=\gw\quad\text{uniformly on }S^{_{N-1}}_+,
\ee
where $\gw$ is the minimal positive solution of $(\ref{Z12})$.
\es
A similar result holds if $\BBR^{_N}_+$ is replaced by a bounded smooth domain $\Gw\subset\BBR^{_N}_+$, which boundary contains $0$. In that case we assume that $T_{\prt\Gw}(0)=\prt\BBR^N_+$ (i.e. $\prt\BBR^N_+$ is the tangent hyperplane to $\prt\Gw$ at $0$ in order to use the spherical coordinates $(r,\gs)$ as above. Finally, if $(p,q)=\left(\frac{N+1}{N-1}, \frac{N+1}{N}\right)$ and $\prt\Gw$ is "very flat" near $0$ in the sense that $\dist(x,T_{\prt\Gw}(0))\leq c|x|^N$ for all $x\in \prt\Gw$ close to $0$, we prove that the function $u_\infty$ defined in the previous theorem 
still satisfies $(\ref{Z14a})$. Note that the above flatness condition is always satisfied if $N=2$ since $\prt\Gw$ is locally the graph of a $C^k$ real valued function ($k\geq 2$) defined on $T_{\prt\Gw}(0)\cap B_\gd$ and degenerate at $0$. \smallskip

When $\frac{2p}{p+1}<q<\min\{2,p\}$, the situation is completely changed and the solutions with strong boundary blow-up are modelized by equation $(\ref{SS15e})$. If $1<q<2$ we set
\begin{equation}\label{beta}\gb=\myfrac{2-q}{q-1},
\end{equation}
and if $1<q<p$
\begin{equation}\label{gamma}\gg=\myfrac{q}{p-q}.
\end{equation}

We prove the following result in the statement of which $\phi_1$ denotes the first eigenfunction of $-\Gd'$ in $W^{1,2}_0(S^{_{N-1}}_+)$. 
\bth{soupscrit} Assume $M>0$ and $\frac{2p}{p+1}<q<\min\{2,p\}$. Then there exists a positive solution $u$ of $(\ref{Z1})$ in $\BBR^{_N}_+$, which vanishes on $\prt\BBR^{_N}_+\setminus\{0\}$ such that 
  \bel{SS21c}
m\phi_1(\gs)r^{-\gg}\leq u(r,\gs)\leq c_4\max\left\{r^{-\ga}, M^{\frac{1}{p-q}}r^{-\gg}\right\}\;\text{ for all }\,(r,\gs)\in (0,r^*)\ti \overline{S^{_{N-1}}_+}.
\ee
for some $m>0$, $r^*\in (0,\infty]$ and where $c_4=c_4(N,p,q)>0$. If $Nq\geq (N-1)p$, $r^*=\infty$. 
\es

Note that our construction which is made by mean of supersolutions and subsolutions does not imply that in the case $\frac{2p}{p+1}<q<\frac{N+1}{N}$, the solution $u_\infty$ obtained in \rth{souscrit-i} satisfies $(\ref{SS21c})$. A similar result holds if $\BBR^{_N}_+$ is replaced by a bounded smooth domain $\Gw\subset\BBR^{_N}_+$, such that $0\in\prt\Gw$, under the flatness condition $\dist(x,T_{\prt\Gw}(0))\leq c|x|^{\gg+1}$ for $x\in\prt\Gw$ near $0$. \medskip

In the sequel $C>0$ denotes a constant the value of which can change from one occurence to another and $c_j$ ($j=0,1,2,...$) a more specific positive constant the value of which depends of more precise elements such as $p,q,N$ or other previous constants $c_i$. \medskip

In a forthcomming article \cite{BVGHV4} we study the isolated singularities of positive solution in a domains. Due to the 
number of parameters even the radial solutions present an amazing rich complexity.   \medskip

\nind{\bf Aknowledgements} The authors are grateful to the anonymous referee for his careful checking of the work and his pertinent observations. This article has been prepared with the support FONDECYT grants 1210241 and 1190102 for the three authors. 
\mysection{Singular boundary value problems}
\subsection{A priori estimates}
We give two series of estimates for solutions of $(\ref{Z1})$ with a boundary singularity according to the sign of $M$.
\bth {dom} Let $\Gw$ be a domain such that $0\in\prt\Gw$, $M\in\BBR$ and $1<q<\min\{p,2\}$. If $u\in C^{1}(\overline\Gw\setminus\{0\})$ is a solution of 
$(\ref{Z1})$ vanishing on $\prt\Gw\setminus\{0\}$, there holds\smallskip

\noindent 1- If $M>0$, there exists $=c_5(N,p,q)>0$ such that  
\begin{equation}\label{1Q1}
u_+(x)\leq c_5\max\left\{M^{\frac{1}{p-q}}|x|^{-\frac{q}{p-q}},|x|^{-\frac{2}{p-1}}\right\}\qquad\text{for all }\,x\in \Gw.
\end{equation}

\noindent 2- If  $M\leq 0$, there exist $c_6=c_6(N,q)>0$ and $c_7=c_7(N,p)>0$ such that 
\begin{equation}\label{1Q2}
u_+(x)\leq \min\left\{c_6|M|^{-\frac1{q-1}}|x|^{-\frac{2-q}{q-1}},c_7|x|^{-\frac{2}{p-1}}\right\}\qquad\text{for all }\,x\in \Gw.
\end{equation}
\es
\Proof We first  assume that $\overline \Gw\subset B_{R_0}$ for some $R_0>0$. Let $\ge>0$, we set 
$$j_\ge(r)=\left\{\BA{lll}
0\qquad&\text{if }\, r\leq 0\\
\frac{r^2}{2\ge}\qquad&\text{if }\, 0\leq r \leq\ge\\
r-\frac{\ge}{2}\qquad&\text{if }\, r \geq\ge.
\EA\right.
$$
If we extend $u$ by $0$ in $\overline\Gw^c\cap B_{2R_0}$ and set $v_\ge= j_\ge(u)$ we have
$$\BA {lll}
-\Gd v_\ge+v_\ge^p-M|\nabla v_\ge |^q=-j'_\ge(u)\Gd u-j''_\ge(u)|\nabla u|^2+(j_\ge(u))^p-M(j_\ge'(u))^q|\nabla u |^q\\[2mm]
\phantom{-\Gd v_\ge+v_\ge^p-M|\nabla v_\ge |^q}
\leq Mj'_\ge(u)\left(1-(j_\ge'(u))^{q-1}\right)|\nabla u |^q+(j_\ge(u))^p-j'_\ge(u)u_+^p\\[2mm]
\phantom{-\Gd v_\ge+v_\ge^p-M|\nabla v_\ge |^q}
\leq M\myfrac{u}{\ge}\left(1-\myfrac{u^{q-1}}{\ge^{q-1}}\right)|\nabla v_\ge |^q\chi_{_{\{0<u<\ge\}}}.
\EA$$
Letting $\ge\to 0$, we deduce from the dominated convergence theorem that $\displaystyle v_0=\lim_{\ge\to 0}v_\ge$ is nonnegative (actually it is the extension of $u^+$  by $0$ outside $\overline\Gw\setminus\{0\}$) and satisfies
\begin{equation}\label{1Q3}
Lv_0:=-\Gd v_0+v_0^p-M|\nabla v_0 |^q\leq 0\quad\text{in }\CD'(B_{2R_0}\setminus\{0\}).
\end{equation}

\noindent{\it The case $M>0$}. 
Following the method of Keller \cite{Ke} and Osserman \cite{Os}, we fix $a\in B_{R_0}\setminus\{0\}$, and introduce $U(x-a)=\gl(|a|^2-|x-a|^2)^{-b}$ for some 
$b>0$. Then putting $r=|x-a|$ and $\tilde U(r)=U(x-a)$, we have
$$\BA{lll}
L\tilde U=-\tilde U''-\myfrac{N-1}{r}\tilde U'-M|\tilde U'|^q+\tilde U^p\\[4mm]
\phantom{L\tilde U}
=\gl(|a|^2-r^2)^{-2-b}\left[\gl^{p-1}(|a|^2-r^2)^{2-b(p-1)}+2b(N-2(b+1))r^2-2Nb|a|^2\right.\\[4mm]
\phantom{L\tilde U-------------}
\left.-M2^qb^{q}\gl^{q-1}r^q(|a|^2-r^2)^{2+b-q(b+1)}\right].
\EA
$$
If $M>0$, the two necessary conditions on $b$ to be fulfilled in order $\tilde U$ be a supersolution in $B_{|a|}(a)$ are 
$$\BA {lll}&(i)\qquad\qquad 2-b(p-1)\leq 0\Longleftrightarrow b(p-1)\geq 2,\qquad\qquad\qquad\qquad\\[4mm]
&(ii)\qquad\qquad 2+b-q(b+1)\geq 2-b(p-1)\Longleftrightarrow b(p-q)\geq q.\qquad\qquad\qquad\qquad
\EA$$
The above inequalities are satisfied if
\begin{equation}\label{1Q3a}b= \max\left\{\myfrac{2}{p-1},\myfrac{q}{p-q}\right\}=\max\left\{\ga,\gg\right\}.
\end{equation}
If $q>\frac{2p}{p+1}$ then $b=\frac{q}{p-q}$ and 
$$L\tilde U\geq  \gl\left(|a|^2-r^2\right)^{-\frac{2p-q}{p-q}}\left[\gl^{q-1}\left(\gl^{p-q}-M2^qb^{q}r^q\right)\left(|a|^2-r^2\right)^{\frac{2p-q(p+1)}{p-q}}-(3b+1)N|a|^2\right].
$$
There exists $c_5>0$ depending on $N$, $p$ and $q$ such that if we choose
$$\gl= c_5\max\left\{M^{\frac{1}{p-q}}|a|^{\frac{q}{p-q}},|a|^{\frac{2p(q-1)}{(p-1)(p-q)}}\right\},
$$
there holds
\begin{equation}\label{1Q4}
L\tilde U\geq 0.
\end{equation}
Since $\tilde U(x)\to\infty$ when $|x|\to |a|$, we obtain by the maximum principle (see \cite{PS} for a ) that $v_0\leq \tilde U$ in $B_{|a|}(a)$. In particular 
\begin{equation}\label{1Q5}
u_+(a)= v_0(a)\leq \tilde U(a)=\gl|a|^{-\frac{2q}{p-q}}=c_5\max\left\{M^{\frac{1}{p-q}}|a|^{-\frac{q}{p-q}},|a|^{-\frac{2}{p-1}}\right\}.
\end{equation}
If $q\leq\frac{2p}{p+1}$ then $b=\frac{2}{p-1}$ and
$$\BA{lll}L\tilde U\geq  \gl\left(|a|^2-r^2\right)^{-\frac{2p}{p-1}}\left[\gl^{p-1}+\myfrac{2}{p-1}\left(N-\myfrac{2(p+1)}{p-1}\right)r^2-\myfrac{2N}{p-1}|a|^2\right.\\[4mm]\phantom{------------------}
\left.-M2^q\left(\myfrac{2}{p-1}\right)^q\gl^{q-1}r^q\left(|a|^2-r^2\right)^{\frac{2p-q(p+1)}{p-1}}\right]\\[4mm]\phantom{L\tilde U}
\geq  \gl\left(|a|^2-r^2\right)^{-\frac{2p}{p-1}}\left[\gl^{p-1}-C|a|^2-C'\gl^{q-1}M|a|^{\frac{4p-q(p+3)}{p-1}}\right].
\EA$$
Hence, if $q=\frac{2p}{p+1}$, $(\ref{1Q4})$ holds if for some $c_5=c_5(N,p,q)>0$, 
$$\gl= c_5\max\left\{M^{\frac{p+1}{p(p-1)}},1\right\}|a|^{\frac{2}{p-1}},
$$
which yields
\begin{equation}\label{1Q6}
u_+(a)= v_0(a)\leq \tilde U(a)=\gl|a|^{-\frac{4}{p-1}}=c_5\max\left\{M^{\frac{p+1}{p(p-1)}},1\right\}|a|^{-\frac{2}{p-1}}.
\end{equation}
While if $q<\frac{2p}{p+1}$, we choose 
$$\gl= c_5\max\left\{M^{\frac{1}{p-q}}|a|^{\frac{4p-q(p+3)}{(p-1)(p-q)}},|a|^{\frac{2}{p-1}}\right\},
$$
where $c_5>0=c_5(N,p,q)$, which yields
\begin{equation}\label{1Q6a}
u_+(a)= v_0(a)\leq \tilde U(a)=\gl|a|^{-\frac{4}{p-1}}=c_5\max\left\{M^{\frac{1}{p-q}}|a|^{-\frac{q}{p-q}},|a|^{-\frac{2}{p-1}}\right\}.
\end{equation}

\noindent{\it The case $M\leq 0$}. We first assume that $M< 0$. By \cite[Lemma 3.3]{NPT-V} $v_0$ satisfies
\begin{equation}\label{1Q7}
-\Gd v_0+|M||\nabla v_0|^q\leq 0\qquad\text{in }\,\CD'(B_{2R_0}\setminus\{0\}).
\end{equation}
Therefore, since $1<q<2$, 
\begin{equation}\label{1Q8}
u_+(a)=v_0(a)\leq c_6|M|^{-\frac{1}{q-1}}|a|^{-\frac{2-q}{q-1}}.
\end{equation}
If $M\leq  0$ there also holds 
\begin{equation}\label{1Q9}
-\Gd v_0+v_0^p\leq 0\qquad\text{in }\,\CD'(B_{2R_0}\setminus\{0\}).
\end{equation}
Hence 
\begin{equation}\label{1Q10}
u_+(a)=v_0(a)\leq c_7|a|^{-\frac{2}{p-1}}.
\end{equation}
In the above inequalities $c_6=c_6(q,N)>0$ and $c_7=c_7(p,N)>0$. Combining these estimates we derive
\begin{equation}\label{1Q11}
u_+(a)\leq \min\left\{c_7|a|^{-\frac{2}{p-1}},c_6|M|^{-\frac{1}{q-1}}|a|^{-\frac{2-q}{q-1}}\right\}.
\end{equation}
Since the estimate is independent of $R_0$, the assumption that $\overline \Gw\subset B_{R_0}$ is easily ruled out. This ends the proof.\qeda\medskip

\nind\Remark If $M= 0$, estimate $(\ref{1Q1})$ is just
\begin{equation}\label{1Q1+} u_+(x)\leq c_7|x|^{-\frac 2{p-1}}.
\end{equation}
If $M< 0$, $(\ref{1Q1+})$ is valid what ever is the value of $q$. Furthermore there also holds
\begin{equation}\label{1Q1++} u_+(x)\leq c_6|M||x|^{-\frac {2-q}{q-1}},
\end{equation}
whatever is the value of $p$, provided $1<q<2$.
\medskip

The equation is not invariant by $u\mapsto -u$, hence the lower and upper estimates are not symmetric.
\bcor{dom+} Under the assumptions of \rth{dom}, there holds\smallskip

\noindent 1- If $M>0$
\begin{equation}\label{1Q12}\BA {lll}
-\min\left\{c_6|M|^{-\frac1{q-1}}|x|^{-\frac{2-q}{q-1}},c_7|x|^{-\frac{2}{p-1}}\right\}\leq u(x)\\[2mm]
\phantom{----------}\leq c_5\max\left\{M^{\frac{1}{p-q}}|x|^{-\frac{q}{p-q}},|x|^{-\frac{2}{p-1}}\right\}\quad\text{for all }\,x\in \Gw.
\EA\end{equation}

\noindent 2- If  $M\leq 0$, there exist $c_6=c_6(N,q)>0$ and $c_7=c_7(N,p)>0$ such that 
\begin{equation}\label{1Q13}\BA {lll}
-c_5\max\left\{M^{\frac{1}{p-q}}|x|^{-\frac{q}{p-q}},|x|^{-\frac{2}{p-1}}\right\}\leq u(x)\\[2mm]
\phantom{---.-----}\leq \min\left\{c_6|M|^{-\frac1{q-1}}|x|^{-\frac{2-q}{q-1}},c_7|x|^{-\frac{2}{p-1}}\right\}\quad\text{for all }\,x\in \Gw.
\EA\end{equation}
\es

We infer from \rth{dom} an estimate of the gradient of $u$ near $0$.

\bth{domgrad} Let $\Gw$ be a smooth bounded domain such that $0\in\prt\Gw$ and $T_{\prt\Gw}(0)=\prt\BBR^N_+$, $M>0$, $p>1$ and $1<q<\min\{2,p\}$. If $u\in C^{1}(\overline\Gw\setminus\{0\})$ is a nonnegative solution of 
$(\ref{Z1})$ vanishing on $\prt\Gw\setminus\{0\}$, for any $r_0>0$ there holds there exists $c_8=c_8(N,p,q,\Gw,r_0,M)>0$ such that  
\begin{equation}\label{1Q14}
|\nabla u(x)|\leq c_8\max\left\{|x|^{-\frac{p}{p-q}},|x|^{-\frac{p+1}{p-1}}\right\}\qquad\text{for all }\,x\in \Gw\cap B_{r_0}.
\end{equation}
The restriction that $|x|\leq 1$ is not needed if $q=\frac{2p}{p+1}$.
\es
\Proof We assume first that $B_2^+\subset \Gw$.\smallskip 

\nind{\it Case 1: $1<q\leq\frac{2p}{p+1}$}. For $0<r<1$ we set 
$$u(x)=r^{-\frac{2}{p-1}}u_r(\tfrac{x}{r})=r^{-\frac{2}{p-1}}u_r(y)\quad\text{with }\,y=\tfrac{x}{r}.
$$
If $\frac{r}{2}<|x|<2r$, then $\frac{1}{2}<|y|<2$ and $u_r>0$ satisfies
$$-\Gd u_r+u_r^p-Mr^{\frac{2p-q(p+1)}{p-1}}|\nabla u_r|^q=0\qquad\text{in }\; B^+_2\setminus B^+_{\frac 12},
$$
and vanishes on $\prt (B^+_2\setminus B^+_{\frac 12})$. Since $0<Mr^{\frac{2p-q(p+1)}{p-1}}\leq M$ as $2p-q(p+1)\geq 0$, by the standard regularity theory we have the estimate
\begin{equation}\label{1Q15}\BA {lll}
\max\left\{|\nabla u_r(z)|: \frac{2}{3}<|z|<\frac{3}{2}\right\}\leq c_9\max\left\{|u_r(z)|: \frac{1}{2}<|z|<2\right\},
\EA\end{equation}
where $c_9$ depends on $N,p,q$ and $M$. Now it follows that
$$
\max\left\{|u_r(z)|: \tfrac{1}{2}<|z|<2\right\}\leq 2^{\frac{2}{p-1}}c_5\max\left\{M^{\frac{1}{p-q}}r^{\frac{2p-q(p+1)}{(p-1)(p-q)}},1\right\},
$$
 by $(\ref{1Q1})$. Therefore
\begin{equation}\label{1Q16}\BA {lll}
\max\left\{|\nabla u(y)|: \tfrac{r}{2}<|y|<2r\right\}\leq 2^{\frac{2}{p-1}}c_5c_9r^{-\frac{p+1}{p-1}}\max\left\{M^{\frac{1}{p-q}}r^{\frac{2p-q(p+1)}{(p-1)(p-q)}},1\right\}
\\[2mm]\phantom{\max\left\{|\nabla u(y)|: \tfrac{r}{2}<|z|<2r\right\}}
\leq c_8\max\left\{|x|^{-\frac{p}{p-q}},|x|^{-\frac{p+1}{p-1}}\right\},
\EA\end{equation}
which is $(\ref{1Q14})$. \smallskip 

\nind{\it Case 2: $\frac{2p}{p+1}<q<2$}. For $0<r<1$ we set 
$$u(x)=r^{-\frac{2-q}{q-1}}u_r(\tfrac{x}{r})=r^{-\frac{2-q}{q-1}}u_r(y)\quad\text{with }\,y=\tfrac{x}{r}.
$$
If $\frac{r}{2}<|x|<2r$, then $\frac{1}{2}<|y|<2$ and $u_r>0$ satisfies
$$-\Gd u_r+r^{\frac{q(p+1)-2p}{q-1}}u_r^p-M|\nabla u_r|^q=0\qquad\text{in }\; B^+_2\setminus B^+_{\frac 12},
$$
We notice that $q(p+1)-2p>0$. Then inequality $(\ref{1Q15})$ holds. Now 
$$
\max\left\{|u_r(z)|: \tfrac{1}{2}<|z|<2\right\}\leq c'_{9}r^{\frac{2-q}{q-1}}\max\left\{r^{-\frac{2}{p-1}}, r^{-\frac{q}{p-q}}\right\},
$$
thus
\begin{equation}\label{1Q15'}\BA {lll}
\max\left\{|\nabla u_r(z)|: \frac{2}{3}<|z|<\frac{3}{2}\right\}\leq c''_{9}r^{\frac{2-q}{q-1}-1}\max\left\{r^{-\frac{2}{p-1}}, r^{-\frac{q}{p-q}}\right\},
\EA\end{equation}
which implies
\begin{equation}\label{1Q15''}\BA {lll}
\max\left\{|\nabla u(x)|: \frac{2r}{3}<|x|<\frac{3r}{2}\right\}\leq c_{8}\max\left\{r^{-\frac{p+1}{p-1}}, r^{-\frac{p}{p-q}}\right\}.
\EA\end{equation}
\smallskip

\noindent \nind{\it  The general case; } If $\prt\Gw$ is not flat near $0$ we proceed as in the proof of \cite[Lemma 3.4]{NPT-V}, using the same scaling as in the flat case which transform the domain $B^+_2\setminus B^+_1)$ into $(B_2\setminus B_1)\cap \frac{1}{r}\Gw$, 
the curvature of which is bounded when $0<r<1$. The same estimates holds, up to the value of the constant $c_8$ and we derive $(\ref{1Q14})$. $\phantom{---------}$\qeda
\medskip

As a consequence we have the following.

\bcor{domgrad+} Under the assumptions of \rth{domgrad} the function $u$ satisfies
\begin{equation}\label{1Q117}
u(x)\leq c_8\gr(x)\max\left\{M^{\frac{1}{p-q}}|x|^{-\frac{p}{p-q}},|x|^{-\frac{p+1}{p-1}}\right\}\qquad\text{for all }\,x\in \Gw\cap B_1.
\end{equation}
The restriction that $|x|\leq 1$ is not needed if $q=\frac{2p}{p+1}$.
\es

 \subsection{Removable singularities}
 
 \noindent {\it Proof of \rth{remov}}.  If $M\leq 0$, $u$ is a nonnegative subsolution of 
$-\Gd u+v^p=0$ which vanishes on $\prt\Gw\setminus\{0\}$, hence it is identically zero by \cite{GmVe}. \smallskip

\noindent {\it Step 1. }We assume $M>0$. It is straightforward to verify from estimates $(\ref{1Q14})$ that under conditions (i) or (ii), $|\nabla u(x)|\leq c_8|x|^{-a}$ with $a\leq N$. Since these conditions imply $q<\frac{N+1}{N}$, it follows that $|\nabla u|^q\in L_\gr^1(\Gw)$.

For any $\ge>0$ we denote by $w_\ge$ the solution of 
  \begin{equation}\label{RW131}\BA {lll}
-\Gd w+w^p= M|\nabla u|^q\qquad&\text{in }\,\Gw_\ge:=\Gw\cap \overline B_\ge^c\\
\phantom{-\Gd +w^p}
w=0\qquad&\text{in }\,\prt\Gw\cap \overline B_\ge^c\\
\phantom{o}
\displaystyle \lim_{|x|\to\ge}w(x)=\infty\qquad&\text{on }\,\prt B_\ge\cap\Gw,
\EA\end{equation}
which exists since $|\nabla u|^q\in L_\gr^1(\Gw)$, see \cite{MV2}. Then $u\leq w_\ge$ in $\Gw_\ge$.  Let $z_\ge$ be the solution of 
  \begin{equation}\label{RW132a}\BA {lll}
-\Gd z+z^p=0\qquad&\text{in }\,\Gw_\ge\\
\phantom{-\Gd +w^p}
z=0\qquad&\text{in }\,\prt\Gw\cap \overline B_\ge^c\\
\phantom{o}
\displaystyle \lim_{|x|\to\ge}z(x)=\infty\qquad&\text{on }\,\prt B_\ge\cap\Gw.
\EA\end{equation}
Denote by $\BBG_\Gw[.]$ the Green operator in $\Gw$. Since $z_\ge+M\BBG_\Gw[|\nabla u|^q]\lfloor_{\Gw_\ge}$ is a supersolution of $(\ref{RW131})$ in $\Gw_\ge$  we deduce
  \begin{equation}\label{RW133-}\BA {lll}
u\leq z_\ge+M\BBG_\Gw[|\nabla u|^q]\lfloor_{\Gw_\ge}\qquad\text{in }\,\Gw_\ge.
\EA\end{equation}
When $\ge\to 0$, $z_\ge$ decreases to $z_0$ which satisfies 
  \begin{equation}\label{RW133}\BA {lll}
-\Gd w+w^p= 0\qquad&\text{in }\,\Gw\\
\phantom{-\Gd +w^p}
w=0\qquad&\text{in }\,\prt\Gw\setminus\{0\}.
\EA\end{equation}
 Since $p\geq\frac{N+1}{N-1}$ it is proved in\cite{GmVe} that any solution of  $(\ref{RW133})$ extends as a continuous solution in $\Gw$ with boundary value $0$, hence $z_0=0$ by the maximum principle. Therefore $u\leq M\BBG_\Gw[|\nabla u|^q]$ in $\Gw$ and the boundary trace $Tr_{\prt\Gw}[u]$ of $u$ is zero. By \cite{MVbook} the fact that $|\nabla u|^q\in L_\gr^1(\Gw)$ jointly with $Tr_{\prt\Gw}[u]=0$ implies in turn that $u^p\in L_\gr^1(\Gw)$ and $u$ is a weak solution of 
  \begin{equation}\label{RW132}\BA {lll}
-\Gd u+u^p= M|\nabla u|^q\qquad&\text{in }\,\Gw\\
\phantom{-\Gd +u^p}
u=0\qquad&\text{on }\,\prt\Gw,
\EA\end{equation}
in the sense that there holds
\begin{equation}\label{RW231}\BA {lll}
\myint{\Gw}{}\left(-u\Gd\gz+u^p\gz-M|\nabla u|^q\gz\right)dx=0\qquad\forall \gz\in W^{2,\infty}(\Gw)\cap C^1_c(\overline\Gw).
\EA\end{equation}

\noindent {\it Step 2. } Let us assume that $p>\frac{N+1}{N-1}$. If $u$ is nonnegative and not identically zero, then by the maximum principle it is positive in $\Gw$. We set $u=v^b$ with $0<b\leq 1$. Then 
 \begin{equation}\label{RW1}
\displaystyle
-\Gd v-(b-1)\myfrac{|\nabla v|^2}{v}+\myfrac{1}{b}v^{(p-1)b+1}-Mb^{q-1}v^{(b-1)(q-1)}|\nabla v|^q=0.
\end{equation}
 For $\ge>0$, 
 $$v^{(b-1)(q-1)}|\nabla v|^q\leq \myfrac{q\ge^{\frac{2}{q}}}{2}\myfrac{|\nabla v|^2}{v}+\myfrac{2-q}{2\ge^{\frac{2}{2-q}}}v^{\frac{(2b-1)q-2(b-1)}{2-q}}.
 $$
 Therefore
  \begin{equation}\label{RW2}
\displaystyle
-\Gd v+\left(1-b-M\myfrac{qb^{q-1}\ge^{\frac{2}{q}}}{2}\right)\myfrac{|\nabla v|^2}{v}+\myfrac{1}{b}v^{(p-1)b+1}-Mb^{q-1}\myfrac{2-q}{2\ge^{\frac{2}{2-q}}}v^{\frac{(2b-1)q-2(b-1)}{2-q}}=0.
\end{equation}
We notice that the following relation is independent of $b$
$$\frac{(2b-1)q-2(b-1)}{2-q}\leq (p-1)b+1\Longleftrightarrow q\leq \frac{2p}{p+1},
$$
with simultaneous equality. We take 
  \begin{equation}\label{RW3}(p-1)b+1=\frac{N+1}{N-1}\Longleftrightarrow b=\frac{2}{(N-1)(p-1)},
\end{equation}
 hence $p> \frac{N+1}{N-1}$ if and only if $0<b< 1$. \smallskip
 
 \noindent We first assume that $0<q<\frac{2p}{p+1}$ and choose $\ge>0$ such that
  \begin{equation}\label{RW4}
1-b-M\myfrac{qb^{q-1}\ge^{\frac{2}{q}}}{2}=0\Longleftrightarrow\ge=\left(\frac{2(1-b)}{Mqb^{q-1}}\right)^{\frac q2}=\left(\frac{2((N-1)p-N-1)}{Mqb^{q-1}(N-1)(p-1)}\right)^{\frac q2}\!\!\!.
\end{equation}
This transforms $(\ref{RW2})$ into 
  \begin{equation}\label{RW5}
\displaystyle
-\Gd v+\myfrac{(N-1)(p-1)}{2}v^{\frac{N+1}{N-1}}-\myfrac{(2-q)b^{\frac{2(q-1)}{2-q}}}{2}\left(\myfrac{q}{2(1-b)}\right)^{\frac{q}{2-q}}M^{\frac{2}{2-q}}v^{\frac{(2b-1)q-2(b-1)}{2-q}}\leq 0.
\end{equation}
Then, as  
$$\frac{(2b-1)q-2(b-1)}{2-q}<\frac{N+1}{N-1},$$
there exists $A>0$, depending on $M$ and $b$, such that 
  \begin{equation}\label{RW6}
\displaystyle
-\Gd v+\myfrac{(N-1)(p-1)}{4}v^{\frac{N+1}{N-1}}\leq A.
\end{equation}
Since $v$ vanishes on $\prt\Gw\setminus\{0\}$, $\tilde v=(v-c_{10}A^{\frac{N-1}{N+1}})_+^{\frac{N+1}{N-1}}$ with $c_{10}=\left(\frac{4}{(N-1)(p-1)}\right)^{\frac{N+1}{N-1}}$ satisfies
  \begin{equation}\label{RW7}
\displaystyle
-\Gd \tilde v+\myfrac{(N-1)(p-1)}{4}\tilde v^{\frac{N+1}{N-1}}\leq 0.
\end{equation}
By \cite{GmVe}, $\tilde v=0$ which implies $v\leq c_{10}A^{\frac{N-1}{N+1}}$ and therefore $u(x)\leq c_{11}A^{\frac{2}{(N+1)(p-1)}}$ in $\Gw$.   Since $u$ vanishes 
on $\prt\Gw\setminus\{0\}$ we extend it in a neighborhood of $0$ by odd reflection trough $\prt\Gw$ and denote by $\tilde u$ the new function defined in 
$B_\ga$ where it satisfies 
\begin{equation}\label{RW8}
\BA{lll}
-div A(x,\nabla\tilde u)+\tilde u^p+B(x,\nabla\tilde u)=0&\qquad\text{in }B_\ga \setminus\{0\}.
\EA
\end{equation}
In this expression the operator $A:(x;\xi)\in B_\ga\ti\BBR^{_N}\mapsto A(x,\xi)\in\BBR^{_N}$ is smooth in $x$ in $B_\ga$ and linear in $\xi$ (see e.g. \cite[Lemma 2.5]{BVGHV0} in a more general setting), it and satisfies for 
all $(x;\xi)\in B_\ga\ti\BBR^{_N}$,
$$\BA {lll}
A(x,\xi).\xi\geq 2|\xi|^2\;\;\text{and }\;|A(x,\xi)|\leq 4|\xi| \quad\text{for all }\; (x;\xi)\in B_\ga\ti\BBR^{_N}.
\EA$$
Since we can write $|B(.,\nabla\tilde u)|\leq 2|\nabla\tilde u|^q=2|\nabla\tilde u|^{q-1}|\nabla\tilde u|=C(x)|\nabla\tilde u|$ in $B_\ga$, then 
$B:(x,\xi)\in B_\ga\ti\BBR^{_N}\mapsto B(x,\xi)\in\BBR$ verifies
$$|B(x,\xi)|\leq C(x)|\xi|,
$$
and $C(x)\leq 2c_8|x|^{-\frac{(p+1)(q-1)}{p-1}}$ by \rth{domgrad}. Since $q<\frac{2p}{p+1}$, $\frac{(p+1)(q-1)}{p-1}<1$. Hence $C\in L^{N+\gt}$ 
for some $\gt>0$. By Serrin's theorem \cite[Theorem 10]{SerActa} the singularity at $0$ is removable and $\tilde u$ can be extended as a 
regular solution of $(\ref{RW8})$ in $B_\ga$. Hence $\tilde u\in C^1(B_{\frac\ga 2})$, and as a consequence $u\in C^1(\overline\Gw)$. 
If $u$ is not zero, it is positive in $\Gw$ and achieves its maximum at some $x_0\in\Gw$ where $\Gd u(x_0)\leq 0$ and $\nabla u(x_0)=0$. 
Contradiction. \smallskip

\noindent Next we assume that $q=\frac{2p}{p+1}$. By the choice of $b$ in  $(\ref{RW3})$, inequality  $(\ref{RW2})$ becomes 
  \begin{equation}\label{RW11}
\displaystyle
-\Gd v+\left(1-b-\myfrac{Mpb^{\frac{p-1}{p+1}}\ge^{\frac{p+1}{p}}}{p+1}\right)\myfrac{|\nabla v|^2}{v}+\left(\myfrac{1}{b}-
\myfrac{Mb^{\frac{p-1}{p+1}}}{(p+1)\ge^{p+1}}\right)v^{(p-1)b+1}\leq 0.
\end{equation}
\noindent We need to make both coefficients positive so that we obtain
 \begin{equation}\label{RW13}\BA {lll}
-\Gd v+\gt v^{\frac{N+1}{N-1}}\leq 0\qquad\text{in }\,\Gw\\
\phantom{-\Gd +\gt v^{\frac{N+1}{N-1}}}
v=0\qquad\text{on }\,\prt\Gw\setminus\{0\}.
\EA\end{equation}
We first choose 
$$\ge^{\frac{p+1}{p}}>\left(\myfrac{M}{p+1}\right)^{\frac 1p}b^{\frac{2}{p+1}},
$$
say
 \begin{equation}\label{RW13x}\ge^{\frac{p+1}{p}}=\left(\myfrac{M}{p+1}\right)^{\frac 1p}b^{\frac{2}{p+1}}+\tilde\ge,
\end{equation}
with $\tilde\ge>0$ so that the coefficient of $v^{\frac{N+1}{N-1}}$ is positive, and we can choose $\tilde\ge$ thanks to the assumption $m^{**}>M$: we have
  \begin{equation}\label{RW12}\BA {lll}
1-b-\myfrac{Mpb^{\frac{p-1}{p+1}}}{p+1}\left(\left(\myfrac{M}{p+1}\right)^{\frac 1p}b^{\frac{2}{p+1}}+\tilde\ge\right)
=1-b-\left(\myfrac{M}{p+1}\right)^{\frac {p+1}p}pb-\myfrac{Mpb^{\frac{p+1}{p-1}}}{p+1}\tilde\ge\\[4mm]
\phantom{--------------}
=b\left(\myfrac{1-b}{b}-\left(\myfrac{M}{p+1}\right)^{\frac {p+1}p}p\right)-\myfrac{Mpb^{\frac{p+1}{p-1}}}{p+1}\tilde\ge
\\[4mm]
\phantom{--------------}
=pb\left(\myfrac{(N-1)p-(N+1)}{2p}-\left(\myfrac{M}{p+1}\right)^{\frac {p+1}p}\right)-\myfrac{Mpb^{\frac{p+1}{p-1}}}{p+1}\tilde\ge
\\[4mm]
\phantom{--------------}
=pb\left(\left(\myfrac{m^{**}}{p+1}\right)^{\frac {p+1}p}-\left(\myfrac{M}{p+1}\right)^{\frac {p+1}p}\right)-\myfrac{Mpb^{\frac{p+1}{p-1}}}{p+1}\tilde\ge
\EA\end{equation}
and the right-hand side is positive if $\tilde\ge$ small enough. Hence we obtain $(\ref{RW13})$. By \cite{GmVe}, $v=0$ and the same holds for $u$. This ends the case $p>\frac{N+1}{N-1}$.\smallskip


\noindent {\it Step 3. } Finally we assume $p=\frac{N+1}{N-1}$ and $1<q<\frac{2p}{p+1}=\frac{N+1}{N}$, then  
$$M|\nabla u(x)|^q\leq c_{12}|x|^{-q\frac{p+1}{p-1}}=c_{12}|x|^{-qN}:=c_{13}Q(x).$$
Hence $u\leq u_1:=c_{13}\BBG_\Gw[Q]$. At this point we need the following intermediate result:\smallskip

\noindent{\it Claim.} Assume $w_\ga=\BBG_\Gw[Q_\ga]$ where $Q_\ga(x)=|x|^{-\ga}$ with $\ga<N+1$, then 
 \begin{equation}\label{RW14a}\BA {lll}
w_\ga(x)\leq c_{\ga}|x|^{2-\ga}\qquad\text{for all }\;x\in\Gw.
\EA\end{equation}
If this holds true, then $u(x)\leq c_{13}c_{qN}|x|^{2-qN}$. By the scaling method of \rth{domgrad}, we obtain
 \begin{equation}\label{RW15a}\BA {lll}
|\nabla u(x)|\leq c_8c_{13}c_{qN}|x|^{1-qN}\Longrightarrow |\nabla u(x)|^q\leq c_{14}|x|^{q(1-qN)}:=c_{14}Q_{q(Nq-1)}(x),
\EA\end{equation}
and thus 
 \begin{equation}\label{RW15b}
 w_{q(Nq-1)}(x)=c_{14}\BBG_\Gw[Q_{q(Nq-1)}](x)\leq c_{14}c_{q(Nq-1)}|x|^{2-q(Nq-1)}\qquad\text{for all }\;x\in\Gw.
\end{equation}
Since $q<1+\frac1N$, $q(Nq-1)-2<Nq-2$. Iterating this process, we finally obtain that $u$ is bounded and we end the proof as in Step 2. 
 $\phantom{-----}$\qeda\medskip
 
 \noindent\Remark It is noticeable that the equation exhibits a phenomenon which is characteristic of Emden-Folwer type equations
  \begin{equation}\label{RW16}
\Gd u=u^p\quad\text{in }\;B_1\setminus\{0\}.
\end{equation}
If $u$ is nonnegative then there exists $a\geq 0$ such that 
  \begin{equation}\label{RW16a}
\Gd u=u^p+a\gd_0\quad\text{in }\;\CD'(B_1).
\end{equation}
If $1<p<\frac{N}{N-2}$ then $a$ can be positive, but if $p\geq\frac{N}{N-2}$, then $a=0$. This means that the singularity cannont be seen 
in the sense of distributions, however there truly exist singular solutions, e.g. if $p>\frac{N}{N-2}$,
  \begin{equation}\label{RW17a}
u_s(x)=c_{N,p}|x|^{-\frac{2}{p-1}}.
\end{equation}
A similar phenomenon exists for solutions of 
  \begin{equation}\label{RW18}\BA {lll}
-\Gd u=u^p\quad\text{in }\;B^+_1\\
\phantom{-\Gd}
u=0\quad\text{in }\;\prt B^+_1\setminus\{0\}.
\EA\end{equation}
In such a case the critical value is $\frac{N+1}{N-1}$ since for $p\geq \frac{N+1}{N-1}$ the boundary value is achieved in the sense of distributions in $\prt B^+_1$.
\subsection{Proof of \rth{remov-2}}
 As in  \rth{remov}, the proof differs according to whether $0<q<\frac{2p}{p+1}$ or $q=\frac{2p}{p+1}$, and we first assume that $u>0$. We perform the same change of unknown as in the previous theorem putting $u=v^b$, but now we choose $b$ as follows
  \begin{equation}\label{RW16b}(p-1)b+1=r\Longleftrightarrow b=\frac{r-1}{p-1},
\end{equation}
and we first assume that 
  \begin{equation}\label{RW16aa}1-b-M\frac{qb^{q-1}\ge^{\frac2q}}{2}=0\Longleftrightarrow \ge=\left(\frac{2(1-b)}{Mqb^{q-1}}\right)^
\frac q2=\left(\frac{2(p-r)}{Mq(p-1)b^{q-1}}\right)^
\frac q2.
\end{equation}
Hence $(\ref{RW5})$ becomes  
  \begin{equation}\label{RW17}
\displaystyle
-\Gd v+\myfrac{p-1}{r-1}v^{r}-\myfrac{(2-q)b^{q-1}}{2}\left(\myfrac{q}{2(1-b)}\right)^{\frac{q}{2-q}}M^{\frac{2}{2-q}}v^{\frac{(2r-p-1)q+2(p-r)}{(p-1)(2-q)}}\leq 0.
\end{equation}
The condition $r\geq \frac{(2r-p-1)q+2(p-r)}{(p-1)(2-q)}$ is equivalent to $2p-q(p+1)\leq r(2p-q(p+1))$ since $1<r<p$. \par
\noindent Assuming first that $q<\frac{2p}{p+1}$, we obtain from $(\ref{RW17})$
  \begin{equation}\label{RW18a}
\displaystyle
-\Gd v+\myfrac{p-1}{2(r-1)}v^{r}\leq A.
\end{equation}
for some constant $A\geq 0$. Since $cap^{\prt\Gw}_{\frac{2}{r},r'}(K)=0$ and $v$ vanishes on $\prt\Gw\setminus K$, it follows from \cite{MV1} that 
$v\leq cA^{\frac1r}$ for some $c>0$, hence $u$ is also uniformly bounded above in $\Gw$ by some constant $a$. Next we have to show that $\nabla u\in L^2(\Gw)$. We also denote by $\Phi_1$ the first eigenfunction of $-\Gd$ in $W^{1,2}_0(\Gw)$ normalized by $\sup\Phi_1=1$ and by $\gl_1$ the corresponding eigenvalue. Since $\frac{N+1}{N-1}<r\leq 3$ we infer from \cite[Theorem 5.5.1]{AdHe}, that 
$$\left(cap^{\prt\Gw}_{\frac12,2}(K)\right)^{\frac{1}{N-2}}\leq B\left(cap^{\prt\Gw}_{\frac2r,r'}(K)\right)^{\frac{1}{N-1-\frac{2}{r-1}}}.
$$
Therefore $cap^{\prt\Gw}_{\frac{2}{r}, r'}(K)=0$ implies $cap^{\prt\Gw}_{\frac12,2}(K)=0$ and there exists a decreasing sequence $\{\gz_n\}\subset C^2_0(\prt\Gw)$ such that 
$\gz_n=1$ in a neighborhood of $K$, $0\leq\gz_n\leq 1$ and $\norm{\gz_n}_{W^{1,2}(\prt\Gw)}\to 0$ when $n\to\infty$, furthermore $\gz_n\to 0$ quasi everywhere. Let $\BBP_\Gw: C^2(\prt\Gw)\mapsto C^2(\overline\Gw)$ be the Poisson operator. It is an admissible lifting in the sense of \cite[Section 1]{MV1} in the sense that
$$\BBP_\Gw[\eta]\lfloor_{\prt\Gw}=\eta\quad \text {and }\;\eta\geq 0\Longrightarrow \BBP_\Gw[\eta]\geq 0.
$$
Put $\eta_n=1-\gz_n$. Then, multiplying equation $(\ref{RW15})$ by  $u(\BBP_\Gw[\eta_n])^{2}$ and integrating, we obtain
$$\BA {lll}\myint{\Gw}{}|\nabla u|^2(\BBP_\Gw[\eta_n])^{2} dx+2\myint{\Gw}{}u\BBP_\Gw[\eta_n]\nabla u. \nabla \BBP_\Gw[\eta_n]dx
\\[4mm]\phantom{-----}+\myint{\Gw}{}u^{p+1}(\BBP_\Gw[\eta_n])^{2}dx-M\myint{\Gw}{}|\nabla u|^qu(\BBP_\Gw[\eta_n])^{2} dx=0,
\EA$$
which implies
$$\BA {lll}
\myint{\Gw}{}|\nabla u|^2(\BBP_\Gw[\eta_n])^{2} dx-2\left(\myint{\Gw}{}|\nabla u|^2(\BBP_\Gw[\eta_n])^{2} dx\right)^{\frac12}\left(\myint{\Gw}{}|\nabla \BBP_\Gw[\eta_n]|^2u^{2} dx\right)^{\frac12}\\[4mm]\phantom{---------}+\myint{\Gw}{}u^{p+1}(\BBP_\Gw[\eta_n])^{2}dx-M\myint{\Gw}{}|\nabla u|^qu(\BBP_\Gw[\eta_n])^{2} dx\leq 0.
\EA$$
It is standard that 
$$\myint{\Gw}{}|\nabla \BBP_\Gw[\eta_n]|^2dx\leq c_{12}\norm{\eta_n}^2_{W^{\frac12,2}(\prt\Gw)}=A_n.
$$
Set $X_n=\norm{\BBP_\Gw[\eta_n]|\nabla u|}_{L^2}$, then
$$X_n^2-2A_nX_n-Ma|\Gw|^{\frac{2-q}{2}}X_n^q\leq 0.
$$
Hence there exist two positive real numbers $a_1$ and $a_2$ depending only on $q$, $|\Gw|$ and $a=\norm u_{L^\infty}$ such that
  \begin{equation}\label{RW19a}
\displaystyle
X_n\leq a_1A^\frac{1}{q-1}_n+a_2M^\frac{1}{2-q}.
\end{equation}
Now $A_n\to 0$ and $X_n\to\norm{\nabla u}_{L^2}^2$, therefore by Fatou's Lemma 
$$|\Gw|^{1-\frac{2}{q}}\norm{\nabla u}_{L^q}^2\leq \norm{\nabla u}_{L^2}^2\leq a_2M^\frac{1}{2-q}<\infty. 
$$
Let $\gz\in C^{1}_0(\overline\Gw)$ and $\eta_n$ as above. Since $\eta_n$ vanishes in a neighborhood of $K$ and $\gz$ vanishes on $\prt\Gw$,
$$\myint{\Gw}{}\BBP_\Gw[\eta_n]\nabla u.\nabla\gz dx+\myint{\Gw}{}\gz\nabla u.\nabla\BBP_\Gw[\eta_n] dx+\myint{\Gw}{}u^p\gz\BBP_\Gw[\eta_n]dx=M\myint{\Gw}{}|\nabla u|^q\gz\BBP_\Gw[\eta_n] dx.
$$
Letting $n$ to infty and using the fact that $\nabla u\in L^2(\Gw)$ and $\nabla \BBP_\Gw[\eta_n]\to 0$ in $L^2(\Gw)$, we derive 
$$\myint{\Gw}{}\nabla u.\nabla\gz dx+\myint{\Gw}{}u^p\gz dx=M\myint{\Gw}{}|\nabla u|^q\gz dx.
$$
Hence $u$ is a nonnegative bounded weak solution of 
  \begin{equation}\label{RW19}\BA {lll}
-\Gd u+|u|^{p-1}u-M|\nabla u|^q= 0\qquad\text{in }\,\Gw\\
\phantom{-\Gd +|u|^{p-1}u-M|\nabla u|^q}
u=0\qquad\text{on }\,\prt\Gw.
\EA\end{equation}
It is therefore $C^2$. Again, by the maximum principle we see that $u$ cannot achieve a positive maximum in $\Gw$, this yields a contradiction. \par

\noindent Next we assume $q=\frac{2p}{p+1}$.    We choose $b=\frac{r-1}{p-1}$ and $(\ref{RW11})$  becomes
  \begin{equation}\label{RW20}
\displaystyle
-\Gd v+\left(1-b-\myfrac{Mp b^{\frac{p-1}{p+1}}\ge^{\frac{p+1}{p}}}{p+1}\right)\myfrac{|\nabla v|^2}{v}+\left(\myfrac{1}{b}-\myfrac{Mb^{\frac{p-1}{p+1}}}{(p+1)\ge^{p+1}}\right)v^{r}\leq 0.
\end{equation}
From there the argument is similar to the one of Step 2-Case $q=\frac{2p}{p+1}$ in the proof of \rth{remov}: we claim that for some suitable choices the function $v$ satisfies 
$$\BA {lll} -\Gd v+\gt v^r\leq 0\qquad&\text{in }\;\Gw\\[1mm]
\phantom{-\Gd +\gt v^r}v=0&\text{in }\;\prt\Gw\setminus K.
\EA$$
We first choose $\ge>0$ so that $(\ref{RW13x})$ holds, hence the coefficient of $v$, say $\tau$ is positive. Then the expression
  \begin{equation}\label{RW21}\BA {lll}
1-b-\myfrac{Mp b^{\frac{p-1}{p+1}}\ge^{\frac{p+1}{p}}}{p+1}=\myfrac{p(r-1)}{p-1}\left(\left(\myfrac{m^{**}_r}{p+1}\right)^{\frac{p+1}{p}}-
\left(\myfrac{M}{p+1}\right)^{\frac{p+1}{p}}\right)-\myfrac{Mpb^{\frac{p+1}{p-1}}}{p+1}\tilde\ge
\EA\end{equation}
is positive provided $\tilde\ge>0$ is small enough. 
Since $cap^{\prt\Gw}_{\frac 2r,r'}(K)=0$ it follows from \cite{MV1} that $v=0$. Hence $u=0$, which ends the proof.\qeda
\subsection{Measure boundary data}
Let $\gm$ be a nonnegative Radon measure on $\prt\Gw$. The results concerning the following two types of equations 
  \begin{equation}\label{MB1}\BA {lll}
-\Gd v+v^p=0\qquad&\text{in }\;\Gw\\[0mm]
\phantom{-\Gd +v^p}v=\gm&\text{in }\;\prt\Gw,
\EA\end{equation}
and
  \begin{equation}\label{MB2}\BA {lll}
-\Gd w=M|\nabla w|^q\qquad&\text{in }\;\Gw\\[0mm]
\phantom{-\Gd }w=c\gm&\text{in }\;\prt\Gw,
\EA\end{equation}
allow us to consider the measure boundary data for equation $(\ref{Z1})$. We recall the results concerning $(\ref{MB1})$ and $(\ref{MB2})$.

\noindent 1- Assume  $p>1$. If $\gm$ satisfies 
  \begin{equation}\label{MB3}\BA {lll}
\text{For all  Borel set }E \subset\prt\Gw,\;cap^{\prt\Gw}_{\frac{2}{p},p'}(E)=0\Longrightarrow \gm(E)=0,
\EA\end{equation}
then problem $(\ref{MB1})$ admits a necessarily unique  weak solution $v:=v_\gm$, see \cite{MV1}, i.e.  $v_\gm\in L^1(\Gw)\cap L^p_\gr(\Gw)$ and for any function 
$\gz\in \BBX(\Gw):=\left\{\eta\in C^1_0(\overline\Gw)\text{ s.t. }\Gd \eta \in L^{\infty}(\Gw)\right\}$, 
there holds
  \begin{equation}\label{MB3*}\BA {lll}
\myint{\Gw}{}\left(-v\Gd\gz+v^p\gz\right) dx=-\myint{\Gw}{}\myfrac{\prt\gz}{\prt{\bf n}}d\gm.
\EA\end{equation}
Notice that there is no condition on $\gm$ if $1<p<\frac{N+1}{N-1}$.\smallskip

\noindent 2- Assume $1<q<2$. If there exists $C>0$ such that $\gm$ satisfies 
  \begin{equation}\label{MB4}\BA {lll}
\text{For all  Borel set }E \subset\prt\Gw,\;\gm(E)\leq Ccap^{\prt\Gw}_{\frac{2-q}{q},q'}(E),
\EA\end{equation}
then problem $(\ref{MB2})$ admits at least a positive solution $w$ for $c>0$ small enough, see \cite[Theorem 1.3]{BVHNV}, in the sense that $w\in L^{1}(\Gw)$, $\nabla w \in L_\gr^{q}(\Gw)$ and for any $\gz\in\BBX(\Gw)$, there holds
  \begin{equation}\label{MB4*}\BA {lll}
\myint{\Gw}{}\left(-w\Gd\gz-M|\nabla w|^q\gz\right) dx=-\myint{\Gw}{}\myfrac{\prt\gz}{\prt{\bf n}}d\gm.
\EA\end{equation}
Notice that if $1<q<\frac{N+1}{N}$ there is no capacitary condition on $\gm$. \medskip

We use also the following result.
\blemma{equiv} Let $p>\frac{N+1}{N-1}$ and $\gm\in\mathfrak M_+(\prt\Gw)$. If $\gm\in W^{-\frac{2}{p},p}(\prt\Gw)$, then there exists 
$C>0$ such that 
  \begin{equation}\label{MBlem1}\BA {lll}
 \gm(E)\leq C\,\left( cap^{\prt\Gw}_{\frac 2p,p'}(E)\right)^{\frac{1}{p'}}\quad\text{for all Borel set }E\subset\prt\Gw.
  \EA\end{equation}
  Conversely, if $\gm$ satisfies 
  \begin{equation}\label{MBlem2}\BA {lll}
 \gm(E)\leq Ccap^{\prt\Gw}_{\frac 2p,p'}(E)\quad\text{for all Borel set }E\subset\prt\Gw,
  \EA\end{equation}
  for some $C>0$, then $\gm\in W^{-\frac{2}{p},p}(\prt\Gw)$.
 \es
\Proof Assume $\gm\in W^{-\frac{2}{p},p}(\prt\Gw)\cap \mathfrak M_+(\prt\Gw)$. If $E$ is a compact subset of $\prt\Gw$ and $\gz\in C^2(\prt\Gw)$ with 
$0\leq\gz\leq 1$, with $\gz=1$ on $E$, then 
$$\gm(E)\leq \myint{\prt\Gw}{}\gz d\gm=\langle\gm,\gz\rangle\leq \norm\gm_{W^{-\frac{2}{p},p}}\norm\gz_{W^{\frac{2}{p},p'}}.
$$
Therefore, by the definition of the capacity,
$$\gm(E)\leq \norm\gm_{W^{-\frac{2}{p},p}}\left( cap^{\prt\Gw}_{\frac 2p,p'}(E)\right)^{\frac{1}{p'}}.
$$
Conversely, if $(\ref{MBlem2})$ holds, then there exists $c_{16}$ such that for any $0<c\leq c_{16}$ there exists a $z_{c\gm}$ to 
  \begin{equation}\label{MBX}\BA {lll}
-\Gd z=z^p\qquad&\text{in }\;\Gw\\[0mm]
\phantom{-\Gd}z=c\gm&\text{in }\;\prt\Gw,
\EA\end{equation}
(see \cite[Theorem 1.5]{BVHNV}) in the sense that $z_{c\gm}\in L^1(\Gw)\cap L^p_\gr(\Gw)$ and $c\BBP_\Gw[\gm]\leq z_{c\gm}$. Hence 
$\BBP_\Gw[\gm]\in L^p_\gr(\Gw)$, which implies $\gm\in W^{-\frac{2}{p},p}(\Gw)$ by \cite {MV1}.
\qeda\medskip

Those  weak solutions are characterized by their boundary trace. 
Let $\Gs_\ge=\{x\in\Gw:\gr(x)=\ge>0\}$ and $\Gs_0=\prt\Gw$. For $0<\ge\leq \ge_0$ the hypersurfaces $\Gs_\gd$ defines a foliation of the set $\Gw_{\ge_0}=\{x\in\Gw:\,0<\gr(x)\leq\ge_0\}$. Let $\gp(x)$ be the orthogonal projection of $x\in\Gw_{\ge_0}$ on $\prt\Gw$. Then $|x-\gp(x)|=\gr(x)$ and ${\bf n}_x=(\gr(x))^{-1}(\gp(x)-x)$. The mapping
$$x\mapsto \Gp(x)=(\gr(x),\gp(x)),
$$
from $\Gw_{\ge_0}$ onto $(0,\ge_0]\ti\Gs_0$ is a $C^2$ diffeomorphism and the restriction $\Gp_\ge$ of $\Gp$ to $\Gs_\ge$ is a $C^2$ diffeomorphism from 
 $\Gs_\ge$ onto $\Gs_0$. Let $dS_\ge$ be the surface measure on $\Gs_\ge$, then a continuous function $u$ defined in $\Gw$ has boundary trace the Radon measure $\gm$ on $\prt\Gw$ if 
  \begin{equation}\label{MB5*}\BA {lll}
\displaystyle \lim_{\ge\to 0}\myint{\Gs_\ge}{}uZ dS_\ge=\myint{\Gs}{}Z d\gm\qquad\text{for all }\;Z\in C(\overline\Gw).
\EA\end{equation}
Equivalently, if $\gz\in C(\prt\Gw)$ and $\gz_\ge=\gz\circ\Gp_\ge^{-1}\in C(\Gs_\ge)$, then 
  \begin{equation}\label{MB5**}\BA {lll}
\displaystyle \lim_{\ge\to 0}\myint{\Gs_\ge}{}u\gz_\ge dS_\ge=\myint{\Gs}{}\gz d\gm\qquad\text{for all }\;\gz\in C(\prt\Gw).
\EA\end{equation}
The functions $v_\gm$ solution of $(\ref{MB1})$ and $w$ solution of $(\ref{MB2})$ admit for respective boundary trace $\gm$ and $c\gm$. Furthermore, for the equations in $(\ref{MB1})$ and $(\ref{MB2})$,  the existence of a boundary trace of a positive solution is equivalent to the fact that $v_\gm\in L^1(\Gw)\cap L^p_\gr(\Gw)$ and $w\in L^{1}(\Gw)$ with $\nabla w \in L_\gr^{q}(\Gw)$ respectively. \medskip

\nind {\it Proof of \rth{adm-meas}}. If we assume that $(\ref{MB6-0})$ holds, the measure $\gm$ is  Lipschitz continuous with respect to $cap^{\prt\Gw}_{\frac 2p,p'}$ and $cap^{\prt\Gw}_{\frac {2-q}q,q'}$. By \cite[Theorem 1.3]{BVHNV} there exists $c_{17}>0$ such that for any $0<c\leq c_{17}$ there exists a  weak  solution $w=w_{c\gm}$ to $(\ref{MB2})$ and there holds for some positive constant $c_{18}$ depending on $q$ and $\Gw$
  \begin{equation}\label{MB8}w_{c\gm}\leq c_{18}c\BBP_\Gw[\gm]. 
\end{equation}
By \cite{MV1}  there exists a unique solution $v_{c\gm}$ to $(\ref{MB1})$ with $\gm$ replaced by $c\gm$. The functions $w_{c\gm}$ and 
$v_{c\gm}$ are respectively supersolution and subsolution of $(\ref{MB1})$ with boundary data  $c\gm$ and there holds, 
  \begin{equation}\label{MB8-1}
 v_{c\gm}\leq c \BBP_\Gw[\gm]\leq w_{c\gm}
\end{equation}
Hence there exists a nonnegative function $u$ satisfying $(\ref{Z1})$ and such that
  \begin{equation}\label{MB9}\BA {lll}
0\leq v_{c\gm}\leq u\leq w_{c\gm}\leq c_{18}c\BBP_\Gw[\gm].
\EA\end{equation}
Moreover $v_{c\gm}\in L^p_\gr(\Gw)$ and $\nabla w_{c\gm}\in L^q_\gr(\Gw)$. Because  $v_{c\gm}$ and $w_{c\gm}$ have boundary trace $c\gm$ in the sense of $(\ref{MB5*})$ and $(\ref{MB5**})$, the function $u$ has the same property and we denote it by $u_{c\gm}$. Assuming that 
$c\leq \min\{c_{16}, c_{17}\}$, there exists also $z_{c\gm}$ solution of $(\ref{MBX})$ which satisfies $z_{c\gm}\in L^p_\gr(\Gw)$ and 
$c\BBP_\Gw[\gm]\leq z_{c\gm}$ by the maximum principle. Therefore $w_{c\gm}\in L^p_\gr(\Gw)$ and finally $u_{c\gm}\in L^p_\gr(\Gw)$.
 
Let $\gf=\BBG_\Gw[u_{c\gm}^p]$, then $\gf\geq 0$ and 
$$-\Gd(u_{c\gm}+\gf)=|\nabla u_{c\gm}|^q.
$$
The function $u_{c\gm}+\gf$ is a nonnegative superharmonic function in $\Gw$. By Doob's theorem \cite[Chapter II]{Doob}, $-\Gd(u_{c\gm}+\gf)\in L^1_\gr(\Gw)$. Hence 
$|\nabla u_{c\gm}|\in L^q_\gr(\Gw)$. This implies that $u_{c\gm}$ is a  weak solution of $(\ref{MB7})$. \qeda


\medskip

\nind {\it Proof of \rcor{adm-meas-cor1}}. We use \cite[Theorem 5.5.1]{AdHe},  with the same cases (a), (b), (c) and (d), and we denote by $K$ is any compact subset of $\prt\Gw$ and by $A$ a positive constant the value of which may vary from one case to another. \smallskip

\nind (a) If $q>\frac{2p}{p+1}$ and $p>\frac{N+1}{N-1}$, equivalently $\frac{2-q}{q-1}<\frac{2}{p-1}<N-1$, then 
 \begin{equation}\label{MB8-a}\displaystyle
cap^{\prt\Gw}_{\frac {2-q}q,q'}(K)\leq A\left(cap^{\prt\Gw}_{\frac {2}p,p'}(K)\right)^\frac{N-1-\frac{2-q}{q-1}}{N-1-\frac{2}{p-1}}.
  \end{equation}
Since $\frac{N-1-\frac{2-q}{q-1}}{N-1-\frac{2}{p-1}}>1$, thus 
  \begin{equation}\label{MB8-a1}
  cap^{\prt\Gw}_{\frac {2-q}q,q'}(K)\leq c_{21}cap^{\prt\Gw}_{\frac {2}p,p'}(K).
  \end{equation}
If $(\ref{MB6-1})$ holds, then 
$$\BA {lll}
\gm(E)\leq C cap^{\prt\Gw}_{\frac {2-q}q,q'}(K)= 
C\min\left\{cap^{\prt\Gw}_{\frac {2-q}q,q'}(K),c_{21} cap^{\prt\Gw}_{\frac {2}p,p'}(K)\right\}\\[2mm]
\phantom{\gm(E)\leq C cap^{\prt\Gw}_{\frac {2-q}q,q'}(K)}
\leq C\max\{1,c_{21}\}\min\left\{cap^{\prt\Gw}_{\frac {2-q}q,q'}(K),cap^{\prt\Gw}_{\frac {2}p,p'}(K)\right\}
\EA$$
and the proof follows by \rth{adm-meas}.\smallskip

\nind (b) If $q=\frac{2p}{p-1}$ and $p>\frac{N+1}{N-1}$, then $p>q$, thus 
  \begin{equation}\label{MB8-b}
cap^{\prt\Gw}_{\frac {2-q}q,q'}(K)\leq c_{22}cap^{\prt\Gw}_{\frac {2}p,p'}(K).
  \end{equation}
  The proof follows as in (a). \smallskip

\nind (c) If $p=\frac{N+1}{N-1}$ and $q>\frac{2p}{p+1}$, then for some $A>cap^{\prt\Gw}_{\frac {2-q}q,q'}(\prt\Gw)$, 
  \begin{equation}\label{MB8-c}\displaystyle
\left(\ln \myfrac{A}{cap^{\prt\Gw}_{\frac {2-q}q,q'}(K)}\right)^{-1}\leq A\left(cap^{\prt\Gw}_{\frac {2}p,p'}(K)\right)^{\frac{2}{N-1}}.
  \end{equation}
  Since for any $r\geq 1$
  $$\left(\ln r\right)^{-1}>\myfrac{2}{N-1} r^{-\frac{2}{N-1}},
  $$
  we deduce
    \begin{equation}\label{MB8-c1}
cap^{\prt\Gw}_{\frac {2-q}q,q'}(K)\leq \left(\myfrac{N-1}{2}\right)^{\frac{N-1}{2}}A^{\frac{N+1}{2}}cap^{\prt\Gw}_{\frac {2}p,p'}(K)
:=c_{23}cap^{\prt\Gw}_{\frac {2}p,p'}(K).
  \end{equation}
  The proof follows as in (a). \smallskip

\nind (d) If $p=\frac{N+1}{N-1}$ and $q=\frac{2p}{p+1}=\frac{N+1}{N}$, then as above
  \begin{equation}\label{MB8-d}\displaystyle
\left(cap^{\prt\Gw}_{\frac {2-q}q,q'}(K)\right)^{q-1}\!\!\leq A\left(cap^{\prt\Gw}_{\frac {2}p,p'}(K)\right)^{p-1}\!\!\Longrightarrow
cap^{\prt\Gw}_{\frac {2-q}q,q'}(K)\leq c_{24}cap^{\prt\Gw}_{\frac {2}p,p'}(K).
  \end{equation}
  The proof follows.

\qeda
\medskip

\nind {\it Proof of \rcor{adm-meas-cor2}}. We adapt again the formulation of \cite[Theorem 5.5.1]{AdHe} to our framework permuting the two capacities and only statement (a) and (c) therein apply. \smallskip

\nind (a) If $\frac{N+1}{N}< q<\frac{2p}{p+1}$ there exists a constant $A>0$ such that if $K\subset\prt\Gw$ is a compact set, then
  \begin{equation}\label{MB7a}\BA {lll}
cap^{\prt\Gw}_{\frac 2p,p'}(K)\leq A\left(cap^{\prt\Gw}_{\frac {2-q}q,q'}(K)\right)^{\frac{N-1-\frac{2}{p-1}}{N-1-\frac{2-q}{q-1}}}.
\EA\end{equation}
Since $\frac{N+1}{N}<q<\frac{2p}{p+1}$ is equivalent to  $N-1-\frac{2}{p-1}>N-1-\frac{2-q}{q-1}>0$, we deduce 
  \begin{equation}\label{MB6-4}
  cap^{\prt\Gw}_{\frac 2p,p'}(K)\leq c_{25}cap^{\prt\Gw}_{\frac {2-q}q,q'}(K).
\ee
We end the proof as in the proof of \rcor{adm-meas-cor1}-(a).\smallskip

\nind (c) If $\frac{N+1}{N}= q<\frac{2p}{p+1}$, then for some $A>cap^{\prt\Gw}_{\frac {2}p,p'}(\prt\Gw)$, 
  \begin{equation}\label{MB8-c2}\displaystyle
\left(\ln \myfrac{A}{cap^{\prt\Gw}_{\frac {2}p,p'}(K)}\right)^{-1}\leq A\left(cap^{\prt\Gw}_{\frac {2-q}q,	q'}(K)\right)^{\frac1N}.
  \end{equation}
Since for $r>1$, 
 $$\left(\ln r\right)^{-1}>\frac1N r^{-\frac1N},
  $$
  we infer
    \begin{equation}\label{MB8x-c}\displaystyle
cap^{\prt\Gw}_{\frac {2}p,p'}(K)\leq N^NA^{N+1}cap^{\prt\Gw}_{\frac {2-q}q,q'}(K):=c_{26}cap^{\prt\Gw}_{\frac {2-q}q,q'}(K),
  \end{equation}
  and the proof follows.
 \qeda
\medskip


The proof in the partially sub-critical case is simpler. \medskip

\nind {\it Proof of \rcor{subcritmeas}}. If $1<p<\frac{N+1}{N-1}$ for any $\gm\in\mathfrak M_+(\prt\Gw)$ problem $(\ref{MB1})$ admits a unique solution $v_\gm$ (see \cite{GmVe}). If $1<q<\frac{N+1}{N}$, then 
there exists $a_0>0$ such that for any non-empty Borel set $E\subset \prt\Gw$, $cap^{\prt\Gw}_{\frac {2-q}q,q'}(E)\geq a_0$. Therefore 
$$\gm(E)\leq \norm\gm_{\mathfrak M}\leq \myfrac{\norm\gm_{\mathfrak M}}{a_0}cap^{\prt\Gw}_{\frac {2-q}q,q'}(E).
$$
It follows from \cite[Theorem 1.3]{BVHNV} that problem $(\ref{MB2})$ admits a solution $w_\gm$ whenever $\norm\gm_{\mathfrak M}$ is small enough. By \cite{BVVi}  problem $(\ref{MBX})$ admits a solution $z_\gm$ with $c\gm$ replaced by $\gm$ provided $\norm\gm_{\mathfrak M}$ is small enough. Furthermore 
 \begin{equation}\label{MB13}w_{\gm}\leq \BBP_\Gw[\gm]\leq z_\gm.\end{equation}
Since $z_\gm\in L^p_\gr(\Gw)$, $w_{\gm}\in L^p_\gr(\Gw)$. Hence by the same arument as in \rth{adm-meas}, there exists a solution $u_{\gm}$ of $(\ref{Z1})$ which satisfies $v_\gm\leq u_\gm\leq w_\gm$. Hence $u_{\gm}\in L^p_\gr(\Gw)$ and by the previous argument $\nabla u_\gm\in L^q_\gr(\Gw)$. This implies again that $u_\gm$ is a  weak solution of $(\ref{MB7})$. \smallskip

\noindent  If $1<p<\frac{N+1}{N-1}$ and $\frac{N+1}{N}\leq q<2$, then problem $(\ref{MB1})$ is uniquely solvable for any  $\gm\in\mathfrak M_+(\prt\Gw)$,   while problem $(\ref{MBX})$ admits a solution $z_\gm$ with $c\gm$ replaced by $\gm$ provided $\norm\gm_{\mathfrak M}$ is small enough and since $(\ref{MB4})$ holds, problem $(\ref{MB2})$ admits a  weak solution provided $0<c\leq c_0$. Since $(\ref{MB13})$ holds with $z_\gm\in L^p_\gr(\Gw)$, the result follows as above.\smallskip

\noindent If $p\geq \frac{N+1}{N-1}$, $1<q<\frac{N+1}{N}$ and $\gm\in\mathfrak M_+(\prt\Gw)$ absolutely continuous with respect to $cap^{\prt\Gw}_{\frac {2}p,p'}$, there exists $u_\gm$ solution of  $(\ref{MB1})$ and $w_\gm$ solution of $(\ref{MB2})$ provided $c\norm\gm_{\mathfrak M}$ is small enough. Since 
$|\nabla w_\gm|^q\in L^1_\gr(\Gw)$ the function $w_\gm$ belongs to the Marcinkiewicz space $M^{\frac{N+1}{N-1}}_\gr(\Gw)$ (see eg. \cite{Vebook1}). Since 
$M^{\frac{N+1}{N-1}}_\gr(\Gw)\subset L^p_\gr(\Gw)$ as $1<p<\frac{N+1}{N-1}$, it implies that $w_\gm$ and therefore $u_\gm$, belongs to 
$L^p_\gr(\Gw)$. The end of the proof is as before.$\phantom{----------}$\qeda

\mysection{Separable solutions}
Separable solutions of $(\ref {Z1})$ in $\BBR^{_N}\setminus\{0\}$ are solutions which have  the form 
$$u(x)=u(r,\gs)=r^{-\gk}\gw(\gs)\qquad\text{for }\, (r,\gs)\in \BBR_+\ti S^{_{N-1}}.
$$
This forces $q$ to be equal to $\frac{2p}{p+1}$, $\gk$ to $\frac{2}{p-1}$ (recall that this defines $\ga$) and $\gw$ satisfies
\begin{equation}\label{Z2}
\BA{lll}
-\Gd' \gw+\ga(N-2-\ga)\gw+|\gw|^{p-1}\gw-M\left(\ga^2\gw^2+|\nabla' \gw|^2\right)^{\frac{p}{p+1}}=0\quad&\text{in }\;S^{_{N-1}}.
\EA
\end{equation}
Constant positive solutions are solutions of
\bel{fi1}X^{p-1}-M\ga^{\frac{2p}{p+1}}X^{\frac{p-1}{p+1}}+\ga(N-2-\ga)=0.
\ee
This existence of solutions to $(\ref{fi1})$ and their stability properties will be detailled in a forthcoming article \cite{BVGHV4}. The understanding of boundary singularities of solutions of $(\ref{Z1})$ is conditioned by the knowledge 
of separable solutions in $\BBR^{_N}_+$ vanishing on $\prt\BBR^{_N}\setminus\{0\}$. Then $\gw$ is a solution of 
\begin{equation}\label{R1}
\BA{lll}
-\Gd' \gw+\ga(N-2-\ga)\gw+|\gw|^{p-1}\gw-M\left(\ga^2\gw^2+|\nabla' \gw|^2\right)^{\frac{p}{p+1}}=0\quad&\text{in }\;S^{_{N-1}}_+\\[2mm]
\phantom{-\Gd' +\ga(N-2-\ga)\gw+|\gw|^{p-1}\gw-M\left(\ga^2\gw^2+|\nabla' \gw|^2\right)^{\frac{p}{p+1}}}
\gw=0\quad&\text{in }\;\prt S^{_{N-1}}_+.
\EA
\end{equation}

\subsection{Existence of singular solutions}
We recall the following result proved in \cite[Corollary 1.4.5]{Vebook} is a variant of Boccardo-Murat-Puel's result  \cite[Theorem 2.1]{BMP} dealing with the quasilinear equation in a domain $G\subset\BBR^{_N}$. 
\bel{ChI-3'-3}\BA {ll}
\CQ(u):=-\Gd u+B(.,u,\nabla u)=0\qquad\text{in }\CD'(G),
\EA\ee
where $B\in C(G\ti\BBR\ti\BBR^{_N})$ satisfies, for some continuous increasing function $\Gg$ from $\BBR^+$ to $\BBR^+$,
\bel{I-3'-2}\BA {ll}
\abs{B(x,r,\xi)}\leq \Gg(|r|)(1+|\xi|^2)\quad\text{for all }\,(x,r,\xi)\in G\ti\BBR\ti\BBR^{_N}.
\EA\ee

\bth{BMP} Let $G$ be a bounded domain in $\BBR^{_N}$. If there exist a supersolution $\gf$ and a subsolution $\psi$ of the equation $\CQ v=0$
belonging to $W^{1,\infty}(G)$ and such that $\psi\leq\gf$, then for any $\chi\in W^{1,\infty}(G)$ satisfying 
$\psi\leq\chi\leq\gf$ there exists a function $u\in W^{1,2}(G)$ solution of $\CQ u=0$ such that $\psi\leq u\leq\gf$ 
and $u-\chi\in W^{1,2}_0(G)$.\es

\noindent\Remark {\it Mutatis mutandi}, the same result holds if $\BBR^{_N}$ is replaced by a Riemannian manifold.\medskip

Their result is actually more general since the Laplacian can be replaced by a quasilinear $p$-Laplacian-type operator and $B$ by a perturbation with the {\it natural} $p$-growth.  This theorem has direct applications in the construction of solutions on $S^{_{N-1}}_+$, but also for the construction of singular solutions 
in several configurations.

\bprop {keylem} Let $\Gw$ be a bounded smooth domain containing $0$, $p>1$, $1\leq q\leq 2$ and $M\in\BBR$. Assume that equation 
\bel{K1}
-\Gd u+u^p-M|\nabla u|^q=0,
\ee
admits a radial positive and decreasing solution $v$ in $\BBR^{_N}\setminus\{0\}$ satisfying
\bel{K2}\displaystyle
\lim_{|x|\to 0}v(x)=\infty.
\ee
Then there exists a positive function $u$ satisfying $(\ref{K1})$ in $\Gw\setminus\{0\}$, vanishing on $\prt\Gw$ and such that 
\bel{K3}\displaystyle
\left(v(x)-\max\left\{v(z):|z|=\gd_0\right\}\right)_+\leq u(x)\leq v(x)\quad\text{for all }\,x\in \Gw\setminus\{0\}.
\ee
where $\gd_0=\dist(0,\prt\Gw)$.
\es
\Proof Put $m=\max\left\{v(z):|z|=\gd_0\right\}=v(\gd_0)$. The function $v_m=(v-m)_+$ is a radial subsolution of $(\ref{K1})$ in $\Gw$, positive in 
$B_{\gd_0}\setminus\{0\}$ and vanishing in $\Gw\setminus B_{\gd_0}$. For $\ge>0$ set $\Gw_\ge=\Gw\setminus \overline B_{\ge}$. Since $v_m$ is  dominated by the supersolution $v$,  there exists a solution $u_\ge$ of $(\ref{K1})$ in $\Gw_\ge$ such that $v_m\leq u_\ge\leq v$ and $u_\ge-v_m\in H^{1}_0(\Gw_\ge)$. By standard regularity estimates, $u_\ge$ is $C^2$, hence it solves
\bel{K4}\BA {lll}
-\Gd u_\ge+u_\ge^p-M|\nabla u_\ge|^q=0&\quad\text{in }\,\Gw_\ge\\
\phantom{-\Gd +u_\ge^p-M|\nabla u_\ge|^q}
u_\ge=v_m&\quad\text{on }\,\prt B_\ge
\\
\phantom{-\Gd +u_\ge^p-M|\nabla u_\ge|^q}
u_\ge=0&\quad\text{on }\,\prt \Gw.
\EA\ee
Notice that $u_\ge$ is unique by the comparison principle. If $0<\ge'<\ge$ the function $u_{\ge'}$ solution of (\ref{K4}) in $\Gw_{\ge'}$ with the corresponding boundary data is larger than 
$v_m$ and in particular $u_{\ge'}\lfloor_{\prt B_\ge}\geq v_m\lfloor_{\prt B_\ge}=u_\ge\lfloor_{\prt B_\ge}$. Hence $u_{\ge'}\geq u_\ge$ in 
$\Gw_\ge$. When $\ge\downarrow 0$, $u_\ge$ increase and converges in the $C^{1,\theta}_{loc}(\overline\Gw\setminus\{0\})$-topology toward some function $u$  which satisfies $(\ref{K1})$  in $\Gw\setminus\{0\}$, is larger that $v_m$ and smaller than $v$, vanishes on $\prt\Gw$ and such that $(\ref{K4})$ holds.\qeda
\medskip

The previous result can be adapted to the study of solutions with a boundary singularity in bounded domains which are flat enough near the singular point or in $\BBR^{_N}_+$.

\bprop {keylem2} Let $p>1$, $1\leq q\leq 2$ and $M\in\BBR$. Assume that the equation $(\ref{K1})$
admits a positive solution $w$ in $\BBR^{_N}_+$ belonging to  $C(\overline{\BBR^{_N}_+}\setminus\{0\})$, radially decreasing in $\BBR^{_N}_+$ and satisfying 
\bel{K5}\BA {lll}\displaystyle
 \lim_{t\to 0}w(t\gs)=\infty\quad \text{uniformly on compact sets }\;K\subset S^{_{N-1}}_+.
\EA\ee
Assume also \smallskip

\nind(i)  either $w\lfloor_{\prt\BBR^{_N}_+\setminus\{0\}}$ is bounded, 
\smallskip

\nind(ii)  or $\Gw\subset\BBR^{_N}_+$ is a bounded smooth domain such that $0\in\prt\Gw$ starshapped with respect to $0$
and such that $w\lfloor_{\prt\Gw\setminus\{0\}}$ is bounded. \smallskip

\nind Then there exists a positive function $u$ satisfying $(\ref{K1})$ in $\BBR^{_N}_+$ in case (i), or $\Gw$ in case (ii) , vanishing on $\prt\BBR^{_N}_+\setminus\{0\}$  in case (i), or $\prt\Gw\setminus\{0\}$ in case (ii), and such that 
\bel{K6-}\displaystyle
\left(w(x)-\sup\left\{w(z):z\in \prt\BBR^{_N}_+\setminus\{0\}\right\}\right)_+\leq u(x)\leq w(x)\quad\text{for all }\,x\in \BBR^{_N}_+,
\ee
where $K=\sup\left\{\displaystyle \limsup_{|z|\to\infty}w(z),\sup\left\{w(z):z\in \prt\BBR^{_N}_+\setminus\{0\}\right\}\right\}$ in case (i) or 
\bel{K6}\displaystyle
\left(w(x)-\sup\left\{w(z):z\in \prt\Gw\setminus\{0\}\right\}\right)_+\leq u(x)\leq w(x)\quad\text{for all }\,x\in \Gw.
\ee
in case (ii).
\es
\Proof The proof is a variant of the preceding one, only the geometry of the domains is changed. \smallskip

\noindent In case (ii) set $m=\sup\left\{w(z):z\in \prt\Gw\setminus\{0\}\right\}$. Then the function $z\mapsto w_m:=(w(z)-m)_+$ is  a subsolution of $(\ref{K1})$ in $\Gw$. It vanishes on $\prt\Gw\setminus\{0\}$ and is dominated by $w$. For $\ge<\gd_0$, let $\Gw_\ge$ denote $\Gw\cap\overline {B_\ge}^c$. We consider the problem of finding $u_\ge$ solution of 
\bel{K7}\BA {lll}
-\Gd u_\ge+u_\ge^p-M|\nabla u_\ge|^q=0&\quad\text{in }\,\Gw_\ge\\
\phantom{-\Gd +u_\ge^p-M|\nabla u_\ge|^q}
u_\ge=w_m&\quad\text{on }\,\prt B_\ge\cap\Gw\\
\phantom{-\Gd +u_\ge^p-M|\nabla u_\ge|^q}
u_\ge=0&\quad\text{on }\, B^c_\ge\cap\prt\Gw.
\EA\ee
Again since $u_\ge-w_m\in H^1_0(\Gw_\ge)$ and since $w_m$ is smaller than $w\lfloor_{\Gw_\ge}$, the solution $u_\ge$ exists and it satisfies 
$w_m\leq u_\ge\leq w$ in $\Gw_\ge$. If $0<\ge'<\ge$, $u_{\ge'}\lfloor_{\prt\Gw_\ge}\geq u_\ge\lfloor_{\prt\Gw_\ge}=v_m$. Hence $u_{\ge'}\geq {\ge}$ in
$\Gw_\ge$. As in the proof of \rprop{keylem} the sequence $\{u_\ge\}$ is relatively compact in the $C^{1,\theta}_{loc}(\overline\Gw\setminus\{0\})$-topology, which ends the proof.\smallskip

\noindent In case (i), for $n>0$ set $K_n=\sup\left\{w(z): z\in\prt B^+_n\setminus\{0\}\right\}$ where, we recall it,  $B^+_n=B_n\cap \BBR^{_N}_+$. The function 
$w_{K_n}=(w-K_n)_+$ is a subsolution of $(\ref{K1})$ in $B^+_n$ which vanishes on $\prt B^+_n\setminus\{0\}$ and is smaller than $w$.
For $0<\ge<n$ we denote by $u_{\ge,n}$ the unique function satisfying
\bel{K7+1}\BA {lll}
-\Gd u_{\ge,n}+u_{\ge,n}^p-M|\nabla u_{\ge,n}e|^q=0&\quad\text{in }\,\Gg_{\ge,n}:=B^+_n\setminus\overline B^+_\ge \\
\phantom{-\Gd +u_{\ge,n}^p-M|\nabla u_{\ge,n}e|^q}
u_{\ge,n}=w_m&\quad\text{on }\,\prt B_\ge\cap\BBR^{_N}_+\\
\phantom{-\Gd +u_{\ge,n}^p-M|\nabla u_{\ge,n}e|^q}
u_{\ge,n}=0&\quad\text{on }\, (\prt B^+_n\cap\BBR^{_N}_+)\cup(\overline \Gg_{\ge,n}\cap\prt \BBR^{_N}_+).
\EA\ee
For $\ge'\leq \ge<n\leq n'$ there holds $w_{K_n}\leq u_{\ge,n}\leq u_{\ge',n'}\leq w$ in $\Gg_{\ge,n}$. Letting $n\to\infty$ and $\ge\to 0$ there exists a subsequence still denoted by $\{u_{\ge,n}\}$ which converges to a solution of $u$ of $(\ref{K1})$ in $\BBR^{_N}_+$ vanishing on $\prt\BBR^{_N}_+\setminus\{0\}$ and satisfying $(\ref{K6-})$.
\qeda\medskip

\nind\Remark The assumption that $w\lfloor_{\prt\Gw\setminus\{0\}}$ is bounded is restrictive. For example if $w(t\gs)=t^{-a}\gw(\gs)$ the flatness assumption means that $\dist(x,\BBR^{_N}_+)=O(|x|^{\ga+1})$ for all $x\in \prt\Gw$ near $0$ (remember that $T_{\prt\Gw}(0)=\prt\BBR^{_N}_+$). This assumption is always satisfied if $p\geq 3$ since $\ga\leq 1$, and it can be avoided if there exists a subsolution.

\bprop {keylem2+} Let $p>1$, $1\leq q\leq 2$ and $M\in\BBR$. Assume that the equation $(\ref{K1})$
admits a positive supersolution $w$ in $\BBR^{_N}_+$ belonging to  $C(\overline{\BBR^{_N}_+}\setminus\{0\})$ satisfying $(\ref{K5})$.
Assume also \smallskip

\nind (i) either there exists a positive subsolution 
$Z\in C(\overline{\BBR^{_N}_+}\setminus\{0\})$ vanishing on $\prt\BBR^{_N}_+\setminus\{0\}$, smaller than $w$ and satisfying $(\ref{K5})$, \smallskip

\nind (ii) or $\Gw\subset\BBR^{_N}_+$ is a bounded smooth domain such that $0\in\prt\Gw$ and there exists a positive subsolution 
$Z\in C(\overline{\Gw}\setminus\{0\})$, vanishing on $\prt\Gw\setminus\{0\}$ such that $Z\leq w\lfloor_{\Gw}$ and satisfying $(\ref{K5})$.\smallskip

\nind Then there exists a positive function $u$ satisfying $(\ref{K1})$ in $\BBR^{_N}_+$ (resp. $\Gw$), vanishing on $\prt\BBR^{_N}_+\setminus\{0\}$ (resp. $\prt\Gw\setminus\{0\}$) and such that 
\bel{K6a}\displaystyle
Z(x)\leq u(x)\leq w(x)\quad\text{for all }\,x\in \BBR^{_N}+  \;\,(\text{resp. }x\in \Gw).
\ee
\es

\nind {\it Example}. If $1<p<\frac{N+1}{N-1}$ it is proved in \cite{GmVe} that if $\Gw\subset\BBR^{_N}_+$ is a smooth bounded domain such that $0\in\prt\Gw$, there exists a nonnegative function $Z_\infty\in C(\overline\Gw\setminus\{0\})\cap C^2(\Gw)$ satisfying 
the equation 
\bel{K8}\BA {lll}\displaystyle
-\Gd Z+Z^p=0\qquad\text{in }\;\Gw\\
\phantom{-\Gd +Z^p}
Z=0\qquad\text{on }\;\prt\Gw\setminus\{0\},
\EA\ee
and such that $t^{\frac{2}{p-1}}Z_\infty(t\gs)\to \psi(\gs)$ uniformly on compact sets $K\subset S^{_{N-1}}_+$ as $t\to 0$ where $\psi$ is the unique a positive solution of 
\bel{K10}\BA {lll}\displaystyle
-\Gd'\psi+\ga\left(N-2-\ga\right)\psi+\psi^p=0\qquad&\text{in }\; S^{_{N-1}}_+\\
\phantom{-\Gd'-\ga\left(\ga+2-N\right)\psi+\psi^p}
\psi=0\qquad&\text{on }\; \prt S^{_{N-1}}_+.
\EA\ee
Furthermore, for any $k>0$ there exists a nonnegative function $Z_k\in C(\overline\Gw\setminus\{0\})\cap C^2(\Gw)$ satisfying $(\ref{K8})$ and such that 
$t^{N-1}Z_k(t\gs)\to k\phi_1(\gs)$ where $\phi_1$ has been introduced in \rth{soupscrit}, uniformly on compact subsets of $S^{_{N-1}}_+$. Furthermore $Z_k\uparrow Z_\infty$ when $k\to\infty$.  If the equation $(\ref{K1})$
admits a positive supersolution $w$ in $\BBR^{_N}_+$ belonging to  $C(\overline{\BBR^{_N}_+}\setminus\{0\})$ and such that 
$Z_k\leq w$ in  $\Gw$ for some $0<k\leq\infty$, then there exists a positive function $u$ satisfying $(\ref{K1})$ in $\Gw$, vanishing on $\prt\Gw\setminus\{0\}$ and such that 
\bel{K11}\displaystyle
Z_k(x)\leq u(x)\leq w(x)\quad\text{for all }\,x\in \Gw.
\ee
The same result holds if $\Gw$ is replaced by $\BBR^{_N}_+$.


\subsection{Existence or non-existence of separable solutions }

Since any large enough constant is a supersolution of $(\ref{Z2})$, it follows by \rth{BMP} that if there exists a nonnegative subsolution $z\in W^{1,\infty}_0(S^{_{N-1}}_+)$, there exists a solution in between. 
\subsubsection{Proof of \rth{exist}}

 We recall that $\gf_1$ is the first eigenfunction of $-\Gd'$ in $W^{1,2}_0(S_+^{_{N-1}})$ with corresponding eigenvalue $\gl_1=N-1$. 
 Put
$$H(\gw)=-\Gd' \gw+\ga(N-2-\ga)\gw+|\gf|^{p-1}\gw-M\left(\ga^2\gw^2+|\nabla' \gw|^2\right)^{\frac{p}{p+1}},
$$
then
$$H(\gf_1)=\left(N-1+\ga(N-2-\ga)\right)\gf_1+\gf_1^p-M\left(\ga^2\gf_1^2+|\nabla' \gf_1|^2\right)^{\frac{p}{p+1}}.
$$
If $\gf_1$ is small enough, there holds $\gf_1^p-M\left(\ga^2\gf_1^2+|\nabla' \gf_1|^2\right)^{\frac{p}{p+1}}<0$, hence 
$\gf_1$ is a subsolution.  However the condition $N-1+\ga(N-2-\ga)\leq 0$ is too stringent. We can use the fact that, up to a good choice of coordinates, $\gf_1=\gf_1(\gs)=\cos\gs$ with $\gs\in [0,\frac\gp 2]$. Furthermore the statement 
"$\gf_1$ is small enough" can be achieved by $\gf_1=\gd\cos\gs$ with $\gd>0$ small enough. Then 
$$\BA {lll}\gd^{-1}H(\gd^{\frac{p+1}{p-1}}\cos\gs)\\[2mm]
\phantom{----}=\left(N-1+\ga(N-2-\ga)\right)\cos\gs+\gd^{p+1}\cos^p\gs-M\gd(\ga^2\cos^2\gs+\sin^2\gs)^{\frac{p}{p+1}}.
\EA$$
The problem is to find $\gd>0$ such that for all $\gs\in [0,\frac\gp 2]$ we have $H(\gd^{\frac{p+1}{p-1}}\cos\gs)\leq 0$. 
Put $Z=\cos\gs$ and $\gd^{-1}H(\gd^{\frac{p+1}{p-1}}\cos\gs)=\gd^{-1}H(\gd^{\frac{p+1}{p-1}}Z)=K_\gd(Z)$, then
$$
K_\gd(Z)=\left(N-1+\ga(N-2-\ga)\right)Z+\gd^{p+1}Z^p-M\gd((\ga^2-1)Z^2+1)^{\frac{p}{p+1}},
$$
where $0\leq Z\leq 1$.
We use the fact that 
$$\ga^2\cos^2\gs+\sin^2\gs\geq \min\{\ga^2,1\}(\cos^2\gs+\sin^2\gs):=\gk^2>0,$$
hence
$$\gd^{-1}H(\gd^{\frac{p+1}{p-1}}\cos\gs)\leq \left(N-1+\ga(N-2-\ga)\right)\cos\gs+\gd^{p+1}\cos^p\gs-M\gd\gk^{\frac{2p}{p+1}}.
$$
Then 
\bel {M}K_\gd(Z)\leq \tilde K_\gd(Z):=\left(N-1+\ga(N-2-\ga)\right)Z+\gd^{p+1}Z^p-M\gd\gk^{\frac{2p}{p+1}},
\ee
and
\bel {M0}\tilde K'_\gd(Z)=N-1+\ga(N-2-\ga)+p\gd^{p+1}Z^{p-1}.
\ee

\nind{ If  $N-1+\ga(N-2-\ga)\geq  0$, equivalently $p\geq  \frac{N+1}{N-1}$, then $\tilde K'_\gd\geq 0$ on $[0,1]$, hence 
$$\tilde K_\gd(Z)\leq \tilde K_\gd(1)=N-1+\ga(N-2-\ga)+\gd^{p+1}-M\gd\gk^{\frac{2p}{p+1}}.
$$
The function $\gd\mapsto \tilde K_\gd(1)$ achieves its minimum for $\gd=\gd_0:=\gb^{\frac{2}{p+1}}\left(\frac M{p+1}\right)^{\frac 1p}$ and 
$$\tilde K_{\gd_0}(1)=N-1+\ga(N-2-\ga)-p\gk^2\left(\frac M{p+1}\right)^{\frac {p+1}p}.
$$
Therefore, when $p\geq  \frac{N+1}{N-1}$, $K_{\gd_0}\leq 0$ on $[0,1]$ if 
\bel{M1}  \BA {lll}
\left(\frac M{p+1}\right)^{\frac {p+1}p}\geq\left(\frac {M_{_{N,p}}}{p+1}\right)^{\frac {p+1}p}:=\myfrac{N-1+\ga(N-2-\ga)}{p\min\{1,\ga^2\}}\\[4mm]
\phantom{\left(\frac M{p+1}\right)^{\frac {p+1}p}\geq\left(\frac {M_{_{N,p}}}{p+1}\right)^{\frac {p+1}p}:}=
\myfrac{(p+1)\left(p(N-1)-(N+1)\right)}{p\min\{(p-1)^2,4\}}.
\EA\ee

\smallskip

\nind{\it  If $N-1+\ga(N-2-\ga)\leq 0$}, equivalently $p\leq \frac{N+1}{N-1}$, it is clear from (\ref{M}) that $\tilde K_\gd(Z)\leq 0$
for any $Z\in [0,1]$ as soon as $\gd\leq \gk^{\frac 1{p+1}}M^{\frac{1}{p}}$. \smallskip

\nind{\it Improvement in the case $\ga>1$, equivalently $1<p<3$ }. We set
$$F(Z)=\myfrac{(\ga^2-1)Z^2+1}{Z^{\frac{p+1}{p}}}.
$$ 
Then 
$$\myfrac{F'(Z)}{F(Z)}=\myfrac{(p-1)(\ga^2-1)Z^2-(p+1)}{p((\ga^2-1)Z^2+1)Z}.
$$
Since 
\bel{M1a} K_\gd(Z)\leq 0\Longleftrightarrow \left(N-1+\ga(N-2-\ga)\right)+\gd^{p+1}Z^{p-1}\leq M\gd F^{\frac{p}{p+1}}(Z)\ee
for all $Z\in (0,1]$, it is sufficient to prove 
\bel{M1b}\displaystyle
 \left(N-1+\ga(N-2-\ga)\right)+\gd^{p+1}\leq M\gd\min_{Z\in (0,1]}F^{\frac{p}{p+1}}(Z)
\ee
The function $F$ is minimal on $(0,1]$ at $Z=Z_0=\sqrt{\ga^2-1}$ (remember that $\ga=\frac{2}{p-1}$) and 
$F(Z_0)=(\ga+2)(\ga-1)^{\frac{\ga+1}{\ga+2}}$. \\
If $Z_0\leq 1$, equivalently $\ga\geq 2$, inequality $(\ref{M1b})$ is satisfied if one find $\gd$ such that 
$$ \left(N-1+\ga(N-2-\ga)\right)+\gd^{p+1}\leq M\gd F^{\frac{p}{p+1}}(Z_0),
$$
and a sufficient condition is 
\bel{M1c}\displaystyle
p\left(\myfrac{M}{p+1}\right)^{\frac{p+1}{p}}\geq p\left(\myfrac{M_{_{N,p}}}{p+1}\right)^{\frac{p+1}{p}}:=\myfrac{N-1+\ga(N-2-\ga)}{F(Z_0)}
\ee
If $Z_0> 1$, equivalently $1<\ga< 2$, the minimum of $F$ on $(0,1]$ is achieved at $Z=1$ with value $F(1)=\ga^2$, hence a sufficient condition is 
$$ \left(N-1+\ga(N-2-\ga)\right)+\gd^{p+1}\leq M\gd\ga^{\frac{2p}{p+1}}),
$$
and we obtain the desired inequality as soon as 
\bel{M1d}\displaystyle
p\left(\myfrac{M}{p+1}\right)^{\frac{p+1}{p}}\geq p\left(\myfrac{M_{_{N,p}}}{p+1}\right)^{\frac{p+1}{p}}:=\myfrac{N-1+\ga(N-2-\ga)}{\ga^2}.
\ee
This ends the proof.\qeda\medskip

\noindent\Remark Introducing $m^{**}$ defined in $(\ref{Z6})$, inequality $(\ref{M1})$ takes the form
\bel{M2}  
M\geq \left(\myfrac{2(p+1)}{\min\{(p-1)^2,4\}}\right)^{\frac{p}{p+1}}m^{**},
\ee
in the general case and a more complicated expression in the case $\ga>1$.
\subsubsection{Non-existence}\smallskip

\bth{non-ex} Let $p>\frac{N+1}{N-1}$ and $M\leq m^{**}$, defined by $(\ref{Z6})$. Then equation $(\ref{Z2})$ admits no positive solution.\es

\noindent\Proof If $\gw$ is a positive solution of $(\ref{Z2})$ the function $\eta$ defined by $\gw=\eta^b$ for some $b>0$ satisfies
$$\BA {lll}-\Gd'\eta+(1-b)\myfrac{|\nabla '\eta|^2}{\eta}+\myfrac{\ga(N-2-\ga)}{b}\eta+\myfrac{1}{b}\eta^{1+(p-1)b}\\[4mm]
\phantom{---------------}
-\myfrac{M\eta^{\frac{(b-1)(p-1)}{p+1}}}{b}\left(\ga^2\eta^2+b^2|\nabla '\eta|^2\right)^{\frac {p}{p+1}}=0.
\EA$$
Since for any $\ge>0$ we have by H\"older's inequality,
$$\BA {lll}
\myint{S^{N-1}_+}{}\eta^{1+\frac{(b-1)(p-1)}{p+1}}\left(\ga^2\eta^2+b^2|\nabla '\eta|^2\right)^{\frac {p}{p+1}}dS\\[4mm]
\phantom{------}\leq
\myfrac{\ge^{\frac{p+1}{p}}p}{p+1}\myint{S^{N-1}_+}{}(\ga^2\eta^2+b^2|\nabla'\eta|^2)+\myfrac{1}{(p+1)\ge^{p+1}}\myint{S^{N-1}_+}{}\eta^{2+(p-1)b}dS,
\EA$$
it follows that
\bel{Ne1}\BA {lll}
\left(2-b-M\myfrac{\ge^{\frac{p+1}{p}}pb}{p+1}\right)\myint{S^{N-1}_+}{}|\nabla'\eta|^2dS+\myfrac{\ga}{b}\left(N-2-\ga-M\myfrac{\ge^{\frac{p+1}{p}}\ga p}{p+1}\right)\myint{S^{N-1}_+}{}\eta^2dS\\
[4mm]\phantom{-------}
+\myfrac{1}{b}\left(1-\myfrac{M}{(p+1)\ge^{p+1}}\right)\myint{S^{N-1}_+}{}\eta^{2+(p-1)b}dS\leq 0.
\EA\ee
If $b\in (0,2)$, $\ge>0$ and $M>0$ are linked by the relation 
\bel{Ne2}
2-b-M\myfrac{\ge^{\frac{p+1}{p}}pb}{p+1}\geq 0\Longleftrightarrow M\ge^{\frac{p+1}{p}}\leq \myfrac{(2-b)(p+1)}{bp},
\ee
inequality $(\ref{Ne1})$ turns into 
\bel{Ne3}\BA {lll}
\left((2-b)(N-1)+\myfrac{\ga(N-2-\ga)}{b}-\myfrac{M\ge^{\frac{p+1}{p}}p}{p+1}\left((N-1)b+\myfrac{\ga^2}{b}\right)\right)\myint{S^{N-1}_+}{}\eta^2dS\\
[4mm]\phantom{--------------}
+\myfrac{1}{b}\left(1-\myfrac{M}{(p+1)\ge^{p+1}}\right)\myint{S^{N-1}_+}{}\eta^{2+(p-1)b}dS\leq 0.
\EA\ee
Next we choose
\bel{Ne4}\ge^{p+1}=\myfrac{M}{p+1},
\ee
and we define the function $b\mapsto L(b)$ by 
\bel{Ne5}
L(b):=(2-b)(N-1)+\myfrac{\ga(N-2-\ga)}{b}-p\left(\myfrac{M}{p+1}\right)^{\frac{p+1}{p}}\left((N-1)b+\myfrac{\ga^2}{b}\right).
\ee
Because $N-1$ is the first eigenvalue of $-\Gd'$ in $W^{1,2}_0(S^{^{_{N-1}}})$, $(\ref{Ne3})$ combined with  $(\ref{Ne4})$ yields
\bel{Ne6}
L(b)\myint{S^{N-1}_+}{}\eta^2dS\leq 0.
\ee
Furthermore, if inequality $(\ref{Ne2})$ is strict, and since $\eta$ is not a first eigenfunction, inequality $(\ref{Ne6})$ is also strict. Then 
$L(b)\geq 0$ if 
\bel{Ne7}
p\left(\myfrac{M}{p+1}\right)^{\frac{p+1}{p}}\leq f(b):=\myfrac{b(2-b)(N-1)+\ga(N-2-\ga)}{(N-1)b^2+\ga^2}.
\ee
Now
$$f'(b)=\myfrac{-2(N-1)^2}{\left((N-1)b^2+\ga^2\right)^2}\left(b+\ga\right)\left(b-\myfrac{\ga}{N-1}\right).
$$
Notice that 
$$\myfrac{\ga}{N-1}\leq 1\Longleftrightarrow p\geq\frac{N+1}{N-1}.
$$
If $1<p\leq \frac{N+1}{N-1}$, then $f'\geq 0$  and in such a case the maximum of $f$ over $(0,1]$ is achieved at 
$b=1$ and for such a value $L(b)\leq 0$. \smallskip

\noindent If $p> \frac{N+1}{N-1}$, then $f$ is increasing on $[0,\frac{\ga}{N-1})$ and decreasing on $(\frac{\ga}{N-1},1]$, hence the maximum is achieved at $b=\frac{\ga}{N-1}$, which gives 
\bel{Ne8}f\left(\frac{\ga}{N-1}\right)=\myfrac{N-1-\ga}{\ga}=\myfrac{(N-1)p-(N+1)}{2}. 
\ee
Therefore there exists no solution if $p\geq \myfrac{N+1}{N-1}$ and
\bel{Z4b}
\left(\myfrac{M}{p+1}\right)^{\frac{p+1}{p}}\leq\left(\myfrac{m^{**}}{p+1}\right)^{\frac{p+1}{p}}:= \myfrac{(N-1)p-(N+1)}{2p}.
\ee
\qeda
\medskip

\noindent\Remark Using \rth{remov} we can prove the previous result in the case $M< m^{**}$ . Indeed, if $\gw$ is a positive solution of $(\ref{Z2})$, $u_\gw(r,.)=r^{-\frac{2}{p-1}}\gw(.)$ is a positive solution of $(\ref{Z1})$ in $\BBR^{_N}_+$ vanishing on
$\prt\BBR^{_N}_+\setminus\{0\}$. Let $\Gw\subset\BBR^{_N}_+$ be any smooth domain such that $0\in\prt\Gw$ and  $\prt\Gw$ is flat near $0$. Then 
$u_\gw\leq K$ on $\prt\Gw$ for some $K>0$. Put $v=(u_\gw- K)_+$, then it is a nonnegative subsolution of $(\ref{Z2})$. For any $\ge>0$ small enough there exists a solution 
$u_\ge$ of 
\bel{M31}\BA {lll}
-\Gd u+u^p-M|\nabla u|^{\frac{2p}{p+1}}=0\qquad&\text{in }\;\Gw_\ge:=\Gw\cap \overline B_\ge^c\\
\phantom{-\Gd +u^p-M|\nabla |^{\frac{2p}{p+1}}}
u=v&\text{on }\;\prt B_\ge\cap \Gw\\
\phantom{-\Gd +u^p-M|\nabla |^{\frac{2p}{p+1}}}
u=0&\text{on }\;B^c_\ge\cap \prt \Gw.
\EA\ee
Then $v\leq u_\ge\leq u_\gw$. Furthermore, for $0<\ge'<\ge$, $u_\ge\leq u_{\ge'}$ in $\Gw_\ge$. Hence $\{u_\ge\}$ converges, when $\ge \to 0$ to a solution $u_0$ of $(\ref{NZ2})$, which satisfies $v\leq u_0\leq u_\gw$ and therefore vanishes on $\prt\Gw\setminus\{0\}$, which is a contradiction.

\section{Solutions with an isolated boundary singularity}

\subsection{Uniqueness of singular solutions}
In this section we use scaling transormations to prove uniqueness of singular solutions.
\bth{uniq} Assume $N\geq 2$, $p>1$, $1<q<2$ and $M>0$. Let $a$ such that
\bel{uniq1}\BA{lll}
(i) \qquad\qquad  0\leq a<\gb\quad&\text{if }\,q\leq\frac{2p}{p+1}\qquad \qquad\qquad \qquad\qquad \qquad   \\[2mm]
(ii) \qquad\qquad  \gb<a\quad&\text{if }\,q>\frac{2p}{p+1}.\qquad \qquad \qquad \qquad \qquad \qquad 
\EA
\ee
Let $\gth\in C^1(\overline{S^{N-1}_+})$ be a nonnegative function, positive in  $S^{N-1}_+$, vanishing on $\prt S^{N-1}_+$, and $\tilde a$ be a real number. Then there exists at most one nonnegative solution of $(\ref{Z8})$ satisfying 
 \bel{uniq2}\BA{lll}
u(r,\gs)=r^{-a}|\ln r|^{\tilde a}\gth(\gs)(1+o(1))\qquad\text{as }\,r\to 0.
\EA
\ee
\es
\Proof The proof is an adaptation to the configuration where the  singularity lies on the boundary of \cite[Theorem 1.15]{BVGHV4}. If $u$ is a solution of $(\ref{Z1})$ in $\BBR^N_+$, $\ell>0$ and $b>0$, define $u_{\ell,b}$ by 
$$u_{\ell,b}(x)=\ell^bu(\ell). 
$$
Then 
$$-\Gd u_{\ell,b}+\ell^{2-b(p-1)}u_{\ell,b}^p-M\ell^{2-q-b(q-1)}|\nabla u_{\ell,b}|^q=0\quad\text{in }\BBR^N_+.
$$
If $\ell>1$, $u_{\ell,b}$ is a supersolution of $(\ref{Z1})$ in $\BBR^N_+$ if and only if 
$$\ga\leq b\leq\gb.
$$
These conditions are compatible if and only if $1<q\leq\frac{2p}{p+1}$. 
Then we take $b\in(a,\infty)\cap[\ga,\gb]$, then 
$$u_{\ell,b}(r,\gs)=\ell^{b-a}r^{-a}|\ln r|^{\tilde a}\gth(\gs)(1+o(1))\quad\text{as }r\to 0.
$$
By $(\ref{1Q1})$ all $u(x)$ tends to $0$ when $|x|\to\infty$. Hence, for any $\ge>0$ the super solution $u_{\ell,b}+\ge$ which is larger than $u$ for $|x|$ small enough and large enough is larger than another solution $\tilde u$ in $\BBR^N_+$. Letting $\ell\downarrow 1$ and $\ge\to 0$ yields $u\geq\tilde u$. In the same way  $\tilde u\geq u$.

\nind If $\ell<1$, $u_{\ell,b}$ is a supersolution of $(\ref{Z1})$ in $\frac{1}{\ell}G$ if and only if 
$$\gb\leq b\leq\ga,
$$
and these conditions are compatible if and only if $\frac{2p}{p+1}\leq q<2$. If $\ga>\gb$ we choose  $b\in (0,a)\cap[\gb,\ga]$ and we conclude as in the first case. \qeda\medskip

\nind\Remark In the case $a=\gb$ a more precise expansion of the singular solution $u$ at $x=0$ yields uniqueness as it is proved in \cite{BVGHV4} in the case of an internal singularity. Since the proof of the next result is based upon a easy adaptation of the ideas in \cite[Theorem 4.4]{BVGHV4}, we omit it. 
\bth{uniq3} Assume $N\geq 2$, $p>1$, $1<q\leq\frac{2p}{p+1}$, $M>0$ and $a\in [0,\gb]$. Assume $\gth$ and $\tilde\gth$ are $C^1(\overline {S^{N-1}_+})$ functions positive in $S^{N-1}_+$ and vanishing on $\prt S^{N-1}_+$ and  ${\tilde a}$ is a real smaller than $a$. Then there exists at most one nonnegative solution of $(\ref{Z8})$ satisfying 
 \bel{uniq4}\BA{lll}
u(r,\gs)=r^{-a}\gth(\gs)+r^{-{\tilde a}}\tilde\gth(\gs)(1+o(1))\qquad\text{as }\,r\to 0.
\EA
\ee
\es

When problem $(\ref{Z8})$ is replaced by $(\ref{NZ1})$ the scaling method becomes much more delicate to apply. However we give below an easy extension when $\prt\Gw$ is flat near  $x=0$.
\bth{uniq5} Assume $N\geq 2$, $p>1$, $1<q<2$, $M>0$ and $\Gw$ is a bounded smooth domain such that 
$0\in\prt\Gw$ and there exists $\gd>0$ such that $\prt\Gw\cap B_\gd=T_{\prt\Gw}(0)\cap B_\gd$. Let $a$ such that
\bel{uniq6}\BA{lll}
(i) \qquad\qquad  0\leq a<\gb\quad&\text{if }\,q\leq\frac{2p}{p+1}\qquad \qquad\qquad \qquad\qquad \qquad   \\[2mm]
(ii) \qquad\qquad  \gb<a\quad&\text{if }\,q>\frac{2p}{p+1}.\qquad \qquad \qquad \qquad \qquad \qquad 
\EA
\ee
Let $\gth\in C^1(\overline{S^{N-1}_+})$ be a nonnegative function, positive in  $S^{N-1}_+$, vanishing on $\prt S^{N-1}_+$, and $\tilde a$ be a real number. If $\Gw$ is starshaped with respect to $0$, then there exists at most one nonnegative solution of $(\ref{Z8})$ satisfying 
 \bel{uniq7}\BA{lll}
u(r,\gs)=r^{-a}|\ln r|^{\tilde a}\gth(\gs)(1+o(1))\qquad\text{as }\,r\to 0.
\EA
\ee
\es
\Proof We use the same change of scale as in \rth{uniq}. In case (i) with $\ell>1$ and $b\in (a,\infty)\cap[\ga,\gb]$ , $u_{\ell,b}$ is a supersolution in $\Gw_\ell=\frac{1}{\ell}\Gw\subset \Gw$ and $\prt\Gw_\ell\cap B_{\gd/\ell}=T_{\prt\Gw}(0)\cap B_{\gd/\ell}$. If $\tilde u$ 
is another solution, $\tilde u\lfloor_{\prt\Gw_\ell}=h_\ell$ and $h_\ell\to 0$ uniformly as $\ell\downarrow 1$ since $u\in C^1(\overline\Gw\cap B^c_\ge)$ for any $\ge>0$. The function 
$\displaystyle u_{\ell,b}+\max_{\prt\Gw_\ell} h_\ell$ is a supersolution of $\ref{Z1}$ in $\Gw_\ell$ larger than $\tilde u$ on $\prt\Gw_\ell\setminus\{0\}$ and near $x=0$, hence it is larger than $\tilde u$ in $\Gw_\ell$. Letting $\ell\downarrow 1$ yields $u\geq \tilde u$. \\
In case (ii), for $\ell<1$, $\Gw\subset \Gw_\ell$ and with $b\in (0,a)\cap [\gb,\ga]$, $u_{\ell,b}$ is a supersolution in $\Gw$ larger that $\tilde u$ on $\prt\Gw\setminus\{0\}$ and near  $x=0$ it is larger than $\tilde u$ in $\Gw$. We conclude as in case (i).\phantom{---------------}\qeda
\subsection{Construction of fundamental solutions}
Let $\Gw$ be either $\BBR^{_N}_+$ or a bounded domain with $0\in\prt\Gw$. A function $u$ satisfying $(\ref{NZ1})$
is a {\it fundamental solution} if it has  a singularity of potential type, that is
\bel{P0-}\BA {lll}\displaystyle
\lim_{x\to 0}\myfrac{|x|^Nu(x)}{\gr(x)}=c_{_N}k,
\EA\ee
for some $k>0$. The function $u$ can also be looked for as a solution of 
\bel{P1}\BA {lll}
-\Gd u+u^{p}-M|\nabla u|^{q}=0\qquad&\text{in }\Gw\\
\phantom{-\Gd +u^{p}-M|\nabla u|^{q}}
u=k\gd_0&\text{in }\prt\Gw,
\EA\ee
in the sense that $u\in L^p_\gr(\Gw\cap B_r)$, $\nabla u\in L^{q}_{\gr,loc}(\Gw\cap B_r,)$ for any $r>0$, and for any $\gz\in C^{1}_c(\overline\Gw)\cap W^{2,\infty}(\Gw)$ there holds
\bel{P2}\BA {lll}
\myint{\Gw}{}\left(-u\Gd\gz+u^p\gz-M|\nabla u|^{q}\gz\right) dx=-k\myfrac{\prt\gz}{\prt{\bf n}}(0).
\EA\ee


We first consider the problem in $\BBR^{_N}_+$. \medskip
 
\nind{\it Proof of \rth{weak-sing}}. The scheme of the proof is surprising since we first show that, in the case $q=\frac{2p}{p+1}$, there exists $M_1>0$ such that for any $k>0$ and any $0<M<M_1$ there exists a solution. Using this result we prove that if $1<q<\frac{2p}{p+1}$, then for any $M>0$ and $k>0$ there exists a solution. Then  we return to the case $q=\frac{2p}{p+1}$ and using the result in the previous case, we prove that when $q=\frac{2p}{p+1}$ we can get rid of the restriction on $M>0$ and $k>0$ for the existence of solutions. \smallskip

\noindent {\it I- The case $q=\frac{2p}{p+1}$ and $M$ upper bounded.} \smallskip
 
\noindent For $\ell>0$ the transformation $T_\ell$  defined by 
\bel{T1}
T_\ell[u(x)]=\ell^{\frac{2}{p-1}}u(\ell x),
\ee
leaves the operator $\CL_{_{\frac{2p}{p+1},M}}$ invariant. We can therefore {\it write}
$$T_\ell[u_k]=u_{k\ell^{\frac{2}{p-1}+1-N}}, $$
in the sense that if $u_k$ satisfies $(\ref{P0-})$ then $T_\ell[u_k]$ satisfies the same limit with $k$ replaced by $k\ell^{\frac{2}{p-1}+1-N}$. However this identity to hold needs some uniqueness for the solutions under consideration satisfying 
$(\ref{P0-})$. This is achieved if  $u_k$ is the minimal solution satisfying $(\ref{P0-})$ 
in which case $T_\ell[u_k]$ is the minimal solution satisfying $(\ref{P0-})$  with $k$ replaced by $k\ell^{\frac{2}{p-1}+1-N}$. Therefore if there exists a solution to $(\ref{NZ1})$ in $\BBR^{_N}_+$, vanishing on $\prt\BBR^{_N}_+\setminus\{0\}$ satisfying $(\ref{P0})$ for some $k>0$, then there exists such a solution for any $k>0$.
\smallskip

\nind {\it Step 1- Construction of a subsolution}. For $k>0$ we denote by  $v_k$ the solution of 
\bel{P7e}\BA {lll}
-\Gd v+v^{p}=0\qquad&\text{in }\BBR^{_N}_+\\
\phantom{-\Gd +v^{p}}
v=k\gd_0&\text{on }\prt\BBR^{_N}_+\setminus\{0\}.
\EA\ee
Such a solution exists thanks to \cite{GmVe} if $\BBR^{_N}_+$ is replaced by a bounded domain $\Gw$. If case  of a half-space the problem is first solved in $B_n^+$ and by letting $n\to\infty$, we obtain the solution in  $\BBR^{_N}_+$. 
Clearly $v_k$  is a subsolution of problem $(\ref{NZ1})$, and it satisfies
\bel{P0}\BA {lll}\displaystyle
\lim_{x\to 0}\myfrac{u(x)}{P_{_N}(x)}=k,
\EA\ee
for some $c'_{_N}>0$, where $P_{_N}(x)=c_{_N}\frac{x_{{_N}}}{|x|^N}$ is the Poisson kernel in $\BBR^{_N}_+$.\smallskip

\nind {\it Step 2- Construction of a supersolution}.  
It is known that
  \bel{P7}|\nabla P_{_N}(x)|^2=|x|^{-2N}c^2(x),
\ee
where $c(.)$ is smooth and verifies
$$0<\tilde c_1\leq c(x)\leq \tilde c_2\quad\text{for some }\tilde c_1,\tilde c_2>0.
$$
We construct $w_k$ in $\BBR^{_N}_+$ of the form
  \bel{Q1e}\BA {lll}
w_k=kP_{_N}+w,
\EA
\ee
where $w$ satisfies 
  \bel{Q2e}\BA {lll}
-\Gd w+w^p=a\gg_2|x|^{-\frac{2Np}{p+1}}\quad&\text{in }\,\BBR^{_N}_+\\
\phantom{-\Gd +w^p}
w=0&\text{on }\,\prt\BBR^{_N}_+,
\EA
\ee
for some $a>0$ to be chosen later on. Then
$$\BA {lll}\CL_{_{\frac{2p}{p+1},M}} w_k=-\Gd w+(kP_{_N}+w)^p -M\left(|k\nabla P_{_N}+\nabla w|^{2}\right)^{\frac{p}{p+1}}\\
\phantom{\CL_{_{\frac{2p}{p+1},M}} w_k}
=(kP_{_N}+w)^p-w^p+a\gg_2|x|^{-\frac{2Np}{p+1}}-M\left(|k\nabla P_{_N}+\nabla w|^{2}\right)^{\frac{p}{p+1}}\\
\phantom{\CL_{_{\frac{2p}{p+1},M}} w_k}
\geq pkP_{_N}w^{p-1}+a\gg_2|x|^{-\frac{2Np}{p+1}}-2M\left(k^{\frac{2p}{p+1}}\gg_2^{\frac{2p}{p+1}}|x|^{-\frac{2Np}{p+1}}+|\nabla w|^{\frac{2p}{p+1}}\right).
\EA$$
Now it is easy to check using Osserman's type construction as in \cite[Lemma 2.1]{Ve84} and scaling techniques that 
$$w(x)\leq \gg_3 \min\left\{a^{\frac1p}|x|^{-\frac{2N}{p+1}},a|x|^{2(1-\frac{Np}{p+1})}\right\},
$$
and 
$$\BA {c}|\nabla w(x)|\leq \gg_4 \min\left\{a^{\frac1p}|x|^{-\frac{2N}{p+1}-1},a|x|^{1-\frac{2Np}{p+1}}\right\}\\
\Longrightarrow\\
|\nabla w(x)|^{\frac{2p}{p+1}}\leq \gg_5\min\left\{a^{\frac2{p+1}}|x|^{-\frac{2p(2N+p+1))}{(p+1)^2}},a^{\frac{2p}{p+1}}|x|^{\frac{2p(p+1-2Np)}{(p+1)^2}}\right\}.
\EA$$
Therefore, if we put
$$\gt=\myfrac{p^2-1}{2p(N+1-p(N-1))},
$$
then $\gt>0$ since $N+1>p(N-1)$ and 
\bel{Q3}\BA {lll}|x|^{\frac{2Np}{p+1}}\CL_{_{\frac{2p}{p+1},M}} w_k\geq 
\gg_2\left(\!a-2Mk^{\frac{2p}{p+1}}\gg_2^{\frac{p-1}{p+1}}\right)\!-2M\gg_5a^{\frac{2p}{p+1}}|x|^{\frac{2p(N+1-p(N-1))}{(p+1)^2}}\\
[2mm]\phantom{|x|^{\frac{2Np}{p+1}}\CL_{_{\frac{2p}{p+1},M}} w_k}
\geq 
\gg_2\left(\!a-2Mk^{\frac{2p}{p+1}}\gg_2^{\frac{p-1}{p+1}}\right)-2M\gg_5a\quad \text{in }\;B^+_{a^{\gt}},
\EA\ee
and similarly, 
\bel{Q4}\BA {lll}|x|^{\frac{2Np}{p+1}}\CL_{_{\frac{2p}{p+1},M}} w_k\geq 
\gg_2\left(a-2Mk^{\frac{2p}{p+1}}\gg_2^{\frac{p-1}{p+1}}\right)-2M\gg_5a\;\text{in }\;(B^+_{a^{\gt}})^c.
\EA\ee
Replacing $\gt$ by its value, we obtain a very simple expression from $(\ref{Q3})$ and $(\ref{Q4})$, valid both in 
$B^+_{a^{\gt}}$ and $(B^+_{a^{\gt}})^c$, namely
\bel{Q5+}\BA {lll}|x|^{\frac{2Np}{p+1}}\CL_{_{\frac{2p}{p+1},M}} w_k\geq 
\gg_2\left(a-2Mk^{\frac{2p}{p+1}}\gg_2^{\frac{p-1}{p+1}}\right)-2M\gg_5a\;\,\text{ in }\BBR^{_N}_+.
\EA
\ee
When
\bel{Q6e}\BA {lll}
M<M_1:=\myfrac{\gg_2}{2\gg_5},
\EA
\ee
then for fixed $k$, if we take 
$$a>\myfrac{2M_1\gg_2^{\frac{2p}{p+1}}k^{\frac{2p}{p+1}}}{\gg_2-2M\gg_5},
$$
we infer that the right-hand side of $(\ref{Q5+})$ is nonnegative, hence $w_k$ is a supersolution. 
\smallskip

\nind {\it Step 3-Existence}. For $0<k\leq k_0$ $w_k$ is a supersolution which dominates the subsolution 
$v_k$. Hence, by  \cite[Theorem 1-4-6]{Vebook} there exists a solution $u_k$ to $(\ref{NZ1})$ in $\BBR^{_N}_+$, vanishing on $\prt\BBR^{_N}_+\setminus\{0\}$ and such that $v_k\leq u_k\leq w_k$. Since
$$\displaystyle\lim_{x\to 0}\myfrac{v_k(x)}{P_{_N}(x)}=\lim_{x\to 0}\myfrac{w_k(x)}{P_{_N}(x)}=k,
$$
it follows that $u_k$ inherits the same asymptotic behaviour. Since $k<k_0$ can be replaced by any $k>0$, the existence of a solution follows.\smallskip

\noindent{\it II- The case $1<q<\frac{2p}{p+1}$.} Assume $M<M_1$, $k>0$ and  $\tilde u_k$ is the minimal solution of  $(\ref{NZ1})$ in $\BBR^{_N}_+$ with $q=\frac{2p}{p+1}$, vanishing on 
$\prt\BBR^{_N}_+\setminus\{0\}$ and such that $(\ref{P0})$. Since $|\nabla \phi|^{\frac{2p}{}2p}\geq |\nabla \phi|^{q}-1$, there holds
$$-\Gd \tilde u_k+\tilde u_k^p+M-M|\nabla \tilde u_k|^{q}\geq 0.
$$ 
Hence $\tilde u^*_k=\tilde u_k+M^{\frac1p}$ is a supersolution $(\ref{NZ1})$ in $\BBR^{_N}_+$ and it dominates $v_k$ defined in $(\ref{P7e})$. By  \cite[Theorem 1-4-6]{Vebook} there exists a solution $u_k$ of $(\ref{NZ1})$, vanishing on $\prt\BBR^{_N}_+\setminus\{0\}$ and satisfying $(\ref{P0})$ under the following weaker form 
\bel{P0x}\BA {lll}\displaystyle
\lim_{t\to 0}\myfrac{u_k(tx)}{P_{_N}(tx)}=k\quad\text{uniformly on compact subsets of }\,\BBR^{_N}_+.
\EA\ee
Since $|x|^{N-1}u_k(x)$ is uniformly bounded and vanishes on $\prt\BBR^{_N}_+\setminus\{0\}$, it is bounded in the $C^1_{loc}(\overline{\BBR^{_N}_+})$-topology. 
Hence $(\ref{P0})$ holds. This proves the result when $M<M_1$. \smallskip

Next let $M>0$ arbitrary and $k>0$. In order to find a solution $u:=u_k$ to $(\ref{NZ1})$, we set $u(x)=\ell^{-\frac{2}{p-1}}U_\ell(\frac x{\ell})$. Then $\CL_{q,M}u=0$ is equivalent to 
$$\CL_{_{q,M_\ell}}U_\ell:=-\Gd U_\ell+U_\ell^p-M_\ell|\nabla U|^q=0\quad\text{with }\;M_\ell=M\ell^{\frac{2p-q(p+1)}{p-1}},
$$
and $(\ref{P0})$ is equivalent to
$$\displaystyle \lim_{x\to 0}\myfrac{U_\ell(x)}{P_{_N}(x)}=\ell^{\frac{2}{p-1}+1-N} k.
$$
Since $2p-q(p+1)>0$ it is enough to choose $\ell>0$ such that $M\ell^{\frac{2p-q(p+1)}{p-1}}<M_1$, and we end the proof using the result when
$M<M_1$. \smallskip

\noindent{\it III- The case $q=\frac{2p}{p+1}$ revisited.} Let $p<\tilde p<\frac{N+1}{N-1}$. Then $\frac{2p}{p+1}<\frac{2 \tilde p}{\tilde p+1}$. This implies that for any $M>0$ and $k>0$ there exists a positive solution $\tilde u_k$ to 
$$-\Gd \tilde u_k+\tilde u_k^{\tilde p}-M|\nabla \tilde u_k|^{\frac{2p}{p+1}}= 0\qquad\text{in }\;\BBR^{_N}_+,
$$ 
vanishing on $\prt\BBR^{_N}_+\setminus\{0\}$ and such that 
$$\displaystyle\lim_{x\to 0}\myfrac{\tilde u_k(x)}{P_{_N}(x)}=k. 
$$
Since $\tilde p>p$ we have $\tilde u_k^{\tilde p}>\tilde u_k^{ p}-1$ and therefore 
\bel{QXe}-\Gd \tilde u_k+\tilde u_k^{p}-M|\nabla \tilde u_k|^{\frac{2p}{p+1}}\geq 1> 0\qquad\text{in }\;\BBR^{_N}_+.
\ee
The function $\tilde v_k$ solution of
\bel{P7x}\BA {lll}
-\Gd v+v^{ p}=0\qquad&\text{in }\BBR^{_N}_+\\
\phantom{-\Gd +v^{p}}
v=k\gd_0&\text{on }\prt\BBR^{_N}_+\setminus\{0\},
\EA\ee
is a subsolution of $(\ref{QXe})$, hence the exists a solution  $u_k$ of  such that $\tilde v_k<u_k< \tilde u_k$ of $(\ref{NZ1})$ in $\BBR^{_N}_+$, vanishing on 
$\prt\BBR^{_N}_+\setminus\{0\}$ and such that $(\ref{P0-})$ holds. 
\smallskip

\nind {\it IV- The case $\frac{2p}{p+1}<q<\frac {1+N}N$}. 
We follow the ideas of Case I. We look for a supersolution $w_k$ of the form $(\ref{Q1e})$ where $w_k$ satisfies 
\bel{QYe}\BA {lll}
-\Gd w+w^p=a\gg_2|x|^{-Nq}\quad&\text{in }\,\BBR^{_N}_+\\
\phantom{-\Gd +w^p}
w=0&\text{on }\,\prt\BBR^{_N}_+,
\EA
\ee
for some $a>0$. Then 
$$\BA {lll}\CL_{_{q,M}} w_k=-\Gd w+(kP_{_N}+w)^p -M\left(|k\nabla P_{_N}+\nabla w|^{2}\right)^{\frac{q}{2}}\\[1mm]
\phantom{\CL_{_{q,M}} w_k}
=(kP_{_N}+w)^p-w^p+a\gg_2|x|^{-Nq}-M\left(|k\nabla P_{_N}+\nabla w|^{2}\right)^{\frac{q}{2}}\\[1mm]
\phantom{\CL_{_{q,M}} w_k}
\geq pkP_{_N}w^{p-1}+a\gg_2|x|^{-Nq}-2M\left(k^{q}\gg_2^{q}|x|^{-Nq}+|\nabla w|^{q}\right).
\EA$$
As in Case I,  by scaling techniques, 
$$w(x)\leq \gg_3 \min\left\{a^{\frac1p}|x|^{-\frac{Nq}{p}},a|x|^{2-Nq}\right\}
$$
and
$$|\nabla w(x)|\leq \gg_4 \min\left\{a^{\frac1p}|x|^{-\frac{Nq}{p}-1},a|x|^{1-Nq}\right\}.
$$
Hence 
$$|\nabla w(x)|^q\leq \gg_5 \min\left\{a^{\frac qp}|x|^{-\frac{Nq^2}{p}-q},a^q|x|^{q(1-Nq)}\right\}.
$$
We set 
$$\gt=-\myfrac{1}{2p'-Nq}=-\myfrac{p-1}{2p-Nq(p-1)}.
$$
Then, by the definition of $\gt$,
\bel{Q3+}\BA {lll}|x|^{Nq}\CL_{_{q,M}} w_k\geq 
\gg_2\left(a-2Mk^{q}\gg_2^{q-1}\right)-2M\gg_5a^q|x|^{q(N+1-Nq)}\\
[2mm]\phantom{|x|^{Nq}\CL_{_{q,M}} w_k}
\geq 
\gg_2\left(a-2Mk^{q}\gg_2^{q-1}\right)-2M\gg_5a^{\frac{1+N-p(N-1)}{\frac{2p}{q}-N(p-1)}}\quad \text{in }\;B^+_{a^{\gt}},
\EA\ee
and  
\bel{Q4+}\BA {lll}|x|^{\frac{2Np}{p+1}}\CL_{_{\frac{2p}{p+1},M}} w_k\geq 
\gg_2\left(a-2Mk^{q}\gg_2^{q-1}\right)-2M\gg_5a^{\frac{1+N-p(N-1)}{\frac{2p}{q}-N(p-1)}}\;\text{in }\;(B^+_{a^{\gt}})^c.
\EA\ee
We obtain a very simple expression from $(\ref{Q3+})$ and $(\ref{Q4+})$, valid both in 
$B^+_{a^{\gt}}$ and $(B^+_{a^{\gt}})^c$, hence
\bel{Q5+a}\BA {lll}|x|^{\frac{2Np}{p+1}}\CL_{_{\frac{2p}{p+1},M}} w_k\geq 
\gg_2\left(a-2Mk^{q}\gg_2^{q-1}\right)-2M\gg_5a^{\frac{1+N-p(N-1)}{\frac{2p}{q}-N(p-1)}}\;\,\text{ in }\BBR^{_N}_+.
\EA
\ee
Using the scaling transformation $T_\ell$ defined in $(\ref{T1})$, the problem  of finding $u_k$ solution of $(\ref{P1})$ is equivalent to looking for a solution of 
\bel{P1-ell}\BA {lll}
-\Gd u+u^{p}-M\ell^{\frac{2p-q(p+1)}{p-1}}|\nabla u|^{q}=0\qquad&\text{ in }\BBR^{_N}_+\\
\phantom{-\Gd +u^{p}-M\ell^{\frac{2p-q(p+1)}{p-1}}|\nabla u|^{q}}
u=k\ell^{\frac{p+1}{p-1}-N}\gd_0&\text{ in }\prt\BBR^{_N}_+.
\EA\ee
If we replace $M$ by $M_\ell:=M\ell^{\frac{2p-q(p+1)}{p-1}}$ and $k$ by $k_\ell:=k\ell^{\frac{p+1}{p-1}-N}$, the inequality $(\ref{Q5+})$ turns into
\bel{Q6}\BA {lll}|x|^{\frac{2Np}{p+1}}\CL_{_{\frac{2p}{p+1},M}} w_{k,\ell}\geq 
\gg_2\left(a-2M_\ell k_\ell^{q}\gg_2^{q-1}\right)-2M_\ell\gg_5a^{\frac{1+N-p(N-1)}{\frac{2p}{q}-N(p-1)}}\;\,\text{ in }\BBR^{_N}_+,
\EA
\ee
where $w_{k,\ell}=w+k_\ell R$ instead of $(\ref{Q1e})$. Notice that $M_\ell k_\ell^{q}=M\ell^{\frac{2p}{p-1}-Nq}k^q$. We choose $\ell>0$ such that  
$M_\ell k_\ell^{q}\gg_2^{q-1}=\frac{a}{4}$, hence 
\bel{Q7}|x|^{\frac{2Np}{p+1}}\CL_{_{\frac{2p}{p+1},M}} w_{k,\ell}\geq 
\frac{a\gg_2}{2}\left(1-\gg_5\gg_2^{-q}k^{-q}a^{\frac{1+N-p(N-1)}{\frac{2p}{q}-N(p-1)}}\right)\;\,\text{ in }\BBR^{_N}_+.
\ee
It is now sufficient to choose $a>0$ such that the right-hand side of $(\ref{Q7})$ is nonnegative and thus $w_{k,\ell}$ is a supersolution. Since $\tilde v_{k,\ell}$ is a subsolution smaller that $w_{k,\ell}$, we end the proof as in Case I. 
\smallskip

\nind {\it V- Uniqueness or existence of a minimal solution.} 
If $1<q\leq \frac{2p}{p+1}$, uniqueness follows from \rth{uniq} applied with $a=N-1<\gb=\frac{2-q}{q-1}$. If $\frac{2p}{p+1}<q<\frac{N+1}{N}$ and if $u_{k,1}$ and $u_{k,2}$ are solutions, they are larger than $v_k$ and the function $u_{k,1,2}=\inf\{u_{k,1},u_{k,2}\}$ is a supersolution larger than $v_k$. Hence there exists a solution $\tilde u_k$ such that 
$$v_k\leq \tilde u_k\leq u_{k,1,2}.
$$
Let $\CE_k$ be the set of nonnegative solutions of $(\ref{NZ1})$ in $\BBR^{_N}_+$, vanishing on 
$\prt\BBR^{_N}_+\setminus\{0\}$ and such that $(\ref{P0-})$ and put
$$u_k=\inf\{\gu:\gu\in \CE_k\}.$$
Then there exists a decreasing sequence $\{\gu_j\}$ such that $\gu_j$ converges to $u_k$ on a countable dense subset of $\BBR^{_N}_+$. By standard elliptic equation regularity theory, $\gu_j$ converges to $u_k$ on any compact subset of 
$\overline\BBR^{_N}_+\setminus\{0\}$. Hence  $u_k$ is a solution of $(\ref{NZ1})$ in $\BBR^{_N}_+$, it vanishes on 
$\prt\BBR^{_N}_+\setminus\{0\}$ and $(\ref{P0})$ since $u_k\geq v_k$. Hence $u_k$ is the minimal solution.
\qeda\medskip

Next  of we consider the same problem in a bounded domain $\Gw$. 
\medskip
 
\nind{\it Proof of \rth{weak-sing2}}.  
We give first proof when $\Gw\subset\BBR^{_N}_+$. We adapt the proof of \rth{weak-sing}. The solution $v_k$ of
\bel{R1-}\BA {lll}
-\Gd v+v^{p}=0\qquad&\text{in }\Gw\\
\phantom{-\Gd +v^{p}}
v=k\gd_0&\text{on }\prt\Gw,
\EA\ee
is a subsolution for $(\ref{NZ1})$ in $\Gw$ and satisfies (\ref{Z10}).  The solution $u_k$ of $(\ref{NZ1})$ in $\BBR^{_N}_+$ vanishing on $\prt\BBR^{_N}_+\setminus\{0\}$ and satisfying $(\ref{P0-})$ is larger than $v_k$ in $\Gw$. Hence the result follows by \rprop{keylem2+}. 
\smallskip

When $\Gw$ is not included in $\BBR^{_N}_+$, 
estimates (\ref{P7}) is valid with the same type of bounds on $c$. We also consider separately the cases $q=\frac{2p}{p+1}$ and $M$ upper bounded, $q<\frac{2p}{p+1}$ and $M>0$ arbitrary and $q=\frac{2p}{p+1}$ and $M>0$ arbitrary and finally $\frac{2p}{p+1}<q<\frac{N+1}{N}$. 
As supersolution we consider the function $w_k:=kP^\Gw+w$ where $w$ satisfies
\bel{R3}\BA {lll}
-\Gd w+w^{p}=a\gg_2|x|^{-\frac{2Np}{p+1}}&\qquad&\text{in }\Gw\\
\phantom{-\Gd +w^{p}}
w=0&\qquad&\text{on }\prt\Gw,
\EA\ee
for some $a>0$. The estimates on $w$ endow the form
$$w(x)\leq \gg_3a^\frac1p|x|^{2(1-\frac{Np}{p+1})},
$$
and 
$$|\nabla w(x)|\leq \gg_4a^\frac1p|x|^{1-\frac{2Np}{p+1}},
$$
where $\gg_3$ and $\gg_4$ depend on $\Gw$. Hence $(\ref{Q5+})$ holds in $\Gw$ instead of $\BBR^{_N}_+$, and we have existence for $M<M_1$, where 
$M_1$ is defined by $(\ref{Q6e})$. Then we prove existence for any $M>0$ and $k>0$ when $q<\frac{2p}{p+1}$ then for any $M>0$ when $q=\frac{2p}{p+1}$ and finally when $\frac{2p}{p+1}<q<\frac{N+1}{N}$ as  in \rth{weak-sing}.\qeda



\subsection{Solutions with a strong singularity}

\subsubsection{The case $1< q\leq \frac{2p}{p+1}$}

If $p= \frac{N+1}{N-1}$ and $1<q<\frac{N+1}{N}$ and if $p> \frac{N+1}{N-1}$ and either $1< q<\frac{2p}{p+1}$ and $M>0$ or $q=\frac{2p}{p+1}$ and $M>m^{**}$ defined in $(\ref{Z6})$, the singularity is removable by \rth{remov}. 
Thus the ranges of exponents that we consider  are the following,
\bel{SS20}\BA{lll}
(i)\qquad &1< q\leq \frac{2p}{p+1}\,\text{ and }\; 1<p< \frac{N+1}{N-1},\qquad \qquad\qquad \qquad \qquad \qquad   \\[2mm]
(ii) \qquad &(p,q)=\left(\frac{N+1}{N-1},\frac{N+1}{N}\right).
\EA\ee
If $(\ref{SS20})$-(i) holds, $q<\frac{N+1}{N}$, and in this range the limit of the fundamental solutions $u_k$ when $k\to\infty$ is a solution with a strong singularity with an explicit blow-up rate. In the case of a bounded domain our construction requires a  geometric {\it flatness} condition
of $\prt\Gw$ near $0$. 
We consider first the case $\Gw=\BBR^{_N}_+$.

\bth{souscrit} Assume $(\ref{SS20})$-(i) holds, then for any $M\geq 0$ there exists a positive solution $u$ of $(\ref{Z1})$ in $\BBR^{_N}_+$ vanishing on $\prt\BBR^{_N}_+\setminus\{0\}$ such that 
\bel{SS1}
\lim_{x\to 0}\myfrac{u(x)}{P_{_N}(x)}=\infty.
\ee
Furthermore, \smallskip

\nind (i) If $1<q<\frac{2p}{p+1}$, 
\bel{SS2}\displaystyle
\lim_{r\to 0}r^{\frac2{p-1}}u(r,.)=\psi\qquad\text{uniformly in }\, S^{_{N-1}}_+,
\ee
where $\psi$ is the unique positive solution of $(\ref{K10})$.\smallskip

\nind (ii) If $q=\frac{2p}{p+1}$,
\bel{SS3}\displaystyle
\lim_{r\to 0}r^{\frac2{p-1}}u(r,.)=\gw \qquad\text{uniformly in }\, S^{_{N-1}}_+,
\ee
where $\gw$ is the minimal positive solution of $(\ref{Z12})$.
\es
\Proof If $k>0$, we denote by $u=u_{k,M}$ the solution of 
\bel{SS6}\BA{lll}
-\Gd u +u^p=M|\nabla u|^q\qquad&\text{in }\,\BBR^{_N}_+\\
\phantom{-\Gd  +u^p}
u=k\gd_0&\text{in }\,\prt\BBR^{_N}_+.
\EA
\ee 
The mapping $k\mapsto u_k$ is  increasing. We set $T_\ell [u]=u_\ell$, where $T_\ell$ is defined in $(\ref{T1})$.\smallskip

\nind  Since $1<q\leq \frac{2p}{p+1}$, 
$$T_\ell [u_{k,M}]=u_{k\ell^{\frac{2}{p-1}+1-N},M\ell^{\frac{2p-q(p+1)}{p-1}}}.$$
It follows from \rth{dom}  and \rth{domgrad} that the sequences $\{u_{k,M}\}$ and $\{\nabla u_{k,M}\}$ converge locally uniformly in $\BBR^{_N}_+$, when $k\to\infty$, to a function $u_{\infty,M}$ which satisfies  $(\ref{Z1})$ in $\BBR^{_N}_+$. Furthermore
\bel{SS7}T_\ell [u_{\infty,M}]=u_{\infty,M\ell^{\frac{2p-q(p+1)}{p-1}}}\qquad\text{for all }\,\ell>0.\ee
\nind In the case $q=\frac{2p}{p+1}$ the function $u_{\infty,M}$  is self-similar, hence 
$$u_{\infty,M}(r,\gs)=r^{-\ga}\tilde\gw(\gs),
$$
where $\tilde\gw$ is a nonnegative solution of $(\ref{Z12})$. Inasmuch $u_{k,M}\geq u_{k,0}=v_k$ (already defined by $(\ref{P7e})$), it follows that 
\bel{SS9}
u_{\infty,M}(r,\gs)\geq u_{\infty,0}(r,\gs)=r^{-\ga}\psi(\gs)\Longrightarrow \tilde\gw\geq\psi\quad\text{in }\;S^{_{N-1}}_+.
\ee
Since $u_{k,M}$ is dominated by any self-similar solution of $(\ref{Z1})$, it implies that $\tilde\gw$ is the minimal positive solution of $(\ref{Z12})$ that we denote by $\gw$ hereafter.  Up to a subsequence, 
$\{T_{\ell_n} [u_{\infty,M}]\}$ converges locally uniformly in $\overline{\BBR^{_N}_+}\setminus\{0\}$ to $u_{\infty,M}$. Consequently  
$$\lim_{\ell_n\to 0}\ell_n^\ga u_{\infty,M}(\ell_n,\gs)=\gw(\gs)\quad\text{uniformly in }\;S^{_{N-1}}_+.
$$
Because of uniqueness, the whole sequence converges, which implies $(\ref{SS3})$.\smallskip

\nind In the case $q<\frac{2p}{p+1}$, using the a priori estimates from \rth{dom} and \rth{domgrad}, we obtain that 
$T_{\ell_n} [u_{\infty,M}](1,\gs)=\ell_n^\ga u_{\infty,M}(\ell_n,\gs)$ converges locally uniformly in $\overline{S^{_{N-1}}_+}$ to 
$u_{\infty,0}(1,\gs)$. Since $u_{\infty,0}(1,.)\geq \psi$, it follows that 
$$\lim_{\ell_n\to 0}\ell_n^\ga u_{\infty,M}(\ell_n,\gs)=\psi(\gs)\quad\text{uniformly in }\;S^{_{N-1}}_+.
$$
Hence $(\ref{SS3})$ follows by uniqueness of the function $\psi$. \\
Uniqueness of positive solution of $(\ref{Z8})$ satisfying $(\ref{SS1})$ follows from \rth{uniq} applied with $a=\ga=\frac{2}{p-1}<\frac{2-q}{q-1}=\gb$.

\qeda\medskip

As a consequence of \rth{exist}-(ii) we have

\bth{souscrit-2} Assume $(\ref{SS20})$-(ii) holds, then for any $M> 0$ there exists a positive separable solution $u$ of $(\ref{Z1})$ in $\BBR^{_N}_+$ vanishing on $\prt\BBR^{_N}_+\setminus\{0\}$
\es

When $\BBR^{_N}_+$ is replaced by a bounded domain there holds.
\bth{souscrit-bd-1} Assume $\Gw\subset\BBR^{_N}_+$ is a bounded smooth domain such that $0\in\prt\Gw$ and $T_{\prt\Gw}(0)=\prt\BBR^N_+$,  and $(p,q)$ satisfies $(\ref{SS20})$-(i). Then for any $M\geq 0$ there exists a positive solution $u$ of $(\ref{Z1})$ in $\Gw$ vanishing on $\prt\Gw\setminus\{0\}$ such that 
\bel{SS13}
\lim_{x\to 0}\myfrac{u(x)}{P_\Gw(x)}=\infty,
\ee
where $P_\Gw$ is the Poisson kernel in $\Gw$. Furthermore \smallskip

\nind (i) If $1<q<\frac{2p}{p+1}$, then
\bel{SS13a}\displaystyle
\lim_{r\to 0}r^{\ga}u(r,.)=\psi\qquad\text{locally uniformly in }\, S^{_{N-1}}_+,
\ee
where $\psi$ is the unique positive solution of 
$$\BA {lll}
-\Gd'\psi+\ga(N-2-\ga)\psi+\psi^p=0\qquad&\text{in } S^{N-1}_+\\
\phantom{-\Gd'\psi+\ga(N-2-\ga)+\psi^p}
\psi=0&\text{in } \prt S^{N-1}_+.
\EA$$
\smallskip

\nind (ii) If $q=\frac{2p}{p+1}$, then
\bel{SS13b}\displaystyle
\psi\leq \liminf_{r\to 0}r^{\ga}u(r,.)\leq\limsup_{r\to 0}r^{\ga}u(r,.)\leq\gw \qquad\text{locally uniformly in }\, S^{_{N-1}}_+.
\ee
\es
\Proof As in the proof of \rth{souscrit}, the sequence $\{u_k\}$ of the solution of $(\ref{NZ1})$ which satisfy  $(\ref{P0-})$ is increasing. Since it is bounded from above 
by the restriction to $\Gw$ of the solutions of the same equation in $\BBR^{_N}_+$, vanishing on $\prt\BBR^{_N}_+\setminus\{0\}$ and satisfying $(\ref{SS1})$, it admits a limit $u_\infty$ which is a solution of $\ref{NZ1}$ which vanishes on $\prt\Gw\setminus\{0\}$ and satisfies $(\ref{SS13})$. In order to have an estimate of the blow-up rate, we recall that the solution $v_k$ of $(\ref{R1-})$ is a subsolution of $(\ref{Z1})$ and $u_k\geq v_k$ Furthermore $\{v_k\}$ converges to $\{v_\infty\}$ which is a positive solution of $(\ref{Z1})$ in $\Gw$, vanishing on $\prt\Gw\setminus\{0\}$ and such that 
\bel{SS14}\displaystyle
\lim_{r\to 0}r^{\ga}v_\infty(r,\gs)=\psi(\gs)\quad\text{locally uniformly in }\,S^{_{N-1}}_+.
\ee
Combined with  $(\ref{SS2})$ and $(\ref{SS3})$ it implies $(\ref{SS13a})$ and $(\ref{SS13b})$ since the solution $u_k$ in $\Gw$ is bounded from above by the solution in $\BBR^{_N}_+$.
\bel{SS15}\displaystyle
\liminf_{r\to 0}r^{\ga}u_\infty(r,\gs)\geq\psi(\gs)\quad\text{locally uniformly in }\,S^{_{N-1}}_+.
\ee
\qeda
\bth{souscrit-bd} Assume $\Gw\subset\BBR^{_N}_+$ is a bounded smooth domain such that $0\in\prt\Gw$ and $T_{\prt\Gw}(0)=\prt\BBR^N_+$,  
$p=\frac{N+1}{N-1}$ and $q=\frac{2p}{p+1}=\frac{N+1}{N}$. If 
 \bel{SS15a}\displaystyle
\dist (x,\BBR^{_N}_+)\leq c_{27}|x|^{N}\qquad\text{for all }\; x\in\prt\Gw\cap B_\gd,
\ee
for some constants $\gd,c_{27}>0$, then there exists a positive solution $u$ of $(\ref{Z1})$ in $\Gw$, vanishing on $\prt\Gw\setminus\{0\}$ such that 
\bel{SS15b}\displaystyle
\lim_{r\to 0}r^{\ga}u(r,\gs)=\gw(\gs)\quad\text{locally uniformly in }\,S^{_{N-1}}_+.
\ee
\es
\Proof The function $u_\gw(r,.)=r^{1-N}\gw$ satisfies $(\ref{Z1})$ in $\BBR^{_N}_+$ and vanishes on $\prt\BBR^{_N}_+\setminus\{0\}$. Since $\nabla\gw$ is bounded, it satisfies 
$$ u(x)\leq c_{19}\qquad\text{for all }\; x\in\prt\Gw\setminus\{0\},
$$
for some constant $c_{19}>0$. Then the result follows from \rprop{keylem2}.\qeda
 
\subsubsection{The case $\frac{2p}{p+1}<q<p$}

If 
\bel{SS21}
 1<p<\frac{N+1}{N-1}\, \text{ and }\; \frac{2p}{p+1}<q<\frac{N+1}{N},
\ee
there exists fundamental solutions $u_k$  in $\BBR^{_N}_+$ by \rth{weak-sing}, or in $\Gw$ by \rth{weak-sing2}. Since the mapping $k\mapsto u_k$ is increasing and $u_k$
is bounded from above the function $\displaystyle u_\infty=\lim_{k\to\infty}u_k$ is a solution of $(\ref{Z1})$ in $\BBR^{_N}_+$ (resp. $\Gw$) vanishing on $\BBR^{_N}_+\setminus\{0\}$ (resp. $\Gw\setminus\{0\}$) which satisfies $(\ref{SS1})$ (resp. $(\ref{SS13})$). However the blow-up rate of $u_\infty$ is not easy to obtain from 
scaling methods since the transformation $T_\ell$ transform $(\ref{Z1})$ into $(\ref{P1-ell})$ where $M$ is replaced by $M\ell^{\frac{2p-q(p+1)}{p-1}}$ which is not bounded when $\ell\to 0$. When $q>\frac{2p}{p+1}$, the natural exponent is $\gg$ defined by $(\ref{gamma})$
The transformation $S_\ell$ defined for $\ell>0$ by
 \bel{SS15c}
S_\ell[u](x)=\ell^{\gg}u(\ell x),
\ee
transforms $(\ref{Z1})$ into 
 \bel{SS15d}
-\ell^{\frac{q(p+1)-2p}{p-q}}\Gd u+|u|^{p-1}u-M|\nabla u|^q=0.
\ee
When $\ell\to 0$, the limit equation is an eikonal equation (up to change of unknown), 
 \bel{SS15e'}
|u|^{p-1}u-M|\nabla u|^q=0.
\ee
Separable solutions of $(\ref{SS15e})$ in $\BBR^{_N}_+$ are under the form $u_\eta(r,.)=r^{-\gg}\eta$ and $\eta$ satisfies 
 \bel{SS15f}
|\eta|^{p-1}\eta-M(\gg^2\eta^2+|\nabla' \eta|^2)^\frac{q}{2}=0\quad\text{in }\,S^{_{N-1}}_+.
\ee
Clearly this equation admits no $C^1$ solution but for the constant ones. As  limit of solutions with vanishing viscosity, the solutions that we obtain are  viscosity solutions outside the origin. We will 
look for solutions having a strong singularity by the method of sub and supersolutions.
Note that $(\ref{SS15e})$ admits an explicit radial singular solution, 
namely
  \bel{SS21b}
U(x)=\gw_0|x|^{-\gg}:=\gg^\gg M^{\frac{1}{p-q}}|x|^{-\gg}.
\ee

\nind {\it Proof of \rth{soupscrit}}. For $n>0$ set $U_n(r)=nr^{-\gg}$. As
$$\gg(p-1)+2=-q(p+1)+\gg+2=\myfrac{2p-q(p+1)}{p-q},
$$
we have
$$\BA {lll}n^{-1}r^{-2-\gg}\CL_{q,M}U_n=-\gg(\gg+2-N)+n^{q-1}(n^{p-q}-\gg^q M)r^{2-(p-1)\gg}.\EA$$
Since $\gg+2-N>0$ because $q>\frac{2p}{p+1}$ and $p<\frac{N+1}{N-1}$, for any $n>\gw_0$ there exists $r_n>0$ such that 
$$n^{q-1}(n^{p-q}-\gg^q M)r_n^{2-(p-1)\gg}= \gg(\gg+2-N).
$$
It implies that $U_n$ is a super solution of $(\ref{Z1})$ in $B_{r_n}\setminus\{0\}$. Furthermore
  \bel{SS21e}r_n=\left(\myfrac{n^{p-1}}{\gg(\gg+2-N)}\right)^{\frac{1}{(p-1)\gg-2}}(1+o(1))\quad\text{when }\;n\to\infty.
\ee 
For a subsolution we set 
\bel{SS5++}W_m(r,\gs)=mr^{-\gg}\phi_1(\gs),
\ee
where $m>0$. Then
\bel{SS5+++}\BA {lll}
r^{p\gg}\CL_{_{q,M}}W_m=-mr^{\frac{q(p+1)-2p}{p-q}}\left(\gg^2-(N-2)\gg+1-N\right)\phi_1\\[2mm]
\phantom{--------------}+m^q\left(m^{p-q}\phi_1^p-M\left(\gg^2\phi_1^2+|\nabla'\phi_1|^2\right)^{\frac q2}\right),
\EA\ee
and this expression is negative for $m>0$ small enough. 
Set $$P(X)=X^2-(N-2)X+1-N=(X+1)(X+1-N).$$
 Then 
$$P(\gg)=\myfrac{p\left(Nq-(N-1)p\right)}{(p-q)^2}.
$$
{\it We first give the proof when $Nq\geq (N-1)p$}. In such case $P(\gg)\geq 0$. 
Hence there exists $m_0>0$ such that for any  $0<m\leq m_0$, $W_m$ is a subsolution in $\BBR^{_N}_+$, smaller than $U_n$ and it is bounded on $\prt B^+_{r_n}\setminus\{0\}$. When $m\leq m_0$, the function $W_m$ defined in $(\ref{SS5++})$ is a subsolution of $(\ref{Z1})$ in $\BBR^{_N}_+$. Since $W_m$ is bounded on 
  $\prt B^+_{r_n}\setminus\{0\}$ there exists a nonnegative solution $u_n$ of $(\ref{Z1})$ in $B^+_{r_n}$ which vanishes on $B^+_{r_n}\setminus\{0\}$ and there holds
  \bel{SS22}
(W_m(x)-mr_n^{-\gg})_+\leq u_n(x)\leq U_n(x)\quad\text{for all }\,x\in B^+_{r_n}.
\ee
  The fact that $B^+_{r_n}$ is just a Lipschitz domain is easily bypassed by smoothing it in a neighborhood of $\prt B_{r_n}\cap \BBR^{_N}_+$. Furthermore, 
  by $(\ref{1Q1})$ and $(\ref{1Q14})$, 
    \bel{SS23}
  u_n(x)\leq c_5\max\left\{|x|^{-\ga}, M^{\frac{1}{p-q}}|x|^{-\gg}\right\}.
  \ee
  and for any $r_0>0$, there exists $c_8>0$ depending on $r_0$ such that 
      \bel{SS24}
  |\nabla u_n(x)|\leq c_8\max\left\{|x|^{-\ga-1}, M^{\frac{1}{p-q}}|x|^{-\gg-1}\right\}.
  \ee
 By standard local regularity theory, there exists a subsequence $\{u_{n_j}\}$ which converges in the $C^1(K)$-topology for any compact set 
 $K\subset \overline{\BBR^{_N}_+}\setminus\{0\}$ to a positive solution $u$ of $(\ref{Z1})$ in $\BBR^{_N}_+$ which vanishes on $\prt\BBR^{_N}_+\setminus\{0\}$ and satisfies $(\ref{SS21c})$.\smallskip
 
 \nind {\it Next we assume $Nq< (N-1)p$}. Observe that $\gg^2\phi_1^2+|\nabla'\phi_1|^2\geq \gd^2>0$, then 
 $$m^p\phi_1^p-Mm^q\left(\gg^2\phi_1^2+|\nabla'\phi_1|^2\right)^{\frac q2}\leq m^p-Mm^q\gd^q.
 $$
 Thus, from $(\ref{SS5+++})$ we obtain
       \bel{SS24a}
r^{p\gg}\CL_{_{q,M}}W_m\leq -mr^{\frac{q(p+1)-2p}{p-q}}P(\gg) +m^p-Mm^q\gd^q,
  \ee
  and $P(\gg)<0$. If we choose 
  $$m=\gd^{\gg}\left(\frac M2\right)^{\frac{1}{p-q}},
  $$
  then 
  $$m^p-Mm^q\gd^q= -\myfrac{Mm^q\gd^q}{2}.
  $$
Therefore $\CL_{_{q,M}}W_m\leq 0$ on $B_{r_*}^+$ where 
  $$r_*=\left(\myfrac{Mm^{q-1}\gd^q}{-2P(\gg)}\right)^{\frac{p-q}{q(p+1)-2p}}.
  $$
 If $a=mr_*^{-\gg}$, then $W_m\leq a$ in $\prt B_{r_*}^+$, thus $W_{m,a}=(W_m-a)+$ is nonnegative in $B_{r_*}^+$ and it is a subsolution of $(\ref{Z1})$ in $B_{r_*}^+$ which vanishes on $\prt B_{r_*}^+\setminus\{0\}$. If we extend it by 
 $0$ in $\BBR^{_N}_+$, the new function is a a subsolution of $(\ref{Z1})$ which belongs to $W^{1,\infty}_{loc}(\overline{\BBR^{_N}_+}\setminus\{0\})$. We end the proof using \rprop{keylem2+} as in the previous case.$\phantom{------}$
 \qeda\medskip


 If $\BBR^{_N}_+$ is replaced by a bounded domain we have the following result. 
 \bth{supcr-2} Let $M>0$ and $\frac{2p}{p+1}<q<p$. If $\Gw\subset \BBR^{_N}_+$ is a bounded smooth domain such that $0\in\prt\Gw$ and $T_{\prt\Gw}(0)=\prt\BBR^N_+$. If 
 \bel{SS25}
\dist (x,\prt\BBR^{_N}_+)\leq c_{20}|x|^{\frac{p}{p-q}}\quad\text{for all }\,x\in \prt\Gw\text{ near }0,
\ee
for some constant $c_{23}>0$. Then there exists a positive solution $u$ of $(\ref{Z1})$ in $\Gw$ vanishing on $\prt\Gw\setminus\{0\}$ satisfying, for some $m>0$,
\bel{SS26}
m\phi_1(\gs)\leq \liminf_{r\to 0}r^{\gg}u(r,\gs)\leq \limsup_{r\to 0}r^{\gg}u(r,\gs)\leq \gw_0,
\ee
uniformly on any compact set $K\subset S^{_{N-1}}_+$.
\es
\Proof 
We recall that $\phi_1$ is the first eigenfunction of $-\Gd'$ in $W^{1,2}_0(S^{N-1})$. Let $R>0$ and $B:=B_R(a)\subset\Gw$ be an open ball tangent to $\prt\Gw$ at $0$. Up to rescaling and since the result does not depend on the value of $M$ we can assume that $R=1$. We set $w_m(x)=m|x|^{-\gth}P_B(x)$ where 
$\gth=\gg+1-N$ and $P_B$ is the Poisson kernel in $B$ expressed by
$$P_B(x)=\myfrac{1-|x-a|^2}{\gs_N|x|^N},
$$
where $\gs_N$ is the volume of the unit sphere in $\BBR^{_N}$. Then 
\bel{SS27}\BA {lll}
m^{-1}\CL_{_{q,M}}w_m\\[2mm]\phantom{m^{-1}}
=-(\gth^2+(2-N)\gth)|x|^{-\gth-2}P_B(x)+2\gth|x|^{-\gth-1}\langle\nabla P_B(x),\frac{x}{|x|}\rangle+m^{p-1}|x|^{-p\gth}P^p_B(x)\\[2mm]
\phantom{m^{-1}=}
-Mm^{q-1}\left(\gth^2|x|^{-2(\gth-1)}P_B^2(x)+|x|^{-2\gth}|\nabla P_B(x)|^2-2\gth|x|^{-2\gth-1}\langle\nabla P_B(x),\frac{x}{|x|}\rangle\right)^{\frac q2}.
\EA\ee
 Since 
$$\BA {lll}\nabla P_B(x)
=-\myfrac{1}{\gs_N}\left(\myfrac{N(1-|x-a|^2)}{|x|^{N+1}}\myfrac{x}{|x|}+\myfrac{2(x-a)}{|x|^N}\right),
\EA$$
then
$$\BA {lll}\langle\nabla P_B(x),\frac{x}{|x|}\rangle=-\myfrac{1}{\gs_N|x|^{N+1}}\left((N-1)(1-|x-a|^2)+|x|^2\right)\\[4mm]
\phantom{\langle\nabla P_B(x),\frac{x}{|x|}\rangle}
=-\myfrac{N-1}{|x|}P_B(x)-\myfrac{1}{\gs_N|x|^{N-1}},
\EA$$
which implies in particular
$$|\nabla P_B(x)|\geq \myfrac{N-1}{|x|}P_B(x)+\myfrac{1}{\gs_N|x|^{N-1}}.
$$
If $q\geq \frac{N-1}{N}p$, equivalently $\gth\geq 0$, we have
\bel{SS28}\BA {lll}|\nabla w_m|^2=\gth^2|x|^{-2(\gth+1)}P_B^2(x)+|x|^{-2\gth}|\nabla P_B(x)|^2-2\gth|x|^{-2\gth-1}\langle\nabla P_B(x),\frac{x}{|x|}\rangle\\[4mm]
\phantom{|\nabla w_m|^2}
\geq \gth^2|x|^{-2(\gth+1)}P_B^2(x)+|x|^{-2\gth}\left(\myfrac{N-1}{|x|}P_B(x)+\myfrac{1}{\gs_N|x|^{N-1}}\right)^2\\[4mm]
\phantom{|\nabla w|^2-----------}
+2\gth|x|^{-2\gth-1}\left(\myfrac{N-1}{|x|}P_B(x)+\myfrac{1}{\gs_N|x|^{N-1}}\right)
\\[4mm]
\phantom{|\nabla w_m|^2}
\geq (\gth^2+(N-1)^2)|x|^{-2(\gth+1)}P_B^2(x).
\EA\ee
Hence
\bel{SS29}\BA {llll}
 m^{-1}\CL_{_{q,M}}w_m\leq -(\gth^2+N\gth)|x|^{-\gth-2}P_B(x)+m^{p-1}|x|^{-p\gth}P^p_B(x)\\[2mm]
 \phantom{-------------}-m^{q-1}M(\gth^2+(N-1)^2)^{\frac{q}{2}}|x|^{-q(\gth+1)}P_B^q(x)\\[2mm]
 \phantom{ m^{-1}\CL_{_{q,M}}w_m}
 \leq m^{q-1}|x|^{-p\gth}P_B^q(x)\left(m^{p-q}P_B^{p-q}(x)-M(\gth^2+(N-1)^2)|x|^{(p-q)\gth-q}\right).
\EA \ee
Now
$$P_B(x)\leq \myfrac{2}{\gs_N|x|^{N-1}}\Longrightarrow P_B^{p-q}(x)\leq \left(\myfrac{2}{\gs_N}\right)^{p-q}|x|^{(1-N)(p-q)}.
$$
Since $(1-N)(p-q)=(p-q)\gth-q$, we obtain finally that, 
$$\BA {llll}
 m^{-1}\CL_{_{q,M}}w_m\leq m^{q-1}|x|^{-q\gth}P_B^q(x)\left(m^{p-q}\left(\myfrac{2}{\gs_N}\right)^{p-q}-M(\gth^2+(N-1)^2)\right).
\EA $$
Choosing $m$ small enough we deduce that $w_m$ is a subsolution in $B$. If we extend it by $0$ in $\Gw\setminus B$, the new function 
denoted by $\tilde w$ is a nonnegative subsolution of $(\ref{Z1})$ in $\Gw$ which vanishes on $\prt\Gw\setminus\{0\}$ and satisfies $(\ref{SS26})$. The proof follows from \rprop{keylem2+}.\smallskip

\nind If $q< \frac{N-1}{N}p$, then $\gth<0$. Since  $\langle\nabla P_B(x),\frac{x}{|x|}\rangle\leq 0$,  $(\ref{SS28})$ is replaced by 
\bel{SS30}\BA {lll}|\nabla w_m|^2=\gth^2|x|^{-2(\gth+1)}P_B^2(x)+|x|^{-2\gth}|\nabla P_B(x)|^2-2\gth|x|^{-2\gth-1}\langle\nabla P_B(x),\frac{x}{|x|}\rangle\\[4mm]
\phantom{|\nabla w_m|^2}
\geq \gth^2|x|^{-2(\gth+1)}P_B^2(x)+|x|^{-2\gth}\left(\myfrac{N-1}{|x|}P_B(x)+\myfrac{1}{\gs_N|x|^{N-1}}\right)^2\\[4mm]
\phantom{|\nabla w_m|^2-----------}
+2\gth|x|^{-2\gth-1}\left(\myfrac{N-1}{|x|}P_B(x)+\myfrac{1}{\gs_N|x|^{N-1}}\right)
\\[4mm]
\phantom{|\nabla w_m|^2}
\geq (\gth^2+(N-1)^2)|x|^{-2(\gth+1)}P_B^2(x)+\left(\myfrac{1}{\gs_N^2|x|^{2(N+\gth-1)}}+\myfrac{2\gth}{\gs_N|x|^{N+2\gth}}\right)\\[4mm]
\phantom{|\nabla w_m|^2----------}
+2(N-1)\left(\myfrac{1}{\gs_N|x|^{N+2\gth}}+\myfrac{\gth}{|x|^{2\gth+2}}\right)P_B(x).
\EA\ee
Set 
\bel{SS31}
\tilde r=\min\left\{2,\left(\myfrac{1}{2\gs_N|\gth|}\right)^{\frac{1}{N-2}}\right\}.
\ee
If $x\in B\cap B_{\tilde r}(0)$, the two last terms in  $(\ref{SS30})$ are nonnegative, hence 
\bel{SS32}\BA {lll}|\nabla w|^2
\geq (\gth^2+(N-1)^2)|x|^{-2(\gth+1)}P_B^2(x)\qquad\text{for all }\,x\in B\cap B_{\tilde r}(0).
\EA\ee
Note that $B\cap B_{\tilde r}(0)=B$ if $\tilde r=2$. Choosing $m>0$ small enough we infer that $w_m$ is a subsolution of $(\ref{Z1})$ in $B\cap B_{\tilde r}(0)$. Denoting by $\hat {m}$ the maximum of $w_{ m}$ on $\prt(B\cap B_{\tilde r}(0))\setminus\{0\}$, then $(w_m-\hat {m})_+$ is a subsolution in $\Gw$. Since the restriction to $\Gw$ of the solution constructed in \rth{soupscrit} dominates $(w_m-\hat {m})_+$, the proof follows as in the first case. 
\qeda

\subsubsection{Open problems}

\nind{\it Problem 1}. Under what conditions are the positive solutions of problem $(\ref {Z12})$ unique ? If instead of separable solutions in $\BBR^{_N}_+$ vanishing on 
 $\prt\BBR^{_N}_+\setminus\{0\}$ one looks for separable radial solutions of $(\ref {Z1})$ in $\BBR^{_N}\setminus\{0\}$ (with $q=\frac{2p}{p+1}$) , then they are under the form
 \bel{SS33}\BA {lll}
 U(x)=A|x|^{-\ga}
\EA\ee
and $A$ is a positive root of the polynomial
 \bel{SS34}\BA {lll}
P(X)=X^{p-1}-M\ga^\frac{2p}{p+1}X^{\frac{p-1}{p+1}}+\ga(N-2-\ga).
\EA\ee
A complete study of the radial solutions of $(\ref{Z1})$ is provided in \cite{BVGHV3}, however it is straightforward to check that if $1<p<\frac{N}{N-2}$, there exists a unique positive root, hence a unique positive separable solution, while if $p>\frac{N}{N-2}$, there exists a unique positive root (resp. two positive roots)
 if 
  \bel{SS35}\BA {lll}
M=(p+1)\left(\myfrac{p(N-2)-N}{2 p}\right)^{\frac{p}{p+1}}:=m^*,
\EA\ee
(resp. $M>m^*$). Uniqueness of solution plays a fundamental role in the description and classification of all the positive solutions with an isolated singularity at $0$.\medskip

\nind{\it Problem 2}. It is proved in \cite{BVGHV3} that if $\max\{\frac{N}{N-1},\frac{2p}{p+1}\}<q<\min\{2,p\}$ and $M>0$, there exist infinitely many local radial solutions of of $(\ref {Z1})$ in $\BBR^{_N}\setminus\{0\}$ which satisfies 
  \bel{SS36}\BA {lll}
u(r)=\xi_Mr^{-\gb}(1+o(1))\quad\text{as }\;r\to 0
\EA\ee
where 
  \bel{SS37}\BA {lll}
\gb=\myfrac{2-q}{q-1}\quad\text{and }\;\xi_M=\myfrac{1}{\gb}\left(\myfrac{(N-1)q-N}{M(p-1)}\right)^{\frac{1}{p-1}}.
\EA\ee
These solutions present the property that there blow-up is smaller than the one of the explicit radial separable solution. It would be interesting to construct such solutions of $(\ref {Z1})$ in $\BBR^{_N}_+$ (or more likely $B^+_R$), vanishing on $\prt\BBR^{_N}\setminus\{0\}$.\medskip

\nind{\it Problem 3}. Is it possible to define a boundary trace for any positive solution of $(\ref {Z1})$ in $\BBR^{_N}_+$, noting the fact such a result holds separately for positive solutions of $(\ref {Z3})$ and $(\ref {Z4})$ ? A related problem would be to define an initial trace for any positive solution of the parabolic equation
  \bel{SS38}
\prt_tu-\Gd u+u^p-M|\nabla u|^{q}=0,
\ee
in $(0,T)\ti\BBR^{_N}$. Initial trace of semilinear parabolic equations ($M=0$ in $(\ref {SS38})$) are studied in \cite{MV0}, \cite{GkVe}. \medskip

\nind{\it Problem 4}. Are the positive solutions of $\ref{NZ1}$ satisfying $(\ref{Z10})$ or  $(\ref{SS13a})$ unique without the flatness and the starshapedness assumptions  of \rth{uniq3}. More generaly, are the weak solutions of the Dirichlet problem with measure boundary data $(\ref{NY1})$ unique ?

\end{document}